%% file: residues-fin.tex
\documentclass[12pt]{amsart}
\usepackage{amsfonts,amssymb,amscd,amsmath,color}
\usepackage{latexsym,amsbsy,bbold,stmaryrd,epsfig}
\usepackage{tocvsec2}
\usepackage{longtable}
\usepackage{mathrsfs}
\usepackage[all]{xy}

\usepackage[
hyperindex=true,pagebackref=true,bookmarks=true,
colorlinks=true,linkcolor=blue,citecolor=red]{hyperref}

\newcommand{\comment}[1]{}
\newcommand*{\medcup}{\mathbin{\scalebox{1.5}{\ensuremath{\cup}}}}%
\title [Integral formulas for DAHA inner products] 
{Integral formulas for \\ DAHA inner products}
\author[Ivan Cherednik]{Ivan Cherednik $^\dag$}
\author[Bradley Hicks]{Bradley Hicks}
\thanks{$^\dag$ \today.
\ \ \ Partially supported by NSF grant
DMS--1901796
}

\address[I. Cherednik]{Department of Mathematics, UNC 
Chapel Hill, North Carolina 27599, USA\\
chered@math.unc.edu}

\address[B. Hicks] {Department of Mathematics, UNC
Chapel Hill, North Carolina 27599, USA\\
hicksb@unc.edu}

\input{macro}

\def\mm{{\mathsf m}}
\def\aa{{\mathsf a}}

\begin{document}
\par
{\centering 
Dedicated to Eric Opdam  \\
on the occasion of his 60th birthday
\medskip
\par} 

\begin{abstract}
The main aim is to obtain
integral formulas for DAHA coinvariants
and the corresponding inner products for any values
of the DAHA parameters. In the compact case, our approach
is similar to the procedure of ``picking up  
residues" due to Arthur, Heckman, Opdam and others;
the resulting formula is a sum of integrals over double affine
residual subtori.  
A single real integral provides the required formula in the
noncompact case.
As $q$ tends to $0$, our integral formulas result in
the trace formulas for the corresponding AHA, which 
calculate the Plancherel measures for the 
spherical parts of the  regular AHA modules.
The paper contains a systematic
theory of DAHA coinvariants, including various results
on the affine symmetrizers  and induced DAHA modules.
\end{abstract}

\maketitle
{\small
Key words: Affine Hecke algebras; spherical functions; residual subtori;
Plancherel  measure; trace formulas; nonsymmetric Macdonald polynomials.

MSC (2010): 20C08, 22E35, 33C52, 33C67, 33D52, 33D67
}
{\footnotesize
\tableofcontents
}




\comment{
A method is suggested for obtaining the Plancherel measure 
for Affine Hecke Algebras as a limit of 
integral-type formulas for inner products in the polynomial 
and related modules of Double Affine Hecke Algebras.
The analytic continuation necessary here is a 
generalization of ``picking up 
residues" due to Arthur, Heckman, Opdam and others, which
can be traced back to Hermann Weyl. Generally, it is a
finite sum of integrals over double affine 
residual subtori; a complete formula is presented for $A_1$
in the spherical case. 
}

\section{\sc Introduction}
This paper is partially based on the talk by the first author
at the conference
``From E6 to double affine E60" in the
honor of Eric Opdam's 60th birthday. The main
aim is to obtain
integral formulas for DAHA {\sf coinvariants}
and the corresponding inner products for any values
of the DAHA parameter $t=q^k$, where $0<q<1$. 
As $q\to 0$ and upon the restriction to symmetric functions,
our integral formula results in
the {\sf trace formula} for the corresponding
Affine Hecke Algebra, AHA for short. This formula 
calculates the Plancherel measure for the decomposition
of the spherical part of the regular
representation of AHA in terms of irreducible unitary modules.
The standard AHA {\sf trace} is the limit $q\to 0$ of 
the DAHA {\sf coinvariant} for the
anti-involution $\Diamond$ associated with the basic
inner product in the polynomial representation.
\vfil

There are two 
directions of this paper: algebraic theory of DAHA coinvariants, 
including the affine symmetrizers and norm-formulas, and
integral formulas for  DAHA coinvariants and
inner products in the compact and noncompact settings. The
corresponding identities are {\sf DAHA trace formulas}.
The integral formulas in the compact case 
are obtained in a way similar to ``picking up  
residues" due to Arthur, Heckman, Opdam and others
(can be traced back to Hermann Weyl); they are sums over
{\sf double affine residual subtori}. The DAHA-invariance
of our formulas is an important new tool, but the combinatorial
aspects are involved so far in the $q,t$-theory. However,
a single real integral 
provides the required meromorphic continuation  
in the noncompact case.
\vfil

A challenge here is to upgrade this approach to 
{\sf global fields}: with the $c$-functions expressed
via  the completed Dedekind zeta-function: 
Kazhdan-Okounkov \cite{KO} and
 De Martino-Heiermann-Opdam \cite{MHO}. The trace formula
becomes Langlands' formula for the inner product
of two {\sf pseudo-Eisenstein series}, $(\th_\phi,\th_\psi)$ in the
notation from \cite{MHO}; see there for the definitions and
justifications. One of the key points of these two papers is
that Dedekind's zeta can be replaced by other functions
satisfying the functional equation (to ensure the cancelation
of the ``unwanted" residues). We expect that the {\sf adelic} product
of DAHA trace formulas can serve global fields, where adding $q$
provides new and interesting deformations of Langlands' formulas.

\vskip 0.2cm
\vfil
Our starting point is that
the integral formulas for the {\sf level-zero and level-one} 
coinvariants  are relatively straightforward 
for $|t|<1$ ($\Re k>0$). They  generalize 
the Macdonald formula 
in the AHA theory for $|t|>1$ outlined in 
Section \ref{SEC:AHA}; the DAHA $t$-parameter 
used throughout this paper 
corresponds to $1/t$ in the standard AHA setting.

\vfil
These integrals are essentially  
$\int f(x) \mu(q^x; q,t) dx$ for suitable spaces
of functions integrated over $\imath \R^n$ in the {\sf compact case}
and $\R^n$ in the {\sf noncompact case}, 
 where $\mu$ is the measure-function
in DAHA theory defined in (\ref{mutildemu}). The 
problem is to extend them to $|t|>1$.

All basic DAHA facts and references we need can be
found in this paper. We frequently  adjust them,
generalize and develop. 
See \cite{C101,C5,CMa,ChD, ChA} for the main 
features of DAHA inner products and coinvariants.
One of the changes vs. \cite{C101} is that we use the
anti-involution $\Diamond$ that does not involve the
conjugation of $q,t$. 

Beyond the DAHA theory, only ``$q$-calculus" and standard
theory of residues is really necessary to
obtain our  integral formulas, though they appeared involved.
This is similar to \cite{HO,O11}. 
\vskip 0.2cm
\vfil

{\bf Meromorphic continuations.}
In this paper, we mostly consider the spaces of
Laurent polynomials or Laurent 
series $f(q^x)$, which are in 
terms of $q^{x_{\al_i}}$ for $x_{\al_i}=(x,\al_i)$ 
for simple roots $\al_i$. The integration is mostly 
over the (imaginary) periods of $q^{x_{\al_i}}$;
however, the full imaginary
integration, real and 
Jackson integrations play a significant role too. 
\vfil
 
Generally,  the problem is to find the meromorphic
continuation of the imaginary integrals 
to $\Re k\le 0$ ($|t|>1$). Interestingly, a single 
integral provides the required  meromorphic function
for all $k$ (with sufficiently small $\Im k$) in the noncompact 
case:
for the integrations in the real directions.
 This is Theorem \ref{thm:noncom};
 ``picking up residues" is not needed there. We note that 
the reciprocals of theta-functions and their
expansions occur naturally in this case when the Gaussians are added
to the space of functions. See e.g. \cite{Car}.
The corresponding $q,t$-Gauss integrals,
noncompact variants of difference Macdonald-Mehta formulas from
\cite{C5}, 
involve Appell functions and similar ones.
\vfil

In the {\sf compact case}, the integrals
are analytic in terms of $k$ with $\Re k$ from 
a disconnected union of segments between the consecutive
singularities
of $\mu$.  The problem is
to use these integrals to obtain the meromorphic continuation
of the initial integral, the one  serving $\Re k>0$, to all 
$0\ge \Re k>-\infty$. 

The resulting formula is a linear combination of 
finitely many integrals over
{\sf double affine residual subtori}, where their number depends
on $\Re k$. See Theorem \ref{thm:iterint}. The contribution of
{\sf double affine residual points} is very  interesting;
the corresponding residues 
generalize {\sf formal degrees} of the AHA discrete series.
Concerning the latter, let us mention here (at least) 
Kazhdan, Lusztig, Reeder, Shoji, Opdam, Ciubotaru, 
S.Kato; also, see some references 
below.

Similar to the AHA theory,
the leading term of the resulting integral
formula is $\int_{\imath \R^n} f(x) \mu(q^x) dx$, where 
 $\Re k$ is arbitrary negative.
 This functional is AHA-invariant but not 
DAHA-invariant; so ``corrections" are needed, which are integrals
over residual subtori.  It is expected that
the DAHA-invariance of our formula
is sufficient to fix uniquely 
the corresponding ``measures" of residual
subtori and those in the AHA limit $q\to 0$. So the
action of DAHA is a major
``hidden symmetry" of the AHA Plancherel formula, which
is of conceptual importance.

The meromorphic continuation is basically by shifting the
contours of integration in the real directions followed by 
``picking up the residues".  We need the analyticity of $f(x)$ 
to ensure that the contours can be moved and the integrability
in the imaginary directions.  When $k\to 0_-$, the link to
the procedure from \cite{HO} is
discussed in $(a)$ from ``Concluding remarks"
after Theorem \ref{thm:iterint}. 
We note that  the pole decomposition
is the key to our approach in some contrast to AHA.

\vskip 0.2cm 
\vfil
{\bf DAHA aspects.}
Our integral formulas are not directly related to the 
reducibility of the polynomial $\HH$-module $\mathscr{X}$, 
where $\HH$
denotes DAHA. The reducibility is for {\sf singular} $t$,
some special $t$ satisfying $|t|>1$;
our formulas are for any $|t|>1$. 
However, there is an important connection.
When the coefficients (the residues) in our formula
have singularities in terms of $k$, our integral formulas
for the inner product 
result in a certain Jantzen-type
filtration of $\mathscr{X}$ in terms of 
$\HH$-modules. Namely, 
the largest submodule  is the radical of
the leading term of the inner product, the 2{\small nd} is
the radical of the restriction of our
integral formula to the 1{\small st}
and so on. 
For $A_n$, this filtration is essentially sufficient to decompose
the polynomial representation (at least for small $n$).
See e.g. \cite{En,C6} about the so-called
Kasatani decomposition.
 
Moreover, the subquotients here 
are naturally supplied with inner
products, given by some integrals,
that can be {\sf unitary} even if $\mathscr{X}$
is not unitary; see \cite{ChA},
Corollary 6.3 for an example. This is always very interesting.
Generally, the problem
of {\sf unitary dual} is one of the keys in harmonic analysis;
see e.g. \cite{ES} for the case of rational DAHA.
\vskip 0.2cm
\vfil

In contrast to the trace of AHA, it is not immediate to see
that the DAHA coinvariants are meromorphic functions in terms
of $t$. This fact can be proven via $(a)$ the theory of nonsymmetric
Macdonald polynomials, $(b)$ the theory of basic anti-involutions
 or $(c)$ the theory of affine symmetrizers. 

The existence of the affine symmetrizer
$\hat{\mathscr{P}}_+(f)$ and its proportionality
to $\hat{\mathscr{I}}_+(f)$ from Theorem \ref{LEVZERO}
seem the most fundamental here. The origin is 
in the $p$-adic
theory of spherical functions. Basically, we
 generalize
the fact that Matsumoto spherical functions can be identified
with nonsymmetric Hall polynomials in the AHA theory. 
\vskip 0.2cm

We extend in this paper the theory of {\sf basic anti-involutions
and coinvariants}
to $Y$-induced DAHA modules $\i_\xi$, where $\xi\in \C^n$
is considered as a character of the $Y$-subalgebra of $\HH$.
 For instance, Theorem 
\ref{thm:norm} gives the norm-formulas for such representations
and simultaneously proves the uniqueness of the corresponding
coinvariant  up to proportionality
for sufficiently general $q,t,\xi$. 

Modules $\i_\xi$  are important in this paper because of
several reasons. First, they
are related to residual points $\xi$;
the irreducible
quotients of $\i_\xi$ for ``non-Steinberg" $\xi$ are
interesting analogs of $\mathscr{X}$.
Second, generic $\i_\xi$ can be naturally identified with
the full regular representation of AHA, the main subject of
the AHA harmonic analysis. One can define the integration
and obtain integral formulas for $\i_\xi$, but this is beyond 
the present paper. This is related to Jackson integrals 
$J_\xi$ and {\sf global hypergeometric functions} from \cite{C5};
see also \cite{Sto2,SSV}.
\vskip 0.2cm
\vfil

{\bf Some perspectives.}
See also ``Concluding remarks" after Theorem \ref{thm:iterint}.
The decomposition of the regular AHA representation
in terms of irreducible modules 
involves deep geometric methods 
(Kazhdan- Lusztig and others) and a lot of functional
analysis (Opdam and others).  Our approach potentially allows
finding the formal degrees of AHA discrete series
via DAHA without any geometry. Paper \cite{O11}
does this within the AHA theory. The DAHA approach  is
expected to be analytically more transparent 
and with additional rich symmetries, which are not present in AHA.
The $q,t$-generalization of the
discrete series remains to be discovered. Actually,
the whole $\mathscr{X}$ behaves as such for sufficiently large $\Re k<0$; 
the
affine symmetrizer $\hat{\mathscr{P}}$ acts there, which is
an important feature of AHA representations of discrete series.
\vskip 0.2cm

As we already discussed, $\mathscr{X}$
is the spherical quotient
of the regular AHA representation supplied with the structure of $\HH$-module,
 $\i_\xi$ are those for the whole regular representation.
The classical AHA trace becomes
the basic $\HH$-coinvariant. 
The presentation of the 
trace as some integral over {\sf unitary dual} is reduced 
to some combinatorial calculations for $\HH$. They
are far from simple but no DAHA {\sf unitary dual} is needed.  

In the case of $A_n$, we 
calculate explicitly in Theorem  \ref{thm:Anres}
the required meromorphic continuation to $|t|>1$ as 
the {\sf pole decomposition} of the ``slightly shifted" initial
integral. This can be generalized
to any root systems and any orderings of iterated integrations, 
but the combinatorics of the resulting formulas requires 
further analysis. Moreover, non-Steinberg-type residual
points  occur beyond $A_n$.
\vskip 0.2cm

The pole decompositions we obtain converge at any $\Re k<0$, but 
only for  relatively small spaces
of $f(x)$ depending on $\Re k$. Such $f$ are 
basically Laurent polynomials of degrees bounded by
$const[-\Re k]$ or the corresponding 
Paley-Wiener functions. Practically arbitrary analytic functions 
$f(x)$ can be considered when we switch to {\sf finite} sums 
of integrals over certain double affine residual subtori.  

This passage is a combinatorial problem,
but not a simple one. Essentially, we need to 
combine the poles into families corresponding to 
proper residual subtori.
We provide the final finite 
integral formulas 
only for $A_2$ in Section  \ref{sec:a2}; see 
\cite{ChA} for the case of $A_1$. For arbitrary root
systems, 
the calculations are involved even in the
simplest interval   $0>\Re k >-\frac{1}{h}$
for the Coxeter number $h$, where the combinatorics of residual
subtori is similar to that from \cite{HO, O11}. 
\vskip 0.2cm

A natural challenge here
is the case of nonequal parameters $t_{\sht}$ and $t_{\lng}$ 
for the root systems $BCFG$, i.e. for
generic $k_{\sht}$ and $k_{\lng}$.
All main results in this paper are for any sufficiently
general $k_\nu$. For
instance, the  pole decomposition 
is obtained for  any \ $\Re k_{\sht}<0$ and $\Re k_{\lng}<0$.
However, the
explicit combinatorial description of residual points is 
provided  only when 
$k_{\lng}=\kappa k_{\sht}$, where $\kappa=1$ or
$\kappa={(\al_{\lng},\al_{\lng})}/{(\al_{\sht},\al_{\sht})}$.
See Theorems \ref{thm:residues} and \ref{thm:resid} in terms of
the closed root subsystems of maximal rank in affine
root systems.

\vskip 0.2cm
\vfil
The harmonic analysis and unitary dual for DAHA are
open projects by now. However, there are quite 
a few {\sf special theories}, where this paper can be
used as such. They are $(a)$\, the 
AHA limit as $q\to 0$ (the starting point for us),  
$(b)$\, the Kac-Moody limit 
as $t\to \infty$
($0<q<1, k\to -\infty$), $(c)$\, 
the quantum groups as $t=q$,
$(d)$ level-one Demazure characters as $t=0$,  
and $(e)$\,
the Heckman-Opdam limit \cite{HO1}:  $q\to 1, t=q^k$. 

Case $(d)$ and the limit $t\to\infty$ correspond to
the theory of nil-DAHA; see \cite{ChO2, ChK}.
In the case of  $(e)$, 
the variables $X_b=q^{x_b}$ for $b\in P$ are considered
unchanged in the limit (they occur as torus coordinates). For $(c)$,
there are actually two quantum group limits in the twisted case: 
when $t_{\sht}=q$ and
$t_{\lng}=t_{\sht}^\kappa$ for $\kappa$ as above. 
\vskip 0.2cm
\vfil

The simplest ``special theory" is actually for $t=1$; then  
DAHA becomes the {\sf Weyl algebra}. It generalizes the main
feature of the latter, the 
projective action of
 $PSL_2(\Z)$. It is the key feature of DAHA theory,
which collapses
in the limits above unless in the following two cases.
\vfil

First, this action exists for the {\sf reduced category} 
in case $(c)$ when $q$ is a root of unity
and, equivalently (due to Kazhdan-Lusztig
and Finkelberg), for the category of
integrable Kac-Moody modules in case $(b)$. The Grothendieck ring 
of the reduced category becomes then the 
{\sf perfect representation} of DAHA at $t=q$.

The second case is the action
of $PSL_2(\C)$ in
the {\sf rational} Heckman-Opdam theory (with
the Calogero operators instead of the Sutherland ones in physics
literature). This is the limit $q\to 1, t=q^k$, where 
$x_b$ above are taken as the variables. The Fourier transform,
which is the action of {\tiny 
$\left(\begin{array}{cc}0& -1 \\
                        1 & 0 \\
\end{array} \right)$} becomes the non-symmetric Hankel transform
(due to Dunkl for any root systems and due to Hermite for $A_1$). 
When $k=0$, we arrive at the Heisenberg algebra. 
\vfil

We note that the usage of Lie groups only is generally
insufficient to incorporate
the Fourier transform; one needs the Heisenberg-Weyl algebras
and DAHA, their (flat) deformations.  Similar to the classical 
polynomial representations
for Heisenberg-Weyl algebras, 
DAHA provides {\sf nonsymmetric theories}, which were 
new even for $A_1$. The nonsymmetric Macdonald polynomials
generalize the characters and spherical functions in the Lie theory,
which are symmetric (unless for Demazure characters). Our paper
is ``nonsymmetric".
\vskip 0.2cm

{\bf Acknowledgements.} We thank very much Eric Opdam
for his support and help, which is and was well beyond this paper. 

\section{\sc Affine root systems and AHA}\label{SEC:AHA}
Let $R\subset \R^n$  be a reduced irreducible (indecomposable)
 root system,\,
$Q=\oplus_{i=1}^n \Z\al_i, P=\oplus^n_{i=1}\Z \om_i$, 
where $\al_i$ are simple roots and
$\{\om_i\}$ are fundamental weights:
$ (\om_i,\al_j^\vee)=\de_{ij}$ for the
coroots $\al^\vee=2\al/(\al,\al).$
Replacing $\Z$ by $\Z_{\pm}=\{m\in\Z, \pm m\ge 0\}$, we obtain
$P_\pm$ and $Q_\pm$. See  e.g., \cite{Bo} or \cite{C101}. 
The normalization will be {\sf twisted} throughout this paper:
$(\al_{\sht},\al_{\sht})=2$ for short roots. 
Accordingly, $\vth=\vth(R_+)$ will denote the 
maximal {\sf short}
root in $R_+$, the set of positive roots. When necessary, we use
the notation $\th=\th(R_+)$ for the maximal (long) root. 
One has $\vth(R_+)=\th(R^\vee_+)$ due to our normalization 
of $(\cdot,\cdot)$, which means that  
$\vth$ becomes the maximal root in $R^\vee=\{\al^\vee, \al\in R\}$.

\vskip 0.2cm

Setting 
$\nu_\al\equal (\al,\al)/2$,
the vectors $\ \tal=[\al,\nu_\al j] \in
\R^n\times \R \subset \R^{n+1}$
for $\al \in R, j \in \Z $ form the
{\sf twisted affine root system\,}
$\tR \supset R$, where $\al\in R$ are considered as $ [\al,0]$.
We will sometimes use the notation $\nu_{\sht}$ and $\nu_{\lng}$
for short and long roots.

The inner product $(\tal,\tbe)$ is that from $\R^n$. i.e. the affine
components are ignored. However, somewhat abusing
the notation, we set $(\tal,z)=(\al,z)+\nu_\al j$, when
the pairing is between $\R^n\ni z$ and $\tal$ is considered, which
will be obvious from the context. In \cite{C101}, we
used the notation $(\tal,z+d)$ for this pairing.
\vfil

We add $\al_0 \equal [-\vth,1]$ to the simple
roots. 
The corresponding set
$\tR_+$ of positive roots is 
$R_+\cup \{[\al,\nu_\al j],\ \al\in R, \ j > 0\}$. The 
corresponding {\sf affine (extended) Dynkin diagram} will be
the usual extended one for $R^\vee$ where all arrows are
reversed.

Note that
$P\subset P^\vee$ and $Q\subset Q^\vee$ for $P^\vee,Q^\vee$
defined for $R^\vee$. The {\sf minuscule weights} are $\om_r$
such that $(\om_r,\al^\vee)\le 1$ for any $\al\in R_+$. Equivalently,
$\nu_r n_r=1$ where $\vth=\sum_{i=1}^n n_i\al_i$. The usage of
the name ``twisted" is not as in Kac-Moody theory, but there
is a direct connection for the systems $B,C,F,G$.   

\vskip 0.2cm
The twisted setup  is convenient for us
because it is ``self-dual" with respect to the DAHA
Fourier transform. Also, the ``level-one theory" for the $C$-type 
in the untwisted setting is actually ``level-two", much more
difficult than ``level-one" is supposed to be.
There are other advantages, but the untwisted root systems are
generally equally important and quite standard in Kac-Moody
theory. 
\vskip 0.2cm

The set of the indices of the images of $\al_0$ under the
action of 
automorphisms of the affine Dynkin diagram will be denoted by 
$O$. 
Let $O'\equal\{r\in O, r\neq 0\}$. The minuscule $\om_r$ are
those for $r\in O'$. We set $\om_0=0$
for the sake of uniformity. All fundamental weights are minuscule
for $A_n$. There are no minuscule weights and $O'=\emptyset$
for $E_{7,8},F_4,G_2$. 
\smallskip

{\bf Affine Weyl group.}
Given $\tal=[\al,\nu_\al j]\in \tR,  \ b \in P$, let
\begin{align}
&s_{\tal}(\tz)\ =\  \tz-(z,\al^\vee)\tal,\
\ b'(\tz)\ =\ [z,\ze-(z,b)]
\label{ondon}
\end{align}
for $\tz=[z,\ze] \in \R^{n+1}$.

The 
{\sf affine Weyl group\,} $\tW=\lan s_{\tal}, \tal\in \tR_+\ran)$ 
is the semidirect product $W\lsmash Q$ of
its subgroups $W=$ $\lan s_\al,
\al \in R_+\ran$ and $Q$, where $\al\in R\subset Q$ are identified 
with the following elements  from $\tW$: 
\begin{align*}
& R\ni \al\mapsto s_{\al}s_{[\al,\,\nu_{\al}]}=\
s_{[-\al,\,\nu_\al]}s_{\al}\in \tW. 
\end{align*}

The {\sf extended Weyl group\,} $ \hW$ is $W\lsmash P$, which can
be defined via
its action in $\R^{n+1}$ extending that of $\tW$ in $\tR$: 
\begin{align}
&(wb)([z,\ze])\ =\ [w(z),\ze-(z,b)] \for w\in W, b\in P, z\in \C^n.
\label{ondthr}
\end{align}
Notice the minus-sign of $-(z,b)$. 

We need the action of
$\hW\ni \hw$ on the functions $X_{[a,\la]}\equal q^{x_a+\la}$ 
for $x_a=(x,a)$, which is
defined as the {\sf action on the indices}:
$\hw(X_{[a,\la]})=X_{\hw([a,\la])}$. Generally, 
$\hw(f(x))\equal
f(\hw^{-1}(x))$ for any function of $x$. This action is
dual to the following {\sf affine action} on vectors $z\in \C^n$: 
$\hw(\!(z)\!)\equal w(z)+b$ for $\hw=bw$. The corresponding extension
of the pairing $(\cdot, \cdot)$ 
is $(z,[a,\la])=(z,a)+\la$.  Namely, one has:
\begin{align*}
wb(x_a)&= (x,wb(a))=
(x,[w(a),-(b,a)])=(x,w(a))-(b,a) 
\\
&=(w^{-1}(x)-b,a)=(b^{-1}w^{-1}(\!(x)\!),a)=wb(x_a),
\end{align*}
where the former $wb(x_a)$ is the action on indices,
the latter $wb(x_a)$ is the action on functions of $x$.
We will use the notation$(\!(\cdot)\!)$ only when some confusion
is likely; almost always, $\hw(\cdot)$ will denote
either the action on $[a,\la]\in \R^{(n+1)}$ or
on $z\in \C^n$ depending on the context.
Throughout the paper:  $X_a=q^{x_a}, X_{[a,\la]}=q^\la X_a$, 
and we set 
$X_a(q^b)=q^{(a,b)}$ for $X_a$ and other functions
of $\{X_a\}$.

{\sf The 
Gaussian} $q^{x^2/2}$ is defined
for $x^2=\sum x_{\al_i}x_{\om_i}$;
it is $W$-invariant, and $b(q^{x^2/2})=q^{b^2/2}X_b^{-1}q^{x^2/2}.$
It is sometimes used as a symbol, when only the action of $\hW$ on it 
is of importance. However,  $q^{\pm x^2/2}$ will be
considered functions for real and imaginary integrations. 
\vskip 0.2cm

The group $\hW$ is isomorphic to $\tW\lsmash \Pi$ for 
$\Pi\equal P/Q$. 
The latter group consists of $\pi_0=$id\, and the images $\pi_r$
of minuscule $\om_r$ in $P/Q$; also, see (\ref{ururstar}).
We note that 
$\pi_r^{-1}$ is $\pi_{r^\varsigma}$,
 where $\varsigma$ is 
the standard involution 
(sometimes trivial) of the {\sf nonaffine\,} Dynkin diagram 
induced by $\al_i\mapsto -w_0(\al_i)$, where $w_0$ is the 
longest element in $W$.
Generally
$\varsigma(b)=-w_0(b)=b^\varsigma$;
we set $X_b^\varsigma=X_{b^\varsigma}$. 

The group $\Pi$
is naturally identified with the subgroup of $\hW$ of the
elements of zero length; the {\sf length\, } is defined as 
follows:
\begin{align}\label{La-set}
&l(\hw)=|\La(\hw)| \for \La(\hw)\equal\tR_+\cap \hw^{-1}(-\tR_+).
\end{align}
I.e. $\La(\hw)$ is the set of positive affine roots that become
negative upon the application of $\hw$. 
Similarly, let $l_\nu$ be the number of $\tal$ in $\La(\hw)$ 
with $\nu_\al=\nu$.
Setting $\hw = \pi_r\tw \in \hW$ for $\pi_r\in \Pi,\, \tw\in \tW,$
\,$l(\hw)$ coincides with the length of any reduced decomposition
of $\tw$ in terms of the simple reflections
$s_i,\, 0\le i\le n$. Respectively, let $l_\nu$ count the number
of $s_i$ for short and long $\al_i (i\ge 0)$. 

For the sake of completeness, we mention 
that the equivalence of these two definitions is based on the key
property of $\La$-sets:
$$ \La(\hw\hu)=\hu^{-1}\bigl(\La(\hw)\bigr)\cup \La(\hu),\ \, 
\La(\hw^{-1})=
-\hw \bigl(\La(\hw)\bigr).
$$
The union here is disjoint if $l(\hw\hu)=l(\hw)+l(\hu)$; generally,
the cancelation of any pairs $\{\tal, -\tal\}$ must be performed
if they occur in the union. See e.g. \cite{C6}. Also, 
$l(b)=2(\rho^\vee, b)$ for
$b\in P_+$. Here and below $\rho=\frac{1}{2}
\sum_{\al>0}\al=\rho_{\sht}+\rho_{\lng}$, \ 
 $\rho^\vee=\frac{1}{2}\sum_{\al>0} \al^\vee=
\sum_\nu \frac{\rho_\nu}{\nu}=\rho_{\sht}+
\frac{\rho_{\lng}}{\nu_{\lng}}.$ 
For $b=\om_i$, $l(\om_i)$ gives the number of $\al\in R_+$ that
contain $\al_i$. 

One has $\om_r=\pi_r u_r$ for $r\in O'$, where $u_r$ is the 
element $u\in W$ of {\sf minimal\,} length such that 
$u(\om_r)\in P_-$,
equivalently, $w=w_0u$ is of {\sf maximal\,} length such that  
$w(\om_r)\in P_+$. The elements $u_r$ are very explicit.
Let $w^r_0$ be the longest element
in the subgroup $W_0^{r}\subset W$ of the elements
preserving $\om_r$; this subgroup is generated by $s_i$ for 
$0<i\neq r$.
One has:
\begin{align}\label{ururstar}
u_r = w_0w^r_0 \hbox{\,\, and\,\, } (u_r)^{-1}=w^r_0 w_0=
u_{r^\varsigma} \for r\in O.
\end{align}
\smallskip
For instance, $\om_1=\pi_1 s_3s_2s_1$, $\om_2=\pi_2 s_2 s_1s_3s_2$,
$\om_3=\pi_3 s_1s_2s_3$ for $A_{3}$. For $B_n$: $\al_n$ is a unique
short simple root, $\om_n=\al_{n-1}+2\al_n$ is a unique
minuscule weight and $\vth=\al_1+\cdots+\al_n$. Also, 
$\om_n=\pi_n u_n$, where $u_n$ sends 
$\al_i\mapsto -\al_{n-i}$ for 
$1\le i\le n-1$ and 
$\al_n\mapsto -\vth$.

\smallskip

The extended 
{\sf Affine Hecke Algebra} for $\tR$, AHA for short,
 is defined as the span
$\h\equal \lan\Pi,\, T_i(0\le i\le n)\ran$ for the
generators subject to the standard 
homogeneous Coxeter relations for $T_i$ and the quadratic 
relations $(T_i-t_i^{\frac{1}{2}})(T_i+t_i^{-\frac{1}{2}})=0$ for
$1\le i\le n$, where $t_i$ depends only on the length of 
$\al_i$, i.e. on $\nu_i=\nu_{\al_i}$. The ring of coefficients 
will be $\Z[t_i^{\pm 1/2}]$ or $\C$ if $t_i$ are 
 considered in $\C^*$. Concerning $\Pi$, the following
relations are imposed: $ \pi_rT_i\pi_r^{-1}\ =\ T_j$ if 
$\pi_r(\al_i)=\al_j$ for $r\in O', 0\le i\le n$.

In the standard $p$-adic setting, $t=p^\ell$, where $p^\ell$ is the
cardinality of the corresponding residue 
field $\mathbb F$; different $t_{\sht}, t_{\lng}$ are in the so-called  
case of unequal parameters. The DAHA $t$ is actually
$p^{-\ell}$ (below).  
\vskip 0.3cm

We set $T_{\hw}\!=\!\pi T_{i_l}\!\cdots\! T_{i_1}$ for 
any reduced decomposition
$\hw\!=\!\pi s_{i_l}\!\cdots\! s_{i_1}\!\in\! \hat{W},$
i.e. where $l\!=\!l(\hw)$. 
Considering $P$ as a subgroup in $\hat{W}$ we obtain
that $Y_b=T_b$ for $b\in P_+$ (for dominant weights)
are pairwise commutative. Then we extend it to any $b\in P$
using $Y_{b-c}=Y_bY_c^{-1}$ for $b,c\in P_+$. This is due to
Bernstein-Zelevinsky and Lusztig. The defining relations of $\h$
in terms of $Y_b$ are:
$T_iY^{-1}_b=Y_b^{-1} Y_{\al_i} T_i^{-1} \text{ for }
(b,\al^\vee_i)=1, 
0 \le i\le  n$, and  $T_iY_b= Y_b T_i \text{ for } (b,\al^\vee_i)=0,$
where $0 \le i\le  n.$
\vskip 0.2cm
\vfil

The canonical anti-involution, trace and unitary scalar product are: 
$T_{\hw}^\ast\equal T_{\hw^{-1}},\, \lan T_{\hw}\ran_{reg} 
=\de_{id,\hw},$ $\lan f,g\ran_{reg} \equal \lan f^\ast g\ran_{reg}= 
\sum_{\hw\in \hat{W}}\,\bar{c}_{\hw}d_{\hw},\, $ 
 where $f\!=\!\sum c_{\hw} T_{\hw},\ g\!=\! \sum d_{\hw} 
T_{\hw}\, \in\, L^2(\h)=\{f\, \mid\, \sum 
\bar{c}_{\hw} c_{\hw}<\infty \}$.  We assume that $t_i$ are 
real and add the complex conjugation to the definition
of $\ast$, which results in  $\bar{c}_{\hw}$. The complex
conjugation, which is necessary for unitarity, 
will be omitted in the DAHA theory below.
\vskip 0.2cm

In the spherical case, we consider $\h \p_+$ for the symmetrizer 
$\p_+=\frac{\sum_{w\in W}t^{-l(w)/2}T_w^{-1}} 
{\sum_{w\in W}t^{-l(w)}}$. By definition, 
$t^{l(\hw)/2}=\prod_{\nu}  t_\nu^{l_\nu(\hw)/2}=
  t_{\sht}^{l_{\sht}(\hw)/2} t_{\lng}^{l_{\lng}(\hw)/2}$.

This space has a natural left action of $\h$. We have 
$\h \p_+=\C[Y_b]\p_+$, for the algebra of Laurent polynomials 
$\C[Y_b, b\in P]$,  which identification is the key in the theory of 
nonsymmetric {\sf 
Matsumoto spherical functions}; see \cite{Mat,O12}, 
\cite{C101} (Section 2.11.2) and \cite{Ion,CMa}.
For instance, the formulas for 
$\lan P(Y) 
\ran$, where $P(Y)\in \C[Y_b, b\in P]$, are sufficient to
recover the trace for any $T_{\hw}$.  

\vskip 0.2cm
\vfil
According to Dixmier, 
 $\lan f,g\ran_{reg}=
\int_{\pi\in\h^\vee} \hbox{Tr}(\pi(f^\star \bar{g}))d\eta(\pi)$
for some non-negative measure $d\eta$ in the space $\h^\vee$
of irreducible unitary
$h$-modules $\pi$, the unitary dual of $\h$.
We omit here some analytic details concerning
the classes of functions. In the spherical
case (referred to as ``sph" below), one takes
$f,g\in \p_+\h \p_+$. It terms of $Y_b$, we consider the
symmetric ($W$-invariant)  Laurent polynomials, which is
based on the so called {\sf Bernstein Lemma}.  
The measure reduces correspondingly. 

\vskip 0.2cm
Macdonald found that
 $\eta_{sph}(\pi)$ as $t>1$ (in the case of one $t$) 
is proportional to $\frac{ds}{c(s,t)c(s^{-1},t)}$ in terms of the 
corresponding $c$-function, where $s \in \exp(\imath \R^n)$. 
 Its meromorphic extension to $0<t<1$ can be obtained by
``picking up residues"
due to Arthur, Heckman, Opdam and others
\cite{CKK,HO,O11,OS}. The final (spherical) 
formula reads:
 
\noindent
$$\int \{\cdot\}\, d\eta_{sph}^{an}(\pi)= 
\sum C_{s_\circ S}\cdot\int_{s_\circ S}\{\cdot\} \,d\eta_{s_\circ S},
$$

\noindent 
summed over 
{\sf affine residual subtori} $s_\circ S$, where $s_\circ S=
\exp(a_\circ +\imath \mathfrak{a})$ for some $a_\circ\in \R^n\supset
 \mathfrak{a}.$
 See \cite{HO}.
{\sf Residual points} are when dim\,$\mathfrak{a}=0$; they
correspond to
square integrable irreducible modules:
their characters $\chi_\pi$ extend to $L^2(\h)$.
This formula involves deep geometric representation theory; see
\cite{KL1,Lu1,Lu}. 

The key point for us is that $\frac{1}{c(s,t)c(s^{-1},t)}$ 
is a limit $q\to 0$ of 
the corresponding symmetric 
Macdonald's measure-function $\de(s;q,t)$ upon $t\mapsto 1/t$.
This measure makes  the symmetric Macdonald
polynomials pairwise orthogonal. We switch to its nonsymmetric
variant $\mu$ in this paper, the measure-function
that makes the nonsymmetric Macdonald
polynomials pairwise orthogonal.  

In the DAHA approach, the problem is to find 
meromorphic continuations of the
DAHA inner products by presenting them as 
integrals over {\sf double affine residual subtori}. 
The main claims are as follows.   

\vskip 0.3cm


\noindent
{\em The $q,t$-generalization of the 
picking up residues is a presentation of  the standard ($\Diamond$-
invariant) 
inner product in the DAHA polynomial representation as
a finite linear combination of
integrals over double affine residual subtori, where the
measure-function  $\mu$ reduces naturally. This formula
provides the meromorphic continuation of the integral
formula for this inner product from $|t|<1$ for any $|t|\ge 1$.
Its DAHA-invariance and some assumptions about the structure
are expected
to determine 
the corresponding coefficients uniquely.
Upon taking the limit $q\to 0$, we obtain an alternative tool 
for finding the $C_{s_\circ S}$-coefficients for AHA 
including the formal degrees (for the residual points).
}


\section{\sc Basic DAHA theory}
Let $\mm$ be the least natural number
such that  $(P,P)=(1/\mm)\Z.$  Thus
$\mm=|\Pi|$ unless 
$\mm=2 \for D_{2k}$ and $\ \mm=1 \for B_{2k},C_{k}.$

The {\sf double affine Hecke algebra, DAHA\,}, depends
on the parameters
$q, t_\nu\, (\nu\in \{\nu_\al,\al\in R\})\,$ and is defined
over the ring
$\Z_{q,t}\equal\Z[q^{\pm 1/\mm},t_\nu^{\pm 1/2}]$
formed by
polynomials in terms of $q^{\pm 1/\mm}$ and
$\{t_\nu^{\pm1/2}\}.$ It will be convenient to use
$t_\nu =q_\nu^{k_\mu}=q^{\nu k_\nu}$ for $q_\nu=q^\nu$. 

For $\tal=[\al,\nu_\al j] \in \tR,\ 0\le i\le n$, we set
\begin{align}\label{talqal}
&   t_{\tal} =t_{\nu_\al}=q_\al^{k_\al} ,
q_{\tal}=q^{\nu_\al},\, t_i=t_{\al_i}=q_i^{k_i},\,
\rho_k=\frac{1}{2}\sum_{\al\in R_+}k_\al \al.
\end{align}
Using $\rho_\nu\equal \frac{1}{2}\sum_{\nu_\al=\nu}\al$,
we have: $\rho_k=\sum_\nu k_\nu \rho_\nu=k_{\sht}\rho_{\sht}+
k_{\lng}\rho_{\lng}$. 
The standard argument based on the application of $s_i$ for $i>0$
to  $\rho_\nu$ gives that $(\rho_\nu, \al_i^\vee)=1$ for $\nu_i=\nu$ 
and $0$ otherwise for  $i>0$. We obtain that
 $\rho_k=\sum_{i}k_i\om_i$.

For pairwise commutative $X_{\om_1},\ldots,X_{\om_n},$ let
\begin{align}
& X_{\tb}\ \equal\ q^{ j}\prod_{i=1}^nX_{\om_i}^{l_i} 
\iif \tb=[b,j],\ \hw(X_{\tb})\ =\ X_{\hw(\tb)},
\label{Xdex}\\
&\hbox{where\ } b=\sum_{i=1}^n l_i \om_i\in P,\ \, j \in
\frac{1}{ \mm}\Z \ \text{\  and \ } \hw\in \hW.
\notag \end{align}
Technically, $X_b=q^{(x,b)}$ and $X_{\om_i}=q^{(x,\om_i)}$. 
Also, $X_{\al_0}=qX_\vth^{-1}$.
\medskip

Recall that 
$\om_r=\pi_r u_r$ for $r\in O'$ (see above) and 
$\pi_r^{-1}=\pi_{\varsigma(i)}$, where
$\,\varsigma\,$ is the action of $-w_0$ on roots and weights;
we set $X_b^\varsigma=X_{b^\varsigma}.$

\begin{definition}
The double affine Hecke algebra $\HH\ $
is generated over $\Z_{q,t}$ by 
$\h=\lan  \Pi, T_i,\ 0\le i\le n \ran$, subject to the homogeneous
Coxeter relations and the quadratic relations 
$(T_i-t_{i}^{1/2})(T_i+t_{i}^{-1/2})=0$, and by
pairwise commutative $\{X_b, \ b\in P\}$ satisfying
(\ref{Xdex}). The following ``cross-relations" are imposed:

(i)\ \ \  $T_iX_b \ =\ X_b X_{\al_i}^{-1} T_i^{-1}$ if
$(b,\al^\vee_i)=1,\,
0 \le i\le  n$;

(ii)\ \ $T_iX_b\ =\ X_b T_i\ $ if $\ (b,\al^\vee_i)=0
\for 0 \le i\le  n$;

(iii)  $\pi_rX_b \pi_r^{-1} = X_{\pi_r(b)} =
X_{ u^{-1}_r(b)}
 q^{(\om_{\varsigma(r)},b)} \text{ for }  r\in O'$.
\label{double}
\end{definition}

\noindent
The action of $\hW$ in $\R^{n+1}$ is used 
in $(iii)$. Namely:\,
 $\pi_r(b)=\om_r u_r^{-1}(b)=
[u_r^{-1}(b),-(\om_r,u_r^{-1}(b))]$, where $-(\om_r,u_r^{-1}(b))=
(b,-u_r(\om_r))=(b,\om_{\varsigma(r)})$ and
$u_r^{-1}=u_{\varsigma(r)}$.  Recall that
$u_r(\om_r)=w_0(\om_r)=-\om_{\varsigma(r)}$. For instance,
one has: $X_r\pi_r=q^{(\om_r,\om_r)}\pi_r X_{\varsigma(r)}^{-1}$.

The pairwise commutative elements $Y_b$ are as above:
\begin{align}
& Y_{b}\equal
\prod_{i=1}^nY_i^{l_i} \iif
b=\sum_{i=1}^n l_i\om_i\in P,\ 
Y_i\equal T_{\om_i},b\in P.
\label{Ybx}
\end{align}
When acting in the polynomial representation
(see below), they are called {\sf difference
Dunkl operators}. We arrive at the presentation
$\HH\!=\!\lan X_b,T_w,Y_b,q^{\pm 1/\mm},
t_\nu^{\pm 1/2}\ran,\, b\in P,w\in W,\ \,$.
The relations for $\{Y_b\}$ with $\{T_i, X_b\}$ 
result from those for $T_0$ and the relations in 
$\h_X\equal \lan T_i X_b\ran$, where
$1\le i \le n,\, b\in P$. 
The algebra $\h_X$ is isomorphic
to $\h=\h_Y$ under $X_b\mapsto Y_b^{-1}$, $T_w\mapsto T_w$. 
\vskip 0.2cm

{\bf Automorphisms and anti-involutions.}
The following maps can be (uniquely) extended to
automorphisms of $\HH\,$, where 
$q^{1/(2\mm)}$ must be added to $\Z_{q,t}$
(see \cite{C101}, (3.2.10)--(3.2.15))\,:
\begin{align}\label{tauplus}
& \tau_+:\  X_b \mapsto X_b, \ T_i\mapsto T_i\, (i>0),
\ Y_{\om_r} \mapsto X_{\om_r}
Y_{\om_r} q^{-\frac{(\om_r,\om_r)}{2}},
\\
& \tau_+:\ T_0\mapsto  q^{-1}\,X_\vth T_0^{-1},\ \,
\pi_r \mapsto q^{-\frac{(\om_r,\om_r)}{2}}X_{\om_r}\pi_r\
(r\in O'),\notag\\
& \label{taumin}
\tau_-:\ Y_b \mapsto \,Y_b, \ \, T_i\mapsto T_i\, (i\ge 0),\
\ X_{\om_r} \mapsto Y_{\om_r} X_{\om_r}
 q^\frac{(\om_r,\om_r)}{ 2},\\
&\tau_-(X_{\vth})= 
q T_0 X_\vth^{-1} T_{s_{\vth}}^{-1};\ \
\si\equal \tau_+\tau_-^{-1}\tau_+\, =\,
\tau_-^{-1}\tau_+\tau_-^{-1},\notag\\
&\si(X_b)\!=\!Y_b^{-1},\   \si(Y_b)\!=\!
T_{w_0}^{-1}X_{b^{\,\varsigma}}^{-1}T_{w_0},\ \si(T_i)\!=\!T_i (i>0).
\label{taux}
\end{align}
Formally,  $\tau_+^l(H)=q^{\,lx^2/2} H (q^{-\,lx^2/2})$ 
for any $H\in \HH$; this is in the polynomial representation,
which is faithful for generic $q$. For instance,
 $q^{\,lx^2/2} Y_b q^{-\,lx^2/2}=q^{lx^2/2-l(x-b)^2}Y_b=
q^{-b^2/2}X_bY_b$ \ for minuscule $b$, which is direct from the
formula for $Y_b$ (the Gaussian is $W$-invariant). 

In particular, $\si(Y_r)=\tau_{-}^{-1}\tau_+\tau_-^{-1}(Y_r)=
q^{-(\om_r,\om_r)}Y_r^{-1}X_r Y_r$ for $r\in O'$. Also,
$\si(\pi_r)=T_{u_r}^{-1}X_{\varsigma(r)}^{-1}$, which gives that
$\si(Y_r)=T_{u_r}^{-1}X_{\varsigma(r)}^{-1}T_{u_r}$.
See formulas
(3.2.16) and (3.2.22) from \cite{C101}.

The justification is as follows. First,
$\si(\pi_r)=T_{w_0}^{-1}X_{\varsigma(r)}^{-1}T_{w_0}T_{u_r}^{-1}=
T_{u_r}^{-1}X_{\varsigma(r)}^{-1}.$ Second, we represent 
$w_0=v u_{r}$ 
where $v=w^{\varsigma(r)}_0$. Then, $T_{w_0}=T_v T_{u_r}$, where
$v(\om_{\varsigma(r)})=
\om_{\varsigma(r)}$, which gives that
 $T_v$ commutes with $X_{\varsigma(r)}$. 
See (\ref{ururstar}) and  
formula (\ref{diamsi}) below. 
\vskip 0.2cm

We note that 
$T_{w_0}^{-1}T_iT_{w_0}=T_{\varsigma(i)}$
for $i>0$, $T_{w_0}^{-1}T_0T_{w_0}=T_0$ and
 $T_{w_0}^{-1}\pi_r T_{w_0}=\pi_{\varsigma(r)}$. 
Generally,
$\si^2(H)=T_{w_0}^{-1} \varsigma(H) T_{w_0}$, where 
the involution $\varsigma$ is  naturally
extended to an automorphism
of $\HH\ni H$:
$$
X_b\mapsto X_{b^\varsigma},\, Y_b\mapsto Y_{b^\varsigma},\,
T_i\mapsto T_{i^\varsigma},\, \pi_r\mapsto \pi _{r^\varsigma},\ 
b\!\in\! P,\, i\!\ge\! 0,\, r\!\in\! O'.
$$

We obtain that the {\sf projective\,} $PSL_2(\Z)$ due to Steinberg 
acts in $\HH$; it  
generated by $\tau_{\pm}$ subject to the
relations $\tau_+\tau_-^{-1}\tau_+=\si=
\tau_-^{-1}\tau_+\tau_-^{-1}$. This group is isomorphic to the braid 
group $B_3$. W
e note the relation  
$\si \tau_{\pm}=\tau_{\mp}^{-1}\si$.
The automorphism $\si^{-1}$
is the 
{\sf DAHA-Fourier transform}. 

\vskip 0.2cm
All these automorphisms fix $\ t_\nu,\ q$
and their fractional powers, as well as the
{\sf anti-involution\,} $\vph$:
\begin{align}
&\vph:\ 
X_b\mapsto Y_b^{-1},\, Y_b\mapsto X_b^{-1},\,
T_w\mapsto T_{w^{-1}}\, (w\in W),
\label{starphi}
\end{align} 
also sending $\pi_r\mapsto 
T_{u_r}^{-1}X_r^{-1}, T_0\mapsto (X_\vth T_{s_\vth})^{-1}.$

One has for $b\in P$:
\begin{align}\label{vphtaupm}
\vph \tau_+\vph\!=\!\tau_-,\, \vph \si\!=\!\si^{-1}\vph,\,
\vph\si^{-1}(Y_b)\!=\!Y_b, \,\vph\bigl(\tau_+^{-1}(Y_b)\bigr)\!=\!
\tau_+^{-1}(Y_b),
\end{align}
which is direct from the definitions.
Also, for $i\ge 0$ and $r\in O$:
\begin{align}
&\vph(\tau_+(T_i))=\tau_+(T_i),\ 
\vph(\tau_+(\pi_r))= (\tau_+(\pi_r))^{-1}=\tau_+(\pi_{r^\varsigma}).
\label{vphtau}
\end{align} 
For the sake of completeness, let us justify (\ref{vphtau}).
We need  to check the first one only for $i=0$, where
$\tau_+(T_0)=q^{-1}X_\vth T_{s_\vth} Y_\vth^{-1}$ is obviously
$\vph$-invariant. For the 2{\small nd}:  $\pi_r=Y_r T_{u_r}^{-1}=
\pi_{r^\varsigma}^{-1}=T_{u_{r^\varsigma}}Y_{r^\varsigma}^{-1}$.
Applying $\vph$, we obtain the identities 
$T_{u_{r^\varsigma}}^{-1}X_r^{-1}=X_{r^\varsigma}T_{u_r}$, \, 
$X_r T_{u_{r^\varsigma}}=T_{u_r}^{-1} X_{r^\varsigma}^{-1}$ \, and
\begin{align*}
\tau_+(\pi_r)=&\,
q^{-\frac{(\om_r,\om_r)}{2}}X_r\pi_r\,=
q^{-\frac{(\om_r,\om_r)}{2}}X_rY_rT_{u_r}^{-1}\\
=&\,
q^{\frac{(\om_r,\om_r)}{2}}\pi_r X_{r^\varsigma}^{-1}
=q^{\frac{(\om_r,\om_r)}{2}}
Y_rT_{u_r}^{-1}X_{r^\varsigma}^{-1}=
q^{\frac{(\om_r,\om_r)}{2}}
Y_r X_r T_{u_{r^\varsigma}}.
\end{align*}
Therefore  $\vph(X_rY_rT_{u_r}^{-1})=
T_{u_{r^\varsigma}}^{-1}X_r^{-1}Y_r^{-1}$ and we obtain the required.
See formula (3.2.12) in \cite{C101}.

\vskip 0.2cm
The  following anti-involution $\star$ results directly from
the group nature of the DAHA relations. Let 
$
H^\star= H^{-1} \for H\in \{T_{\hw},X_b, Y_b, \pi_r, q, t_\nu\}.
$
To be exact, it is naturally extended to the fractional
powers of $q,t$:
$$
\star:\ t^{\frac{1}{2}}_\nu \mapsto t_\nu^{-\frac{1}{2}},\
q^{\frac{1}{2\mm}}\mapsto  q^{-\frac{1}{2\mm}}.
$$
It commutes with any (anti-)isomorphisms of $\HH$. 
This anti-involution serves the standard inner product in the theory 
of the DAHA polynomial representation $\mathscr{X}$,
 but we will use $\Diamond$ instead. 
For $l\in \Z$, the {\sf anti-involutions} $\Diamond_l$ 
preserve $q,t_\nu$ and send:
\begin{align}\label{diamondef}
&  \Diamond: X_b\!\mapsto\! T_{w_0}^{-1}X_{-w_0(b)}T_{w_0},\,  
Y_b\!\mapsto\! Y_b,\, T_w\!\mapsto\! T_{w^{-1}},\, 
\pi_r\!\mapsto\! T_{w_0}^{-1}
\pi_rT_{w_0},
 \\
\Diamond_l=& q^{l x^2/2}\circ\Diamond\,\circ q^{-lx^2/2}:
X_a\mapsto X_a^{\Diamond},\ \,
Y_b\mapsto q^{l x^2/2}Y_b q^{-lx^2/2}=\tau_+^l(Y_b),\notag
\end{align}
where $b\in P, w\in W, r\in O$. Here, formally $\Diamond(q^{lx^2/2})=
q^{lx^2/2}$; we use
that $x^2$ is $W$-invariant and $\varsigma$-invariant.
Thus, $\Diamond_l$ is the composition
$\tau_+^l\circ \Diamond$. We note that
$\Diamond=\vph\si^{-1}$, 
$\Diamond\circ \tau_{\pm}=\tau_{\pm}^{-1}\circ \Diamond$
and $\Diamond\circ \si=\si^{-1}\circ\Diamond.$

Chapter 3 of \cite{C101} is actually the theory of 
$\vph,\star,\Diamond_{\pm 1}$ and the corresponding
symmetric forms in the polynomial representation and its
Fourier-dual, which is the space generated by delta-functions
at the points $\pi_b(-\rho_k)=b-u_b^{-1}(\rho_k)$ for $b\in P$.
\vskip 0.2cm

Let us provide the counterpart of the symmetries
from (\ref{vphtau}) for $\Diamond$:
\begin{align}
&\Diamond(\si(T_i))=\si(T_i) (i\ge 0),\ 
\Diamond(\si(\pi_r))= (\si(\pi_r))^{-1}=\si(\pi_{r^\varsigma}).
\label{diamsi}
\end{align} 

The first relation is not immediate only for 
$\si(T_0)=T_{s_\vth}^{-1}X_\vth^{-1}.$ One has:
$\Diamond(\si(T_0))=T_{w_0}^{-1}X_\vth^{-1} T_{w_0}T_{s_\vth}^{-1}=
T_{s_\vth}^{-1}X_\vth^{-1} T_{s_\vth}T_{s_\vth}^{-1}=
T_{s_\vth}^{-1}X_\vth^{-1}.$ We use that
$T_{w_0}^{-1}X_\vth^{\pm 1} T_{w_0}=
T_{s_\vth}^{-1}X_\vth^{\pm 1} T_{s_\vth}$, which follows
from (3.2.22) in \cite{C101}, and can be check directly using
that $w_0=u s_\vth$ for $u$ such that $u(\vth)=\vth$; indeed,
$w_0(\vth)=-\vth=s_\vth(\vth)$. We obtain that
$T_{w_0}^{-1}X_\vth^{\pm 1} T_{w_0}=T_{s_\vth}^{-1}T_u^{-1}
X_\vth^{\pm 1} T_u T_{s_\vth}$, where $T_u$ commutes with 
any polynomial of $X_\vth$. 

The second equality is justified as follows. One has: 
$\Diamond(\si(\pi_r))=\Diamond(T_{u_r}^{-1}X_{\varsigma(r)}^{-1})=
T_{w_0}^{-1}X_r^{-1}T_{w_0} T_{u_{\varsigma(r)}}^{-1}=
T_{u_{\varsigma(r)}}^{-1} X_{r}^{-1}$ due to 
$T_{w_0}^{-1}X_{\varsigma(r)}^{\pm 1}
T_{w_0}=T_{u_r}^{-1}X_{\varsigma(r)}^{\pm 1}T_{u_r}.$ 
Alternatively, one can use here and above 
that $\Diamond
=\vph\si^{-1}$.

\section{\sc Polynomial representation}
Its theory is based on the PBW Theorem (actually,
there are $6$ of them for different orderings of $X,T,Y$):
\begin{theorem}[{\sf PBW for DAHA}]\label{PBWDAHA}
Every element in $\HH$ can be uniquely written in the form
\begin{equation}\label{pbwcdaha}
\sum_{a,w,b}C_{a,w,b}\,X_{a}T_{w}Y_b \text{ for } 
C_{a,w,b}\in \C_{q,t},\ a\in P,\, w\in W,\, b\in P^\vee, 
\end{equation}
where $\C_{q,t}=
\C[q^{\pm 1/m},t^{\pm 1/2}]$; actually, 
$\Z_{q,t}$ is sufficient.
 \sq
\end{theorem} 
\vskip 0.5cm

The theorem readily results in the definition of the 
{\sf polynomial 
representation} of $\HH$ in  
$\mathscr{X}\equal\C_{q, t}[X_b]=\C_{q,t}[X_{\om_i}]$. 
Using Theorem \ref{PBWDAHA}, we can  
identify $\mathscr{X}$ with the induced representation
Ind$_\h^{\HH}\,\C_+$, where $\C_+$ is the one-dimensional
module of $\h$ such that $T_{\hw}\mapsto t^{l(\hw)/2}\equal
\prod_\nu t_\nu^{l_\nu(\hw)/2}$. We note that
$t^{l(b)/2}=\prod_\nu t_\nu^{(\rho^\vee_\nu,b)}=
\prod_\nu q^{k_\nu (\nu\frac{\rho_\nu}{\nu},b)}=
q^{(\rho_k,b)}$ for $b\in P_+$.
 
The generators $X_{b}$ act by multiplication; 
$T_i (i\ge 0)$ and $\pi_r (r\in O^*)$ 
act in $\mathscr{X}$ as follows:
\begin{eqnarray}\label{pitpolyn}
&\pi_{r}\mapsto \pi_{r},\ \, 
T_{i}\mapsto t_i^{1/2}s_{i}+\dfrac{t_i^{1/2}-t_i^{-1/2}}
{X_{\al_{i}}-1}(s_{i}-1)\text{ for } t_i=t_{\al_i}. 
\end{eqnarray}
Recall that 
$s_{0}(X_{b})=X_{b}X_{\vth}^{-(b, \vth)}q^{(b, \vth)}$.
The images of $T_i$ for $i>0$ are 
{\sf Demazure-Lusztig operators}. 
\vskip 0.2cm

{\bf DAHA coinvariants.} Generally, they can be defined for 
any anti-involutions of $\HH$ and $\HH$-modules; $\mathscr{X}$
will be considered here.

\begin{definition}\label{def:coinv}
(i) The {\sf Shapovalov
anti-involution} $\varkappa$ of $\HH$ for $Y$
is such that $T_w^{\varkappa}\!=\!T_{w^{-1}}$ and 
the following property holds: for any $H\in \HH$,\,
the decomposition $H=\!\sum c_{awb}\! Y_a^{\varkappa} 
T_{w} Y_b$ exists and 
is unique.

(ii) Given $\varkappa$,
the corresponding {\sf coinvariant} $\varrho$ 
is a functional (a linear map to $\C$)
on $\HH$ such that  $\varrho(H^\varkappa)=\varrho(H)$.
Then 
$\{A,B\}_\varrho\equal \varrho(
A^{\varkappa}\, B)
\,=\,\{B,A\}_\varrho$ and 
$\{HA,B\}_\varrho=\{A,H^\varkappa B\}_\varrho.$ 

(iii) A anti-involution  $\varkappa$ is called {\sf basic} 
if $\varrho$ is unique up to proportionality   
 among the functionals acting 
 via the projection
$\HH\ni H\mapsto H(1)$ onto $\mathscr{X}$. The Shapovalov ones
are basic. 
 \sq
\end{definition}

For Shapovalov $\varkappa$, the coinvariant$\varrho$ 
normalized by the relation
$\varrho(1)=1$ is unique: 
$\varrho(H) =
\,\sum c_{awb}\, \varrho(Y_a)
\varrho(T_{w}) \varrho(Y_b)$, where
$H$ is expanded as in $(i)$. Here 
$\varrho$ is the character
of $\h$ sending $T_{i}\mapsto t_{i}^{1/2}$ for $i\ge 0$
and $\pi_r\mapsto 1$. This formula for $\varrho$ 
works for arbitrary $q,t_\nu$.


An anti-involutions $\varkappa$ is {\sf basic} if and only if
$dim \bigl(\HH/(\j\!+\!\j^\varkappa)\bigr)\!=\!1$, where
$\mathscr{X}\!=\!\HH/\j$ for the left ideal
 $\j\!=\!\{H\in \HH \mid H(1)=0\}$, where 
$1\in \mathscr{X}$ and $H(\cdots)$ is the
action of $H$ in $\mathscr{X}$.
\vskip 0.2cm

The anti-involution $\vph$ from (\ref{starphi})
is a Shapovalov one
due to ``PBW". The corresponding
{\sf evaluation pairing} provides the duality and evaluation 
conjectures practically without calculations; 
see \cite{C6}. We will use sometimes the notation $\{\cdot\}$ 
or $\{\cdot\}_{-\rho_k}$ for it.
The corresponding form $\{A,B\}$ and its restriction to $\mathscr{X}$
are well defined for any $q,t$ and the study of its radical is an 
important tool in the theory of the polynomial representation 
of DAHA.

The anti-involution $\star$, sending  $g\mapsto g^{-1}$ for
$g=X_a,Y_b,T_w,q,t_\nu$, is {\sf basic} for
{\sf generic} $q,t$ but not 
a Shapovalov one. It is
proven in \cite{C101} that
the corresponding  inner product in $\mathscr{X}$ is {\sf unique}
up to proportionality for generic $q,t_\nu$. 

Similarly,  
the anti-involution $\Diamond$ from (\ref{diamondef})
is {\sf basic} for generic $q,t_\nu$ but not a Shapovalov
one. It 
governs the inner product 
in $\mathscr{X}$ making the nonsymmetric Macdonald polynomials
(below) pairwise orthogonal and fixing $q,t_\nu$.
The corresponding bilinear form is the key in the 
DAHA harmonic analysis, including 
the Plancherel formula for $\mathscr{X}$ and its Fourier
image, the representation 
of $\HH$ in delta-functions.  The notation $\lan\cdot\ran$
will be used below for the corresponding coinvariant.

The conjugations $\Diamond_{\pm 1}$
of $\Diamond$ by $q^{\pm x^2/2}$ are
Shapovalov ones.
So the corresponding symmetric form is well-defined for 
any $q,t_\nu$; the notation will be $\lan\cdot \ran_{\pm 1}$.
The radical of the pairing for $\Diamond_1$
is closely related to that 
for $\vph$; they coincide in the rational theory.
The anti-involutions $\Diamond_{\pm 1}$ are the key in the 
difference Mehta-Macdonald formulas and are used to calculate  
the Fourier transform of the DAHA modules $\mathscr{X}q^{\mp x^2/2}$. 
\vskip 0.2cm

{\bf Mu-functions.}
We set
\begin{equation}\label{mutildemu}
\mu(X;q,t_\nu)=\prod_{\tilde{\al}>0}\frac{1-X_{\tilde{\al}}}
{1-t_\al X_{\tilde{\al}}},\ \ \tilde{\mu}(X;q,t_\nu)=
\prod_{\tilde{\al}>0}\frac{1-t_\al^{-1}X_{\tilde{\al}}}
{1-X_{\tilde{\al}}}.
\end{equation}
Recall that
$
\La(\hw)\,\equal\,\tilde{R}_{+}\cap \hat{w}^{-1}(\tilde{R}_{-})
=\{\tal>0\,|\,\hw(\tal)<0\} \for  \hat{w}\in \hat{W};
$
this set consists of $l(\hat{w})$ positive roots.
The following are the key relations
for the functions $\mu,\tilde{\mu}$:
\begin{align}\label{murelations}
&\frac{\hw^{-1}(\mu)}{\mu}=
\frac{\hw^{-1}(\tilde{\mu})}{\tilde{\mu}}=
\prod_{\tilde{\al}\in \Lambda(\hat{w})}
\frac{1-t_\al^{-1}X_{\tilde{\al}}^{-1}}{1-X_{\tilde{\al}}^{-1}}\cdot
\frac{1-X_{\tilde{\al}}}{1-t_\al^{-1}X_{\tilde{\al}}}\\
=&\prod_{\tilde{\al}\in \Lambda(\hat{w})}
\frac{1-t_\al^{-1}X_{\tilde{\al}}^{-1}}{1-t_\al^{-1}
X_{\tilde{\al}}}\cdot
\frac{1-X_{\tilde{\al}}}{\ \,1-X_{\tilde{\al}}^{-1}}= 
\prod_{\tilde{\al}\in \Lambda(\hat{w})}
\frac{t_\al^{-1}-X_{\tilde{\al}}}{1-t_\al^{-1}X_{\tilde{\al}}}.\notag
\end{align}
We see that $\mu/\tilde{\mu}$ is (formally)
a $\hW$-invariant function.
Note that both functions, $\mu$ and $\tilde{\mu}$, are
invariant under the action of $\Pi=\{\pi_r,r\in O\}$ and
under the automorphisms of the {\sf non-affine} Dynkin diagram.
Also, $\frac{\hw^{-1}(\mu)}{\mu}$ is invariant under the
``conjugation" $q\mapsto q^{-1}, t_\nu\mapsto t_\nu^{-1}$,
which sends $X_b\mapsto X_b^{-1}$ and 
$X_{\tal}\mapsto X_{\tal}^{-1}$ (in $\mu$).

The action on functions here and generally is 
$\hw(f(x))=f(\hw^{-1}(x))$; notice $\hw^{-1}$. This results in
the action of  $\hw$ (without $\{\cdot\}^{-1}$) on the indices
of $X_\al$. For instance, $w(X_a)=q^{(w^{-1}(x),a)}=q^{(x,w(a))}=
X_{w(a)}$, $b(X_a)=b(q^{(x,a)})=q^{(-b(x),a)}=
q^{(x-b,a)}=q^{-(b,a)}X_a=X_{[a,-(b,a)]}$.

The $W$-symmetrization of  $\mu$ is essentially
the {\sf Macdonald's function}:
\begin{align}\label{del-fun}
\de\equal&\,\mu\prod_{\al>0}\frac{1-X^{-1}_{\al}}
{1-t_\al X^{-1}_{\al}}:\ 
\sum_{w\in W} w^{-1}\bigl(\frac{\mu}{\de}\bigr)
=\sum_{w\in W}w\Bigl(\prod_{\al>0}\frac{1-t_{\al}X_{\al}^{-1}}
{1-X_{\al}^{-1}}\Bigr)\notag\\
=&
\sum_{w\in W}(-1)^{l(w)}\prod_{\al>0}
\frac{X_{w(\al)}^{1/2}-t_{\al}X_{w(\al)}^{-1/2}}
{X_{\al}^{1/2}-X_{\al}^{-1/2}}=\sum_{w\in W}
t^{l(w)}=P(t_\nu),
\end{align}
where the latter is the Poincare polynomial of $W$.
\vskip 0.2cm

For the sake of completeness, 
let us provide the formula for the constant term ct$(t_\nu)$ 
of $\mu$ (the coefficient of $X^0$):
\begin{align}
&\hbox{ct}(t_\nu)\equal\text{CT}(\mu)=
\prod_{\al\in R_{+}}\prod_{i=1}^{\infty}
\frac{(1-q^{(\al,\rho_k)+i\nu_\al})^{2}}
{(1-t_\al q^{(\al,\rho_k)+i\nu_\al})(1-t_\al^{-1}
q^{(\al,\rho_k)+i\nu_\al})}.\notag
\end{align}
To define ct$(\mu)$ we expand $\mu$ in terms of $t_\nu$ and 
$X_{\tal}$ for $\tal>0$. Then $\hbox{ct}(t_\nu)$ is 
an element in $\Z[t_\nu][[q]]$.
We will use this formula mainly for $t_\nu^{-1}$ instead of $t_\nu$.
Here, as above,
 $q^{(\al,\rho_k)+i\nu_\al}=q_\al^{(\al^\vee, \rho_k)+i}$.
\vskip 0.2cm

{\bf Jackson integrals.}
We mostly follow here \cite{ChSel} and \cite{C101}.
Let us fix $\xi\in \C^n$ and set
$X_a(bw)\equal q^{(b+w(\xi),a)}$ for $bw\in \hW$.
For instance,
$\mu(0)\!=\!\mu(q^\xi)$ and  $\bigl(\hw^{-1}(\mu)/\mu\bigr)(0)\!=\!
\mu(\hw)/\mu(0).$ 

Provided the convergence,
the
{\sf Jackson integral} is defined as $J_\xi(f)=J(f;\xi)\equal 
\sum_{\hw\in \hW}f(\hw)
\mu(\hw)/\mu(0)$. It is a sum,
but can be expressed as a difference of some
integrals (see below).
From formula (\ref{murelations}):
\begin{align}\label{muxi}
&\mu(\hat{w})/\mu(0)=
\prod_{[\al,j\nu_{\al}]\in \Lambda(\hat{w})}
\frac{t_\al^{-1}-q_{\al}^{(\al^\vee,\xi)+j}}
{1-t_\al^{-1}q_{\al}^{(\al^\vee,\xi)+j}}\,.
\end{align}
Recall that
 $q^{(\al,\xi)+j\nu_{\al}}=q_\al^{(\al^\vee, \xi)+j}$. 
The key property of these ``integrals" is
that $J(f;\xi)$ does not depend on $\xi$ up to a coefficient
of proportionality (serving all $f$) for the spaces
$\mathscr{X}$ and $\mathscr{X}q^{x^2/2}$. The coefficient of
proportionality is formula (3.5.14) from \cite{C101}.  This is
due to the uniqueness of coinvariants for {\sf basic} $\varkappa$. 

Also,  $J(f;\xi)=0$ in these spaces if
$(\tal,\xi)=0$ for some $\tal$, where the pairing is affine:
$([\al,j],\xi)=(\al,\xi)+j$.  
More exactly, we have the following lemma, which will be used later.

\begin{lemma}\label{lem:stabxi} For generic $t$ and 
sufficiently general $\xi$:\, $J(f;\xi)=0$
in any spaces of functions 
 if $(\tal,\xi)=0$ for at least one $\tal\in \tR$.
\end{lemma}
{\it Proof.}
Applying a proper $\tu\in \tW$ to $\xi$ (for the affine action)
we can assume that such 
$\tal$ form a root subsystem  with simple 
roots $\al_{i'}$ for $i'$ from a subdiagram of the affine Dynkin
diagram of $\tR$.  
One has: 
$\La(\hat{w}s_i)=\La(\hat{w})\cup \{\al_i\}$ for this $i$, and
$\frac{\mu(\hat{w}s_i)}{\mu(0)}=
\frac{\mu(\hat{w}s_i)}{\mu(0)} \frac{t_i^{-1}-1}{1-t_i^{-1}}=
-\frac{\mu(\hat{w}s_i)}{\mu(0)}.
$
Recall that $bw$ is considered here as the point $q^{w(\xi)+b}$. Thus,
the Jackson summation is identically zero upon the
restriction to any
right coset $\{\hw W'\}$ for 
the Weyl group $W'$ generated  by $s_{i'}$. \sq
\vskip 0.2cm

The following modification of $J(f;\xi)$ 
is needed for  $\xi=-\rho_k$:\, we set $J(f;-\rho_k)\equal
\sum_{\pi_b}f(\hw)
\mu(\pi_b)/\mu(0)$, where $b\in P$. For an explicit formula, let 
$\La'(\hw)=
\{[-\al,j\nu_\al]
\mid \tal=[\al,j\nu_\al]\in \La(\hw)\}$. 
We follow Section 4 of \cite{ChSel}.
Recall that $b=\pi_b u_b$
for minimal $u_b$ such that $b_-=u_b(b)$, $b_+=w_0(b_-)$ and
$-b_-=b_+^\varsigma$. 
From (3.1.17) in \cite{C101}: 
\begin{align}\label{La-pi-prim}
\La'(\pi_b)=\Bigl\{[\al,j\nu_\al] \text { s.t. } \al\in R_+,\, 
&-(b_-,\al^\vee)>j>0 
\text{ if } u_b^{-1}(\al)\in R_-\notag\\
\text{ or }\, &-(b_-,\al^\vee)\ge j>0 
\text{ if } u_b^{-1}(\al)
\in R_+\Bigr\}.
\end{align}
Then  $\mu(0)/\mu(\hw)$ is well-defined
for any $\hw$ and it is non-zero if and only if $\hw=\pi_b$ for 
$b\in P$, which follows from
\begin{align}\label{muphoval}
&\frac{\mu(\hat{w})}{\mu(0)}=
t^{-l(\pi_b)}\!\!\!\!\!\prod_{[\al,j\nu_{\al}]\in \Lambda'(\hat{w})}
\frac{1-t_\al q_{\al}^{(\al^\vee\!,\,\rho_k)+j}}
{1-t_\al^{-1}q_{\al}^{(\al^\vee\!,\,\rho_k)+j}}, \
t^{l(\hw)}\equal\prod_\nu t_\nu ^{l_\nu(\hw)}.
\end{align}

As an application, we
obtain the following {\sf Jantzen-type filtration} in 
$\HH$-modules $\f_\xi$ 
linearly generated by the characteristic
functions at points $\hw$. The action of
$T_i,\pi_r$ is
 dual to that in terms of Demazure-Lusztig
operators; see \cite{C101} and Theorem \ref{thm:normde} below. 

Here $\xi$ can be arbitrary. It is deformation 
$\xi_\ep=-(1+\ep)\xi$ becomes generic for small $\ep$
and we can define $J(f; \xi_\ep)$ and 
find
$\ell_1>0$ such that $J_1(f)=\ep^{\ell_1} J(f; \xi_\ep)$ is 
non-singular and nonzero identically.
The first term of this filtration will be then the
$\HH$-submodule of $\f_\xi$ generated by the characteristic
functions where $J_1$ vanishes. 
Inside this module,
consider  $J_2(f)=\ep^{\ell_2} J(f; \xi_\ep)$ for
some $\ell_2<\ell_1$ that is nonzero; the second module will be
the span of characteristic functions
where $J_2$ vanishes. Continue by induction.

As an example, let $\xi=-\rho_k$. Then for generic $q,t_\nu$  
the last submodule will be
the Fourier transform of $\mathscr{X}$. This procedure can be
applicable to any $q,t_\nu$; then $\mathscr{X}$ will be 
decomposed further. 

\vskip 0.2cm
{\bf Affine symmetrizers.}
We continue to assume that $0<q<1$ and use the notation
$t_{\al}=q_\al^{k_\al}$. 
Let
$$\hat{\mathscr{P}}_+(f)\equal
\sum_{\hat{w}\in \hat{W}}t^{-l(\hat{w})/2}
T^{-1}_{\hat{w}}(f),\ \
\hat{\mathscr{I}}_+(f) \equal\sum_{\hat{w}\in \hat{W}}
\hat{w}(\tilde{\mu}f).
$$

Also, the affine Poincar\'e series,
is defined as $\hat{P}(t_\nu)=\sum_{\hat{w}\in \hat{W}}
t^{l(\hat{w})}$; it is $\frac{|\Pi|}{(1-t)^{n}}
\,\prod_{i=1}^n\frac{1-t^{d_i}}{1-t^{d_i-1}}$
in terms
of the {\sf degrees} $d_i$ in the simply-laced case.
See Theorem 2.8 from \cite{CMa}.
In any module over $\h=\lan T_w, Y_b\ran$, the operator
$\hat{\mathscr{P}}_+/\hat{P}(t_\nu^{-1})$ is a projection onto
the space of {\sf spherical vectors}
defined as follows: $\{v \mid T_{\hw}v=t^{l(\hw)/2}v\}$.
This is provided the convergence of $\hat{\mathscr{P}}_+$
and when $\hat{P}(t_\nu^{-1})\neq 0$.

All constructions below can be extended to the
{\sf minus-symmetrizers} (generally, to arbitrary
characters of $\h_Y$), but we will
stick to $\hat{\mathscr{P}}_+$. Recall that $0<q<1$.


\begin{theorem}\label{DIAZERO}
Let us assume that $\mathscr{X}$ has a nonzero symmetric
form $\lan f,g \ran$ with the anti-involution $\Diamond$
normalized by $\lan 1,1\ran=1$. Given any
$f,g\in \mathscr{X}$, $\lan f,g\ran$ is a rational
function in terms of $q,t$. Provided that  $\Re k_\nu<0$
and $|\Re k_\nu|$ are sufficiently large (depending on $f,g$), 
\begin{align}\label{psymformula}
\lan f,g\ran=
t^{-l(w_0)/2}\,\hat{\mathscr{P}}_{+}(f\, T_{w_0}(g^\varsigma))/
\hat{P}(t_\nu^{-1}),
\end{align}
where $\hat{\mathscr{P}}_{+}(f)$ is a constant for 
$f\in \mathscr{X}$ assuming the convergence. Thus, 
formula (\ref{psymformula}) supplies any $\HH$-quotient of 
$\mathscr{X}$ with a partially defined (when converges!)
bilinear symmetric form associated with $\Diamond$, which
satisfies the
normalization condition 
$\lan1,1\ran=1$. \sq
\end{theorem} 

This is Theorem 2.16 from \cite{CMa}.
The following Theorem is an adjustment of some of the 
claims from Theorems 2.5, 2.6, 2.11 there. 
 
\begin{theorem}\label{LEVZERO}
(i) We assume that $\Re k_\nu<0$ for all $\nu$.
Given $a_+\in P_+$, the sums
$\hat{\mathscr{I}}_{+} (X_{a})$ 
absolutely converge for any $a\in W(a_+)$ if and only if\,
$|t^{-l(a_+)}q^{-(a_+,\om_i)}|<1$ for all $i=1,\ldots,n$,
where $l(a_+)=2(\rho^\vee, a)$. 
Equivalently, the conditions become $\Re(2\rho_k+a_{+},\om_i)<0$
in terms of $k_\nu$.

(ii) Employing the formulas
from (\ref{pitpolyn}) for $T_{\hw}$, 
the expansion  $\hat{\mathscr{P}}_+=\sum_{\hw\in\hW}a_{\hw}\hw$ 
is with $t$-meromorphic functions $a_{\hw}$. This holds
by construction
for $\hat{\mathscr{I}}_+$. As formal series and
as operators acting in $\mathscr{X}$:\ 
$\hat{\mathscr{P}}_+
=\text{ct}(t_\nu^{-1})\hat{\mathscr{I}}_+$,
where ct$(t_\nu^{-1})$ is the constant term 
 and the conditions from $(i)$ are imposed if these
operators act in $\mathscr{X}$.

(iii) 
Let $l>0$. The operators 
$\hat{\mathscr{I}}_+$ and  
$\hat{\mathscr{P}}_+$ 
converge absolutely at any given 
$f\in \mathscr{X}q^{\,lx^{2}/2}$
for any $k$ for $\hat{\mathscr{I}}_+$ and under the
constraints 
$\Re (h_{k}^{\lng}), \Re (h_{k}^{\sht})<1$ for 
the operator $\hat{\mathscr{P}}_+$. Here  
$h_{k}^{\sht}=
(\rho_k,\vth)+k_{\sht}$, $h_{k}^{\lng}=
(\rho_k,\th^\vee)+k_{\lng}$ and $\th$
is the maximal root in $R_+$. For instance,
this constraint is $\Re k<\frac{1}{h}$ for  $h_k=k h$ 
in the simply-laced case, where $h$ is the 
Coxeter number. Then
$\hat{\mathscr{P}}_+=\text{ct}(t_\nu^{-1})\hat{\mathscr{I}}_+.$
\sq
\end{theorem}

The convergence conditions in $(i)$ follow directly from
(\ref{murelations}). 
We note that $\hat{\mathscr{P}}_+(f)$ is regular
by construction for $f\in \mathscr{X}q^{\,lx^{2}/2}$ but is 
well-defined only 
for $\Re k<1/h$;  $\hat{\mathscr{I}}_+(f)$ is well-defined 
for any $k$ but  has poles. For instance, their proportionality
gives that the latter has no poles for $\Re k<1/h$, which is
far from obvious from its definition.

The adelic version of this argument is expected to provide
an alternative approach to the fact that the 
Langlands formula for the inner product
of {\sf pseudo-Eisenstein series} has no singularities
due to the Dedekind zeta-functions. See \cite{KO, MHO}.
\vskip 0.2cm

{\bf E-polynomials.}
One of the key results in the DAHA theory is that the norms of
{\sf nonsymmetric Macdonald polynomials}
under the spherical
normalization are $\frac{\mu(0)}{\mu(\hat{w})}$.
For generic $q,t$, they are defined  as follows:
$$
\e_b\equal
E_b/E_b\bigl(q^{-\rho_k}\bigr),\ Y_a(E_b)=q^{(a,-\pi_b(\rho_k))}E_b=
q^{(a,-b+\u_r^{-1}(\rho_k))}E_b,\ b\in P.
$$
The normalization of $E_b$ is $E_b=X_b+(\text{lower terms}).$
The following formula is based on the technique of intertwiners
and relations (\ref{vphtau}), (\ref{diamsi}) above.
For generic $q,t$ and  $b,c\in P$:
\begin{align}\label{nsymnorm}
t^{-l(w_0)/2}\, \text{\rm ct}\bigl(\e_{b}\,T_{w_0}(\e_{c})\,
\mu(X;q,t_\nu)\bigr)/\text{\rm ct}(t_\nu)=
\de_{b,c}\,\mu(q^{-\rho_k})/\mu(\pi_b).
\end{align}
This is essentially Corollary 3.4.1 from \cite{C101}, 
where the anti-involution $\Diamond$ occur there in
formula  (3.4.22). 

Using (\ref{nsymnorm}), we obtain a direct demonstration
of the fact  that the
coinvariant $\lan \cdot \ran$ associated with $\Diamond$
is a meromorphic function for any $k_\nu$;
see (\ref{def:coinv}). Indeed,
$\lan f \ran$ must be proportional to ct$(f\mu)$ for any
Laurent polynomial $f$ due to the uniqueness of the coinvariant
for $\Diamond$. The coefficient of proportionality
is explicit. Then we express $f$ via $\{E_b\}$ and use
that  $\lan E_b\ran =0$ for $b\neq 0$.  

Actually, the proof of (\ref{nsymnorm}) contains the justification
of the uniqueness of $\lan \cdot \ran$.
Let us  extend this formula and its proof 
 to  general $Y$-induced representations. 
\vskip 0.2cm

\section{\sc Induced modules}\label{sec:indmod}
The technique of intertwiners and the theory of basic coinvariants
can be naturally extended to $Y$-induced $\HH$-modules. We mostly
follow \cite{C101,CMa}.

Given $\xi\in \C^n$, the induced representation
$\i_\xi$ is defined as a unique (up to isomorphisms)
 $\HH$-module induced 
from the character
 $\tilde{\xi}$ of the algebra $\C[Y_b, b\in P]$ 
defined  as follows:
$\tilde{\xi}(Y_b)=q^{-(\xi,b)}$.
 In the main examples, $\xi$ depend
of $q$ and $t_\nu$, which are considered as nonzero numbers or 
as formal parameters.
 
As a vector space, $\i_\xi$ is naturally isomorphic
to the affine Hecke algebra
 $\h_X=\lan T_w,X_b\ran$. It is $Y$-semisimple with simple 
$Y$-spectrum
if and only if $q^{\hw(\rho_k)}\neq q^{\rho_k}$
for any id$\neq \hw\in \hW$.

The module $\mathscr{X}$ is a canonical quotient of $\i_\xi$
for $\xi=-\rho_k$. 
We will mostly assume that $0<q<1$ and
it is {\sf generic} with respect to 
$t_\nu$:\
$q^{m}\not\in t_\nu^{\Z}$ for $m\ge 1$. Then $\i_{-\rho_k}$ is
semisimple  when and only when
$w(\xi)\neq \xi$ modulo $2\pi \imath\aa$  for any $w\in W$. 

The $Y$-spectrum of $\i_\xi$ for any $\xi$ 
is $\{q^{w(\xi)+a}\}$, where $a\in P, w\in W$: 
the spaces of {\sf pure} eigenvectors are 
$\{v \mid Y_b(v)=q^{-(b,a+w(\xi))}v,\  b\in P\}$. They are nonzero 
for any $a,w$ and 
the corresponding {\sf generalized spaces} of eigenvectors
linearly generate $\i_\xi$ for any $\xi$.  
This module is irreducible if and only if
$q_\nu^{(\al^\vee,\xi)}\not\in \{t_\nu 
q_\nu^{\Z}\}$ for
any $\al\in R$ and $\nu=\nu_\al$.

\comment{
We will extend $\tilde{\xi}$ to 
$\h_Y$ as follows: $\tilde{\xi}(T_wY_b)=t^{l(w)/2}q^{-(\xi,b)}$;
it is a homomorphism when $\xi=-\rho_k$, but not in general.
Actually, any ``trace" $\chi$ of $\H=\lan T_w, w\in W\ran$,
can be taken here instead of $\chi(T_w)=t^{l(w)/2}$ in the
definition of the coinvariants for $\i_\xi$. More exactly,
we need the property $\chi(T_{w^{-1}})=\chi(T_w)$ for $w\in W$.
 The main examples are 
the characters of finite-dimensional 
representations of $\H$. So we generally set:
$\tilde{\xi}(T_wY_b)\equal\chi(T_w)q^{-(\xi,b)}$.
}

Given Shapovalov $\varkappa$ and $ \xi\in \C^n$, 
the coinvariants $\varrho$ 
are defined
by the relations $\varrho(H^\varkappa)=\varrho(H)$,\,
$\varrho(H Y_b)=\tilde{\xi}(Y_b)\varrho(H)$ and 
$\varrho(T_w)=\tau(T_w),$ where $\tau: \H\to \C$ is an
arbitrary linear map satisfying the relation
$\tau(T_w)=\tau(T_{w^{-1}}).$ The simplest choice
is $\tau(T_w)=t^{l(w)/2}$ for $w\in W$. 

One has then:\
$ \varrho((Y_a^\varkappa) T_w Y_b)=
\tilde{\xi}(Y_{a+b})\tau(T_{w})=q^{(\xi,a+b)}\tau(T_w).$
We see that given Shapovalov $\varkappa$,$\tau$
and an arbitrary $\xi$, there exists a unique
coinvariant up to proportionality. 

The anti-involutions 
$\star, \Diamond$ do not require a choice of $\tau$
for their definition and the uniqueness.
They are {\sf basic} for generic $q,t_\nu,\xi$,
i.e. the corresponding coinvariant is unique under
the normalization $\varrho(vac)=1$. We will prove this below
in process of obtaining the
norm-formula in $\i_{\xi}$ for generic $\xi$. This will be 
 based on the technique of
intertwiners. 

The notation for
the coinvariants with $\xi$ for $\Diamond_l$ will be
 $\lan \cdot \ran_{l,\xi}$; 
we write $\lan \cdot\ran_\xi$ for $l=0$, and $\lan \cdot \ran_l$
for the polynomial representation. 
\vskip 0.2cm

{\bf The norm-formula}.
We follow 
Theorem 3.6.1 from \cite{C101} and (3.6.23). It was
stated there for the anti-involution $\ast$; we adjust it
accordingly and change the proof. 
The next theorem includes the uniqueness of 
$\lan \cdot \ran_\xi$ above for $\Diamond$ and
for generic parameters. 
We set:
\begin{align}\label{interG}
&\Phi_i\!=\!T_i\!+\!\frac{t_i^{1/2}\!-\!t_i^{-1/2}}{X_{\al_i}-1},\ 
\phi_i\!=\!t_i^{1/2}\!+\!\frac{t_i^{1/2}\!-\!t_i^{-1/2}}{X_{\al_i}-1}
\!=\!
\frac{t_i^{1/2} X_{\al_i}\!-\!t_i^{-1/2}}{X_{\al_i}\!-\!1},\\
S_i=&\,\phi_i^{-1}\Phi_i,\,G_i=\Phi_i\phi_i^{-1},
S_{\hw}=\,\pi_r S_{i_\ell}\cdots
S_{i_1},\ G_{\hw}=\,\pi_r G_{i_\ell}\cdots
G_{i_1}, \notag
\end{align} 
where $0\le i\le n,\, \hw=\pi_r s_{i_\ell}\cdots
s_{i_1}$; recall that $X_{\al_0}=q X_{\vth^{-1}}.$
The decomposition
of $\hw\in \hW$ is not necessarily reduced here
because $S_i^2=1=G_i^2$ for $0\le i\le n$. This relation and
the fact that $S,G$ do not depend on the choice of the
reduced decomposition
follow from
the symmetries $S_{\hw} X_b=X_{\hw(b)} S_{\hw}$ 
and $G_{\hw} X_b=X_{\hw(b)} G_{\hw}$ for $\hw\in \hW$.
We obtain that $S_i^2$ is a rational function in terms of $X_b$ and
$S_i^2(1)=1$ in $\mathscr{X}$, which gives that $S_i^2=1$ and
$G_i^2=\phi_i S_i^2 \phi_i^{-1}=1.$ 

We will need actually   $\hat{S}_{\hw}\equal\si(S_{\hw})$,
$\hat{G}_{\hw}\equal\si(G_{\hw})$.
One has:
$\hat{S}_{\hw} Y_b=Y_{\hw(b)} \hat{S}_{\hw}$ for $\hw\in \hW$,
and the same symmetry holds for $\hat{G}$.  

Accordingly, we set $f_{\hw}\equal \hat{S}_{\hw}(v),$
$e_{\hw}\equal \hat{G}_{\hw}(v),$
where $v=vac$ is the cyclic generator of $\i_\xi$, $\hw\in \hW$.
To obtain explicit formulas for $f_{\hw}, e_{\hw}$
in terms of $\pi_r, T_i$, let
$$
S_i(c)\!=\!
\frac{T_i\!+\!(t_i^{1/2}\!-\!t_i^{-1/2})/(X_{\al_i}(q^c)-1)}
{t_i^{1/2}\!+\!(t_i^{1/2}\!-\!t_i^{-1/2})/(X_{\al_i}(q^{-c})-1)},\ 
G_i(c)\!=\!
\frac{T_i\!+\!\frac{t_i^{1/2}\!-\!t_i^{-1/2}}{X_{\al_i}(q^c)-1}}
{t_i^{1/2}\!+\!\frac{t_i^{1/2}\!-\!t_i^{-1/2}}{X_{\al_i}(q^{c})-1}}.
$$
Here $c\in \C^n$. Using the affine action $bw(\!(z)\!)=w(z)+b$:
\begin{align}\label{e-inter}
f_{\hw}\!=&\si\bigl(\pi_r S_{i_\ell}(c_\ell)
\cdots S_{i_1}(c_1)\bigr)
(v),\ \,
e_{\hw}\!=\!\si\bigl(\pi_r G_{i_\ell}(\!(c_\ell)\!)\cdots 
G_{i_1}(c_1)\bigr)
(v)\notag\\
&\text{for }  \hw=\pi_r s_{i_\ell}\cdots s_{i_1},\ 
c_1\!=\!\xi,\ c_2\!=\!s_{i_1}(\!(c_1)\!), \ldots, 
c_\ell\!=\!s_{i_{\ell-1}}(\!(c_{\ell-1})\!).
\end{align}
These formulas justify that $\{e_{\hw},f_{\hw}\}$
 are well-defined and nonzero for generic $q,t\in \C^\ast$.
One has for $\hw=bw,\, b\in P,\, w\in W$:
$$
Y_a(f_{\hw})=q^{-(a,b+w(\xi))}f_{\hw} \text{\, and \,} 
\i_\xi=\oplus_{\hw\in \hW}\C f_{\hw}.
$$
The same relations hold for $\{e_{\hw}\}.$
\vskip 0.2cm

For $\xi=-\rho_k$ and generic $q,t_\nu$, the module $\i_\xi$ has
a canonical quotient obtained by imposing additional relations 
$T_{w} (vac)=t^{l(w)/2}vac$ for $w\in W$, which is
$\mathscr{X}$. The elements
 $e_{\pi_b},f_{\pi_b}$ and their images in $\mathscr{X}$
 are well-defined generic
$q,t$. We note that we used the following  normalization
in formulas  (3.3.42), (3.3.44) from \cite{C101}:
$$\hat{E}_b=\tau_+
\bigl(\pi_r G_{i_\ell}(c_\ell)\cdots G_{i_1}(c_1)\bigr)(1)
\text{ for } b\in P.
$$
The relation to spherical polynomials is:
$\e_b=q^{(\rho_k+b_+,b_+)}\hat{E}_{\pi_b}$ for $b\in P$,
which results from formula (\ref{vphtau}).
\vskip 0.2cm

\comment{
We will impose the {\sf spherical normalization}: 
$\{e_{\hw}\}_\xi=1$; the 
coinvariant $\{\cdot\}_\xi$ is for the
Shapovalov anti-involution $\vph$ and the character
$\tilde{\xi}: Y_a\mapsto  q^{-(a,\xi)}$:
$\{Y_a^\vph\, T_w\, Y_b\}_\xi=\tilde{\xi}(Y_a)\tilde{\xi}(Y_b)
t^{l(w)/2}$.
}

The following norm formula in $\i_\xi$ is actually 
the fundamental fact that the DAHA-Fourier transform of $\i_\xi$
is the corresponding Delta-representation. The Fourier-images 
of $f_{\hw}, e_{\hw}$ become the corresponding characteristic
and delta-function at $\hw=bw$,
where $\hw$ is considered as the point $q^{w(\xi)+b}$.

Concerning the spherical normalization $\{\cdot\}_\xi=1$ for
the evaluation
coinvariants $\{\cdot\}_\xi$, 
one needs to calculate its change under the 
action of $\tau_+(S_{\hw})$, which follows 
Proposition 6.6 from \cite{C6}. The simplest choice
of $\tau$ is   $\tau(T_w)=t^{l(w)/2}$
for $w\in W$; however $\tau(T_0)$ will then depend on $\xi$,
which makes the final formula somewhat more involved than that for 
$\xi=-\rho_k$ with $\{\cdot\}$ acting via $\mathscr{X}$. 

The next theorem is the calculation
of change of the norms is mostly parallel to 
Theorem 3.6.1 from \cite{C101} and (3.6.23). They were 
for the anti-involution $\ast$; we will do this 
for $\Diamond$ and with some improvements.

\begin{theorem}\label{thm:norm}
For generic $\xi,q,t_\nu$ 
let $\lan f \ran_\xi$ be the coinvariant for $\Diamond$
acting via $\i_\xi$ normalization by the
condition $\lan vac \ran_\xi=1$. Then:
$$
\lan f_{\hu}^\Diamond \,f_{\hw}\ran_\xi=
\de_{\hu,\hw}\,\mu(\hw)/\mu(0),\ 
\lan e_{\hu}^\Diamond \,e_{\hw}\ran_\xi=
\de_{\hu,\hw}\,\mu(0)/\mu(\hw),
$$
for any $\hu, \hw\in \hW$, where $\de$ is the Kronecker 
delta and $0=\text{\rm id}\in \hW$ is
considered as $q^{\xi}$. In particular, such
 $\lan H\ran_\xi$ is unique and its values at 
$H=X_a T_w Y_b\in \HH$ are rational in terms of
$q^{(\xi,\al)}$ for $\al\in R$ and fractional
powers of $q,t_\nu$.
\end{theorem}
{\it Proof.} It is based on the formulas in (\ref{diamsi}) coupled
with the identity $S_i^2=1$ for $0\le i\le n$. We set
$\psi_i=\si(\phi_i)=
t_i^{1/2}+\frac{t_i^{1/2}-t_i^{-1/2}}{Y_{\al_i}^{-1}-1}=
\frac{t_i^{1/2} Y_{\al_i}^{-1}-t^{-1/2}}{Y_{\al_i}^{-1}-1},
$ where
$Y_{\al_0}=q^{-1} Y_\vth^{-1}$. Then $\hat{S}_i=
\psi_i^{-1}(\si(T_i)+\frac{t_i^{1/2}-t_i^{-1/2}}{Y_{\al_i}^{-1}-1})$,
\begin{align}\label{Gsquare}
&\hat{S}_i^\Diamond=\psi_i \hat{S}_i \psi_i^{-1}=
\frac{t_i^{1/2} Y_{\al_i}^{-1}-t_i^{-1/2}}{Y_{\al_i}^{-1}-1}
\Bigl(\frac{t_i^{1/2} Y_{\al_i}-t_i^{-1/2}}{Y_{\al_i}-1}
\Bigr)^{-1}\hat{S}_i\\
&=
\frac{t_i^{1/2}-t_i^{-1/2}Y_{\al_i}}
{t_i^{-1/2}-t_i^{1/2} Y_{\al_i}}\hat{S}_i=\frac{1-t_i^{-1}Y_{\al_i}}
{t_i^{-1}- Y_{\al_i}}\hat{S}_i=\frac{t_i^{-1}-Y_{\al_i}^{-1}}
{1-t_i^{-1}Y_{\al_i}^{-1}}\hat{S}_i.
\end{align}
Also, $(\si(\pi_r))^\Diamond=
\si(\pi_{\varsigma(r)})$ for $r\in O'$. We arrive at the relations
\begin{align*}
&\lan\, (\hat{S}_i f_{\hu})^\Diamond \ \hat{S}_i f_{\hw}\,\ran_\xi=
\lan\, f_{\hu}^\Diamond \ (\hat{S}_i^\Diamond 
\hat{S}_i)\, f_{\hw}\,\ran_\xi=\\
&\frac{t_i^{-1}-Y_{\al_i}^{-1}}
{1-t_i^{-1}Y_{\al_i}^{-1}}(Y\!\mapsto\! q^{-\hw(\xi)})
\lan\, f_{\hu}^\Diamond \, f_{\hw}\,\ran_\xi=
\frac{t_i^{-1}-q^{(\al_i,\hat{w}(\xi))}}
{1-t_i^{-1}q^{(\al_i,\hat{w}(\xi))}}
\lan\, f_{\hu}^\Diamond \, f_{\hw}\,\ran_\xi,
\end{align*}
where $([\al,j\nu_\al],z)=(\al,z)+j\nu_\al$, which is needed
here for $\al_0=[-\vth,1]$. Similarly,
$\lan\, (\hat{G}_i e_{\hu})^\Diamond \ \hat{G}_i e_{\hw}\,\ran_\xi=
\frac{1-t_i^{-1}q^{(\al_i,\hat{w}(\xi))}}
{t_i^{-1}-q^{(\al_i,\hat{w}(\xi))}}
\lan\, e_{\hu}^\Diamond \, e_{\hw}\,\ran_\xi.$ 
Using $\La(\hw)=\{\,\al_{i_1},s_{i_1}(\al_{i_2}),
s_{i_1}s_{i_2}(\al_{i_3})\ldots,
\hw^{-1}s_{i_\ell}(\al_{i_{\ell-1}})\,\}$ for a
reduced decomposition
$\hw=\pi_r s_{i_\ell}\cdots s_{i_1}$ (formula (3.1.10) from
\cite{C101}), we obtain:
\begin{align}\label{ediamu}
&\lan\, f_{\hu}^\Diamond \ f_{\hw}\,\ran_\xi=
\de_{\hat{u},\hat{w}}\!\!
\prod_{[\al,j\nu_{\al}]\in \Lambda(\hat{w})}
\frac{t_\al^{-1}-q_{\al}^{(\al^\vee,\xi)+j}}
{1-t_\al^{-1}q_{\al}^{(\al^\vee,\xi)+j}}=\frac{\mu(\hat{w})}{\mu(0)},
\end{align}
and its reciprocal for $\lan\, e_{\hu}^\Diamond \ e_{\hw}\,\ran_\xi$.
\sq

\vskip 0.2cm
We will interpret this theorem as the Plancherel formula
for the DAHA-Fourier transform of $\i_\xi$. Let 
$\f_\xi=\oplus_{\hw\in \hW} \C\chi_{\hw}$ for the
characteristic functions $\chi_{\hw}$ at $\hw=bw$
considered as points $q^{w(\xi)+b}$. It is a
module over the smash product of $\C[X_a, a\in P]$ and
the group algebra $\C\hW$. The action is $S_{\hu}(\chi_{\hw})=
\chi_{\hu\hw}$ and $X_a(\chi_{\hw})=X_a(\hw)\chi_{\hw}$ 
for $\hu,\hw\in \hW$
and $a\in P$. Here, as above, $X_a(\hw)=q^{(a,w(\xi)+b)}.$ 

The action of $\HH$ in $\f_\xi$ is obtained when we use the
action of $S_{\hw}$ to define that of $T_{\hw}$, namely
the formulas 
$T_i=\phi_i S_{i}-\frac{t_i^{1/2}-t_i^{-1/2}}{X_{\al_i}-1}$,
$S_{\pi_r}=\pi_r$. The resulting action will involve the denominators 
in terms of $X$, so we need to assume that $\xi, q,t_i$ are
in a general position when applying them to $\chi_{\hw}$.
See  formula (3.4.10) from \cite{C101}.

\begin{theorem}\label{thm:normde}
For generic $\xi,q,t_\nu$ and any $\hw\in \hW$,
the $\C$-linear map $F: f\mapsto \hat{f}$ sending
 $f_{\hw}\mapsto \chi_{\hw}$ induces
the automorphism $\si^{-1}$ for $H\in \HH$: $FH=\si^{-1}(H)F$.  
The inner product in $\f_\xi$ given by
the formula $(f,g)=\sum_{\hw\in \hW}f(\hw)g(\hw)\mu(\hw)/\mu(0)$
corresponds to the anti-involution $\si^{-1}\circ\Diamond\circ\si=
\Diamond\circ\si^2$. Here $f(bw)=f(q^{w(\xi)+b})$ etc.
For any $f,g\in \i_\xi$, we have the 
Plancherel formula:
$
\lan f^\Diamond \, g\ran_\xi= (\hat{f},\hat{g}).
$
\end{theorem}
{\it Proof.}  By construction:
 $\si^{-1}(\hat{S}_{\hw})=
S_{\hw}$. This gives $FH=\si^{-1}(H)F$.
Then we use that the pairing
$\sum_{\hw\in \hW}f(\hw)\,T_{w_0}(g^\varsigma)(\hw)\,\mu(\hw)$ 
corresponds
to the anti-involution $\Diamond$ and that 
$\si^2(H)=T_{w_0} H^\varsigma \,T_{w_0}^{-1}$ for $H\in \HH$.
See Corollary 3.4.3 from \cite{C101}. This is a general fact
for any kind of $\hW$-invariant integration with the measure
function $\mu$; the Jackson integration 
$\sum_{\hw\in \hW}f(\hw)\mu(\hw)$ is taken here as such.
\sq
\vskip 0.2cm

Recall that the numerator of $\mu$ is nonzero at $q^{\xi}$ if
and only if the corresponding $\i_\xi$ is $Y$-semisimple with
simple spectrum; 
the denominator of $\mu(q^x)\mu(q^{-x})$
is nonzero at $x=\xi\,$ if and only if $\i_\xi$ is irreducible. 
Equivalently, $\i_\xi$ is irreducible if and only if 
all binomials in the numerators and denominators of (\ref{ediamu}) 
for all $\hw$ are nonzero.

The values  $\mu(\hw)$  are naturally some residues,
which will be used to obtain meromorphic continuations
of the integral formulas for the inner products. Thus, we
interpreted these values  as
norms of $Y$-eigenvectors $f_{\hw}\in\i_\xi$.
For $\xi=-\rho_k$,  the elements
$e_{\pi_b}, f_{\pi_b}$ become 
special normalizations of Macdonald's polynomials
 in $\mathscr{X}$, 
the quotient of $\i_\xi$.

\vskip 0.2cm

\section{\sc  Residues and closed subsystems}
We will use the definition of the residues 
from \cite{GH}, Ch.5. Generally,
 $Res_0\Bigl(\frac{h(x)}{f_1(x)\cdots f_n(x)}\,
d\,x_1\wedge\cdots\wedge
d\,x_n\Bigr)=
h(0)/\det\Bigl(\frac{\partial f_i}{\partial x_j}\Bigr)(0)$,
where the orientation  of the integration domain
$\{x=(x_i)\in \C^n\,\text{\large $\mid\,$}  |f_i(x)|<\ep\}$ 
is by the inequality
$d(\arg(f_1))\wedge\cdots\wedge
d(\arg(f_n)\ge 0$. The assumptions here are that
$h(x)$ is regular at $x=0$,
$f_i(0)=0$ and the determinant 
is nonzero, i.e. $0$ is a nondegenerate singularity.

 We will fix below 
the orientation 
to ensure that
\begin{equation}\label{resxi}
Res(\mu,0)=Res_0(\mu)=\prod_{\al>0}\prod_{i=1}^\infty
\frac{1-q_\al^{i+(\xi,\al^\vee)}}
{1-q_\al^{i-k_\al+(\xi,\al)}}\cdot 
\frac{1-q_\al^{i-(\xi,\al^\vee)}}
{1-q_\al^{i-k_\al-(\xi,\al)}}.
\end{equation}
We have  here $f_i\!=\!(1-t_iX_{\al_i})$.
If the variable $x_{\al_i}, 1\le i\le n$ are naturally ordered, then
the orientation is 
{\sf clockwise} for 
the loops around $f_i=0$. Permuting $\{x_{\al_i}\}$ will not
change the residue, because the orientation will change too.
The orientation and the corresponding wedge forms 
will be extended below to points $\hw$  using the
action of $\hW$. 

We note that the residue of any Laurent series in terms of $X_i$
is its constant term and it does not change if $X_i$ are
changed to variables $X'_i=\prod_{i=1}^n X_i^{c_{i,j}}$ for
$(c_{i,j})\in GL(n,\Z)$. We will use this below. However,
the presentation
of a function as a Laurent series
depends on the domain where the function 
is considered.

Generally, the residues can be complicated to calculate
algebraically. Analytically, they are integrals
of some top wedge forms $\om$ over 
$\Ga=\{x\in \C^n \text{\,\large $\mid$ } |f_i(x)|=\ep, 1\le i\le n\}$ 
and depend
only on the  (middle) homology class of $\Ga$ in $H_n\bigl(\{x 
\text{\large\, $\mid$ }
\prod_{i=1}^n f_i(x)\neq 0\}\bigr)$ and the class of the form $\om$
in the corresponding cohomology. See \cite{GH}. 
\vskip 0.2cm

{\bf Closed subsystems}. We will need {\sf closed root subsystems}
$R' \subset R$ (``closed subroot systems" is used too)
or those in $\tR$ of the same rank as $R$.
By definition, it is required that $\tal+\tbe\in
R'$ for any roots $\tal, \tbe$ in 
$R'$ 
if this sum belongs to $\tR$. Also, we will consider 
{\sf full} affine extensions $\tilde{R}'$ of $R'$,
which are with all $[\al,\nu_\al \Z]$ if $[\al,\cdots]\in R'$.
The positivity there will be induced from that $\tilde{R}$
unless stated otherwise. We will actually use the notation
$R'$ for subsystems in $R$; otherwise (in $\tR$), the notation
$R^\dag$ will be used. 

Let $\tR_{\lng}$ and $\tR_{\sht}$ be the root subsystems formed
by long and short roots in $\tR$ (similarly, for $R$); they
are of rank $n$.
The sum $\tal+\tbe\in \tR$ of 2 long roots $\tal$ and $\tbe$ 
can be only 
long, so $\tR_{\lng}$ is a closed root subsystem of rank $n$.
Indeed, $(\tal,\tbe)<0$ in this case; otherwise,  
 $|\tal+\tbe|^2/2>|\tal|^2/2=\nu_{\lng}$, which is impossible.
  Thus, $\tal+\tbe=
s_{\tal}(\tbe)$, i.e. it is long. Recall that
$(\tal,\tbe)=(\al,\be)$ for the non-affine components for $\tal,\tbe$.

Similarly,  $(\tal,\tbe)<0$ for short $\tal$ 
and long $\tbe$ if $\tal+\tbe\in \tR$. Thus,
$\tbe+\nu_{\lng}\tal=s_{\tal}(\tbe)$ is long and 
$\tbe+\tal=s_{\tbe}(\tal)$ is short in this case.
Similarly, $\tal+\tbe$ can
be a long root for short $\tal$ and $\tbe$ only if $(\tal,\tbe)=
(\al,\be)=0$ unless for $G_2$. In this case,
$s_{\tal}(\tal+\tbe)=
\tbe-\tal$ is a long root too. For $G_2$, $\tbe-\tal$ will be long
if $(\al,\be)=0$ for short $\tal$ and $\tbe$.
\vskip 0.2cm
\vfil

In the finite case, the list of closed {\sf maximal}
 subsystems $R'\subset R$  of rank $n$ is 
essentially
due to Borel- de Siebenthal;
there are no such subsystems for
$A_n$ and they are always reducible unless 
for $B_n, E_{7,8}, F_4, G_2$. Setting 
$\th=\sum_{i=1}^n n_i \al_i$, the key step is that any $\al_i$
with $n_i>1$ (assumed  prime for the maximality) can be
replaced by $-\th$ to generate such an $R'$, 
possibly reducible. 

We note that the usage of $\vth$ here instead of $\th$
leads to root subsystems
of rank $n$ in $R$,  but they can be non-closed. For instance, 
$B_m \oplus B_{n-m}\subset B_n$ for $2\le m\le n-2$  can occur
in this way. It is of rank $n$ but 
non-closed: $\vep_m+\vep_{m+1}$ in the standard notation
is a root, but
not in this subsystem. Here $\al_i=\vep_i-\vep_{i+1}$
for $i<n$, $\al_n=\vep_n$ and $\vth=\vep_1$.

The  Dynkin diagram of $B_n$ extended by $\al_0=[-\vth,1]$ is
that for the usual extended diagram
of $C_n$ with the reversed arrows. So the examples
above are when it splits into two connected segments.
\vfil

In the affine case, the description of {\sf maximal} 
closed subsystems $R^\dag\subset \tR$
is quite similar; see
Theorem 5.6 from \cite{FRT} and \cite{RV} for the maximal
ones. The affine classification is
basically the nonaffine one with the list of $p_\al>0$
such that affine roots $[\al,\nu_\al j]\in \tR'$ are those
for $\{j\}=
\{j_0+p_\al \Z\}$; such $p_\al$ always exist. 
If the {\sf maximal} closed ones in $\tR$ are known, then {\sf all} closed
root subsystems of the same rank as $R$ can be found 
by induction. Basically, the tables of maximal closed root
subsystems in $R$  are sufficient for this.
\vskip 0.2cm
\vfil

We will need below the affine root
subsystems $\tR'=\{[\al',\nu'_\al j],\al'\in R'\}$, where $\nu'_\al=
\nu_{\al'}$ is taken from $R$, 
and other affine definitions
for reduced ($=$decomposible) $R'$.
The corresponding $\vth'$ and $\al'_0$ are not unique then.
They must be defined for each connected
component of the Dynkin diagram of $R'$. The affine
Weyl group $\tW'$ becomes the direct
products of those from the connected
components; the corresponding $P$-lattice $P'$ 
and the extended affine Weyl group
$\hW'$ are the products of those for the connected components.

\begin{theorem}\label{thm:residues}
Let $0<q<1$ be generic: 
$q^m\neq t_\nu^l$ for any integer $l,m\neq 0$ and $\nu$.
Assume that the numerator of $\mu$
is nonzero at $q^{\xi}$, which condition does not depend
on the choice of the positivity in $\tR$. Equivalently,
$\i_\xi$
is semisimple with the simple spectrum.


(i) For $\xi=-\rho_k$, assume that $t_{\sht}^{\nu_{\lng}}\neq 1$;
also, let $t_{\sht}^j\neq -1$
for any $1\le j<n$ in the case of  $C_n$ and 
$t_{\sht}^4\neq 1$ for $F_4$. 
Then there are exactly $n$ binomials in the denominator of 
$\mu$ vanishing at $\hw=bw$, i.e. at $q^{b-w(\rho_k)}$, 
if and only if $\hw=\pi_b$ for $b\in P$, i.e.
when $w=u_b^{-1}$. Given $\pi_b$,  these binomials are
$\bigl(1-t_iX_{\pi_b(\al_i)}\bigr)$ for $1\le i\le n$.
For other $\hw$,
the number of such binomials vanishing at
 $\hw$ is smaller than $n$ and $\mu(\hw)/\mu(0)=0$; see
formula (\ref{muphoval}).

(ii) Let $t_{\sht}=t_{\lng}$ or 
$t_{\sht}^{\nu_{\lng}}=t_{\lng}.$
We continue to assume that  $\xi$ is such that  
the numerator of $\mu$ is
nonzero at $q^{\xi}$ and assume now that
its denominator has exactly $n$
binomials $(1-t X_{\tbe_i})$ 
that vanish at $q^{\xi}$. Let  $\tbe_i=[\be_i,\cdots]$. Then 
$\{\tbe_i,1\le i\le n\}$ is a set of simple roots
in the {\sf closed} root subsystem
$R^\dag=\tR\cap \oplus_{i=1}^n\Z\tbe_i$. Unless there exist
short $\tbe_i,\tbe_j$ for the systems $BCFG$ such that  
$\tbe_i-\tbe_j=[\be,m]$ for long $\be\in R$, where $m$ is not divisible
by $\nu_{\lng}$,  the set $\{\be_i,1\le i\le n\}$ is a set of 
simple roots in 
$R'=R\cap \oplus_{i=1}^n\Z\be_i$.

(iii) Continuing $(ii)$, 
let $\tilde{R}'=\{[\al,\nu_\al j] \mid \al\in R',j\in \Z\}
\subset \tR$, $|t_{\sht}|>1$  and $q$ is such that
$q<t_{\lng}^{-h_\dag}$, where $h_\dag$ is the maximum of 
Coxeter numbers
of the irreducible components of $R^\dag$. 
Then 
$\{(\tw')^{-1}(\tbe_i)\}$ become simple roots of $\tR'$ for
the positivity induced from $\tR_+$
and some $\tw'\in \tW'\subset \tW\subset \hW$, where $\tW'$ is
defined for $\tR'$.  More exactly, for every connected component
of $R'$, exactly one simple root $\al'_{i^\circ}$ 
for $i^\circ$ from the corresponding
{\sf twisted-affine} Dynkin diagram 
is {\sf not}  in $\{(\tw')^{-1}(\tbe_i)\}$.

(iv) For $\xi=-\rho_k$ as in (i), 
$Res(\mu,\pi_b)=\frac{\mu(\pi_b)}{\mu(0)} Res(\mu,0)$,
where the ratio is calculated in (\ref{muphoval}) and 
the residues are as above. Explicitly,
\begin{equation}\label{reszero}
Res(\mu,0)=\prod_{i=1}^n (1-t_i X_{\al_i})
\prod_{\tilde{\al}>0}\frac{1-X_{\tilde{\al}}}
{1-t_\al X_{\tilde{\al}}}.
\end{equation}
For $\xi$ in the setting of $(ii-iii)$, the formulas are
as follows. 
The corresponding residues must be calculated
for $\tR',\mu'$ as for $(i)$ and then
multiplied by
$\mu/\mu'(q^\xi)$, which is assumed nonzero.

\comment{
Res(\mu,0)=
\prod_{\al>0}\prod_{i=1}^\infty
\frac{1-q_\al^{i-(\rho_k,\al^\vee)}}
{1-q_\al^{i-k_\al-(\rho_k,\al)}}\cdot 
\frac{1-q_\al^{i+(\rho_k,\al^\vee)}}
{1-q_\al^{i-k_\al+(\rho_k,\al)}}.
\end{equation}
}
\end{theorem}
{\it Proof.}
{\sf (i)}.
The binomials $\al_i (1\le i\le n)$ are such for $\hw=0$, 
i.e. at the point $q^{-\rho_k}$.
Then   
$\bigl(1-t_iX_{\hw(\al_i)}\bigr)$ for $1\le i\le n$ 
belong to the denominator of $\mu$ if and only if
 $\hw(\al_i)\in \tR_+$ for $1\le i\le n$ 
and $\La(\hw)$ does not contain roots from $R_+$. This is the
defining property of elements $\pi_b$; see 
(\ref{La-pi-prim}). Thus $\hw=\pi_b$ for some $b\in P$ and
the ratio $\mu(\hw)/\mu(0)$ is then
nonzero due to (\ref{muphoval}). Thus, it suffices to
check that the binomials from the denominator of $\mu$
vanishing at $0$ are exactly those for $\{\al_i,1\le i \le n\}$.

Next, if $\bigl(1-t_\al X_{\tal}\bigr)(q^{-\rho_k})=0$ 
for $\tal=[\al,\nu_\al j]>0$, then 
$j=0$ because $q$ is generic, i.e. $\tal$ is nonaffine
and $\al>0$; let $\nu=\nu_\al$. 

For this $\al$ and any $i>0$ such that $\nu_i=\nu$, one has:
$(\al, \al_i)\le 0$.
Otherwise, there exists $\al_i$ such that $\be=
\al-\al_i$ is a positive 
root in $R$ satisfying 
$\bigl(1-X_{\be}\bigr)(q^{-\rho_k})=0$. However, the assumption
is that this is impossible for any $\be$ (positive or negative).
These inequalities give
that  $\al$ and the roots $\al_i$ such that  $\nu_i=\nu$ are linearly
independent, which is impossible in the case of $A,D,E$. 

Let us
consider now  $B,C,F,G$. Then such $\{al_i\}$ are simple roots 
in the root subsystem
$R_{\nu}$ formed by all roots $\be\in R$ such that $\nu_\be=\nu$,
but possibly not all simple roots there. 
The positivity in $R_\nu$ is that induced from $R$;
$\al$ remains positive in $R_\nu$. Let us check that $\al$
is simple in $R_{\nu}$. We will use the notation from 
the tables of \cite{Bo}.



For any 
non-simple positive {\sf short}
root $\al$, there exists $\al_i$ of
the same length such 
that $\be=\al-\al_i\in R$. This gives
that $X_\be(q^{-\rho_k})=1$, which 
contradicts our condition
for the numerator of $\mu$. We conclude that  $\al$
can be only long if it is non-simple.

The same claim (the existence of $\al_i$) holds
for long $\al$ in  $R_{\lng}$ unless
$\al$ is simple in $R_{\lng}$ with one reservation.
In the case of $C_n$, there is no such $\al_i$ for 
$\al=\vep_j$ for  $j<n$ in the 
notation from \cite{Bo}. For such $\al$, 
$\be=(\al-\al_n)/2=\vep_j-\vep_n$ is a short root in $R$ 
and  $X_\be(q^{-\rho_k})=\pm 1=t_{\sht}^{n-j}$. The latter
relation is  excluded and we can omit $C_n$ in the next
considerations.
\vfil

Let us consider now $B,G$. Since long $\al$ is linearly independent
with $\al_i$, the dimension of the space generated by $\al_i$ 
is $(n-1)$. Thus $\al$ must be the unique
simple root of $R_{\lng}$ that is not one of $\al_i$. Recall that
$\{\al_i\}$ are simple in $R$ and remain simple in $R_{\lng}$,
but the latter system  contains other simple roots (unless
for $ADE$). We obtain that 
$(\al-\al_m)/\nu_{\lng}$ is a short root in $R$ for some $\al_i$,
which contradicts the condition  
$t_{\sht}^{\nu_{\lng}}\neq 1$.  
 
\vfil
 
So the simplicity of long $\al$ remains to be checked 
only in $R_{\lng}$ for $F_4$. Then $(R_{\lng})_+=\{\vep_i-\vep_j\}$
for $i<j$ and $\al_1=\vep_2-\vep_3, \al_2=\vep_3-\vep_4$ for $F_4$
in the notation from  
Plate 8 of \cite{Bo}. Then for any non-simple root $\be>0$
in $(R_{\lng})_+$ either $\be-\al_i$ or $(\be-\al_i)/2$ belongs
to $R$, where $i=1$ or $i=2$. We come to a contradiction.
\vfil

Thus, we obtain that $\al$ must be 
simple in $R_{\nu}$ (but not simple 
in the whole $R$). One has: $\al=\al_{m}+\nu_{\lng}\al_{l}$ and
$\al=\vep_1-\vep_2$ unless possibly for $F_4$.
Thus, $t_{\sht}^\nu=1$ or 
$\al=2\al_4+\al_2+2\al_3$ for $F_4$,  which results in
$t_{\sht}^4=1$. These two relations were excluded. They can really
occur as well as the relation $t_{\sht}^{n-j}=\pm 1$ in the
case of $C_n$. 
\vfil

We note that, actually, the
classification is not strictly necessary for the
last step. One can use that there exists at least one short
$\al_i$ such that $(\al,\al_j)>0$. Indeed, the rank
of $R$ would be $>n$ otherwise. Thus, $\al-\nu \al_i$ is a long 
positive
root. It can be only simple in $R$ if $\al_i$
is neighboring to long simple roots in the Dynkin diagram.
The case of $\al_i=\al_4$
for $F_4$ is exceptional and must be considered separately.  
\vskip 0.2cm

{\sf (ii).} Similarly to the considerations above,
one has: $(\be_i,\be_j)\le 0$ for $1\le i<j\le n$.
Indeed, $\tbe_i-\tbe_j\in \tR$ otherwise
and, additionally,  $\tbe_i-\nu_{\lng}\tbe_j\in \tR$
if $\nu_j<\nu_i$. This gives that one of these differences 
will make the numerator of $\mu$ vanishing, which is impossible.
Then the required claims result from the following lemma. 

\begin{lemma} \label{lem-dif}
Let $(\be_i,\be_j)\le 0$ for $\tbe_i=[\be_i,\cdots]\in \tR$
and $1\le i<j\le n$. Then $\be_i$ can assumed in $R_+$ upon the
action of some $w\in W$. Provided this, assume that  
$\be=
\sum_{i=1}^n m_i\be_i\in R$ with $m_i\in \Z$ such that
$m_im_j<0$ for at least 
one pair $(i,j)$. Then there exist $i,j$
such that $\be_i-\be_j\in R$ and, additionally, 
$\be_i-\nu_{\lng} \be_j\in R$ if $\nu_j<\nu_i$. Moreover, 
$\tbe_i-\nu_{\lng} \tbe_j\in \tR$, including the
$ADE$ systems. For $BCFG$, $\tbe_i-\tbe_j\in \tR$ unless
$\be_i,\be_j$ are short and  $\be_i-\be_j$ is long or (always)
if $\tbe=\sum_{i=1}^n m_i\tbe_i\in \tR$. 
\end{lemma}
{\it Proof.} The positivity condition making
$\{\be_i\}$ positive is $(\eta,\be)>0$ for 
$\eta=-\sum_{i=1}^n \eta_i\be_i$ for
sufficiently general $\eta_i>0$ and they are linearly independent,
which is standard. 
We set $x=\be-\sum_{m_j<0}m_j\be_j=\sum_{m_i>0}m_i\be_i$. Then 
$(x,x)>0$ and $(\be,\sum_{m_i>0} m_i \be_i)>0$. Therefore,
$(\be,\be_i)>0$ for at least one $\be_i$ with $m_i>0$, and
$\be'=\be-\be_i\in \tR$ has the corresponding  
$\sum_{m'_i>0} m'_i$ smaller by $1$ than that for
$\be$. Similarly, we can diminish $-\sum_{m_j<0}m_j$ by $1$ and 
continue diminishing the sums $\sum_i$ or $\sum_j$ 
until we obtain $\be_i-\be_j\in R$. If $\nu_j<\nu_i$
here, then $(\be_i,\be_j)>0$; otherwise, $|\be_i-\be_j|>|\be_i|$.
This results in 
$\be_i-\nu_{\lng}\be_j=s_{\be_j}(\be_i)\in R$. Moreover, then
$\tbe_i-\nu_{\lng}\tbe_j-s_{\tbe_j}(\tbe_i)\in \tR$,

The argument above used for the relation
$\tbe=\sum_{i=1}^n m_i\tbe_i\in \tR$
provides (formally) that $\tbe_i-\tbe_j\in \tR$.
Generally, $\be_i-\be_j$ does not imply 
 $\tbe_i-\tbe_j\in \tR$  only if
$\be_i,\be_j$ are short, $\tbe_i-\tbe_j=[\be, m]$ is long and
$m$ is not divisible by $\nu_{\lng}$. This proves 
the last claim.
\sq 
\vskip 0.2cm

{\sf (iii).} The roots $\tbe_i$ are positive with respect
to the following positivity condition $|X_{\tbe}(q^\xi)|<1$ 
due
to the inequalities $|t_\nu|>1$. Recall that $X_{\tbe}=
X_\be q^{\nu_{\be}j}$
for $\tbe=[\be,\nu_\be j]$ and 
$X_{\tbe}(q^\xi)=t_\be^{-1}$
for $\tbe =\tbe_i$.
It is possible that
$|X_{\tal}|=1$ for some $\tal\in \tR$, so we may need
to deform $\xi$ a little to ensure that this is really 
some {\sf positivity} in $\tR$. For $\tR'$, it suffices 
to assume that  $q$ is sufficiently
small, which will be checked together with the simplicity of
$\tbe_i$.

For $\tbe\in \tR'$, the range of the values  $|X_{\tbe}(q^\xi)|$
is a union of $V_j=\{q^j |t_{\sht}|^m\}$,  where $j\in \Z$ and
$-C<m<C$ for some constant $C$ calculated
in terms of the Coxeter numbers of the irreducible
components of $R^\dag$. 
One can assume that $V_i\cap V_j=\emptyset$
for $i\neq j$ for sufficiently small $q$. Then
 $V_0=\{\,|X_{\tbe}(q^\xi)\,|\text{ s.t. }\, \tbe\in R^\dag\}$ 
by construction.

If $|X_{\tbe}(q^\xi)|<1$ is not a {\sf positivity condition}
 for the root system
$\tR'$ or if $\{\tbe_i\}$ are not simple for this positivity, then there
exists $\tbe\in \tR'$ such that $|X_{\tbe}(q^\xi)|\in V_0$, which 
can be only if $\tbe\in R^\dag$. 
Using Lemma \ref{lem-dif},
we obtain that then there exists
$\tbe=\tbe_i-\tbe_j\in \tR'$ such that 
$X_{\tbe}(q^\xi)=q^{\pm m}$ for $m>0$. However, this
is impossible for sufficiently small $q$. A more exact analysis
shows that the inequality for $q$ from $(iii)$ is sufficient here.
Alternatively,
one can use 
$(ii)$, which states that $\tbe_i$ are simple in $R^\dag$ for
some {\sf positivity}. 

Then we find $\tw'$ in the affine Weyl group $\tW'\subset \tW$ 
of $\tR'$ 
transforming the standard {\sf affine Weyl chamber} for $\tR'$ to
that for the positivity above. 
The roots $(\tw')^{-1}(\tbe_i)$ then become
those described in $(iii)$. 
\vskip 0.2cm

Part $(iv)$ is actually a reformulation
of $(i)$. We set $f_i=(1-t_i X_{\al_i})$
and choose the orientation clockwise. The corresponding
residue is
obtained from $\mu$ by the deletion of these binomials from the
denominator and the evaluation of
the rest at $0$, which is the point $q^{-\rho_k}$. The extension
to the setting of $(ii-iii)$ is straightforward.\sq
\vskip 0.2cm

{\bf Comments.}
Given $t_\nu$, the inequality for $q$ in $(iii)$  means that 
$-\Re k_\nu$ must be sufficiently large. Recall that $q\to 0$
is the limit to AHA. Actually, 
the condition for $q$ needed here is entirely algebraic.
This inequality provides it, but this claim holds for generic $q$,
which is similar to Lemma \ref{lem-dif}. 

{\sf Residues.} Let us provide a variant
of  formula (\ref{reszero}) in $(iv)$ for $b=0$ 
when $\{\tbe_i\}=\{\al_i, 0\le i\le n\}\setminus
\{\al_j\}$ for some $j\ge 0$. I.e. the formula below
will be its (minor) generalization. 
Then  
$\rho^{\dag}_k=
\frac{1}{2}\sum_{\tal\in R^{\dag}_+} \nu_\al \tal$\  and 
$$
Res(\mu,0)=(1\!-\!t_j X_{\al_j})\,
\prod_{i=0}^n(1\!-\!X_{\al_i})\!\!\!\!
\prod_{\tilde{\al}>0,\tal\neq \al_i}\frac{1\!-\!X_{\tilde{\al}}}
{1\!-\!t_\al X_{\tilde{\al}}}\,
\bigl(X\!=\!q^{-\rho^{\dag}_k}\bigr),
$$
which is for a suitable choice of the orientation.
\comment{
\begin{equation*}
Res(\mu,0)=\frac{1-q_\vth^{1-k_\vth+(\rho^{\dag}_k,\vth)}}
{1-q_{j}^{-k_{j}-(\rho^{\dag}_k,\al_j)}}
\prod_{\al>0}\prod_{i=1}^\infty
\frac{1-q_\al^{i-(\rho^{\dag}_k,\al^\vee)}}
{1-q_\al^{i-k_\al-(\rho^{\dag}_k,\al)}}
\frac{1-q_\al^{i+(\rho^{\dag}_k,\al^\vee)}}
{1-q_\al^{i-k_\al+(\rho^{\dag}_k,\al)}}.
\end{equation*}
}
The extension to arbitrary $\xi$ when the
numerator of $\mu$ has no zeros
and the corresponding $R'$,
 including the case of different $k_\nu$,
is quite similar. 
\vskip 0.2cm

\vskip 0.2cm

{\sf Parts (ii-iii).} The conditions $t_{\lng}=t_{\sht}$, or 
$t_{\lng}=t_{\lng}^{\nu_{\lng}}$ there are the two cases of 
{\sf equal parameters} in the twisted setting.
Actually, the latter relation is more common; for instance,
it is compatible with the usage of DAHA for quantum group
invariants of links. 
We obtained that the classification of $\xi$ under these conditions
can be reduced to that of  closed {\sf finite} root
subsystem $R^\dag$ of rank $n$ in a closed {\sf affine}
subsystem $\tR'\subset \tR$, for a
closed root subsystems $R'\subset R$ of rank $n$. 
The classification of the
latter up to the action of $W$ follows from the
Borel - de Siebenthal theory. We note that this theory actually
uses ``affine tools", so passage from AHA to DAHA seems natural
from the perspective of classification the residual points. 

Generally, the simplest case is $R'=R$ when 
$\{\tw^{-1}(\tbe_i)\}=\{\al_i,i\neq i^\circ\}$ for some $\tw\in \tW$.
 Furthermore,
If $i^\circ\neq 0$, then we can assume that 
$\nu_{i^\circ} n_{i^\circ}>1$ modulo the action of $\hW$,
where $\th=\sum_{i=0}^{n} n_i\al_i$.
If $R'=R$ and $i^\circ=0$, 
then we arrive at $(i)$: 
$\{\tw^{-1}(\tbe_i)\}=\{\al_i,1\le i\le n\}$.
\vskip 0.2cm

{\sf The case of $A_n$.} This always holds for
$A_n$ because the only closed root subsystem in
$R$ of rank $n$ is $R$ and all $n_i$ are $1$.
Thus, we can take
$i^\circ=0$ for $A_n$ modulo $\hW$ and 
$\xi$ from $(ii-iii)$ are
$-k\rho$ and their images under the action of 
$\pi_b$ for $b\in P$.
This is parallel  to the ``orbit" of the {\sf Steinberg
representation} in the AHA theory.
\vskip 0.2cm

{\sf Induced modules.}
Recall that the numerator of $\mu(q^{\xi})$ is nonzero if and
only if the $\HH$-module $\i_{\xi}$ is $Y$-semisimple with 
simple spectrum.
Generally, $(ii)-(iii)$ give some class of $\xi$ for generic $q,t$
where 
$\i_\xi$  are direct counterparts
of $\i_{-\rho_k}$. Their canonical irreducible
 $\HH$-quotients described in Theorem 3.6.1 from \cite{C101}
generalize $\mathscr{X}$. In the notation there: $\Upsilon_0=\hW$
and $\Upsilon_*=\Upsilon_+$. 
We note  some links to \cite{VV}.

\vskip 0.2cm
\vfil

{\bf Some examples.}
The closed subsystems $R'\not\subset R$ of 
rank $n$ from $(ii)$ can
be ``even" $A_1^n=A_1\oplus \cdots \oplus A_1$
(n times) for $C_n,D_{2m\ge 4},E_{7,8},F_4,G_2$, the most
reducible.  We use the
notation 
$X\oplus Y$ for the root system $X\cup Y$ in the direct sum
of the corresponding $\R$-spaces. 
For instance,  
$R'=\{\be_i=2\vep_i, 1\le i \le n\}$ is such for $C_n$. 
In this case, we must have
$X_{\al_i}(q^\xi)=X_{\vep_i-\vep_{i+1}}(q^\xi)=-1$ for $1\le i<n$
to ensure that the numerator of the corresponding
$\mu$ is nonzero at $q^\xi$.

For $D_{4}$, the closed subsystem $R'=A_1^4$ is
as follows: $\be_{1,2}=\vep_1\pm \vep_2, \be_{3,4}=\vep_3\pm \vep_4$. 
Accordingly, $\xi=-{\rho'}_k=
-k(\vep_1+\vep_3)$. One has:  
 $X_{\be_i}(q^\xi)=t^{-1}$ for $1\le i\le 4$  
and 
$X_{\vep_2-\vep_3}(q^\xi)=-t$; notice the minus-sign.
Another variant is for
$\be_{1}=\vep_1\pm \vep_2, \be_2=\vep_2-\vep_3,
\be_{3,4}=\vep_3\pm \vep_4$. i.e. for $A_1\oplus D_3$ in
$D_4$ ($D_3=A_3$); then $X_{\vep_i\pm \vep_j}=
t^{(4-i)\pm (4-j)}$
can be taken. These two examples can be readily extended to 
$D_{n-m}\oplus D_m$ in $D_n$ for any $1\le m\le n, n\ge 4$. 
As above, 
the notation is from \cite{Bo}. 

Not all closed
subsystems of rank $n$ can really occur in $(ii)$;
say, $A_1^6$ in $D_6$ will have zeros in the numerator
of $\mu$ if we follow the above construction for 
$A_1^4\subset D_4$.
\vskip 0.2cm

Let us give an example when not all $\be\in R'$ can be
lifted to $R^\dag$ and $\be_i$ are not all simple in $R'$.
For the root system $B_n$, we take
$\tbe_i=\al_i=\vep_i-\vep_{i+1}$ for $1\le i\le n-2$, $\tbe_{n-1}=
\vep_{n-1}, \tbe_n=[-\vep_n,1]$. 
Then
 $X_{\vep_n}(q^\xi))=q\, t_{\sht}$,   $X_{\vep_{n-1}}(q^\xi))=
t_{\sht}^{-1}$,  $X_{\vep_{n-i}}(q^\xi))=
t_{\sht}^{-1}t_{\lng}^{1-i}$ for $2\le i \le n-1.$ Thus,
$X_{\vep_n+\vep_{n-1}}(q^\xi)=q$ and all other $X_{\al}(q^\xi)$
for $\al\in R$ contain powers of $t_\nu$. For
$|t_\nu|>1$ and generic $q$,
the numerator of $\mu(q^{\xi})$ is nonzero.
It is used here 
that $\tbe_{n-1}-\tbe_n=[\vep_{n-1}\!+\!\vep_{n-2},-1]$ is
not from $\tR$ because  $\vep_{n-1}\!+\!\vep_{n-2}$ is long.
\vskip 0.2cm

\comment{
The corresponding
$\xi$ form ``reduced packets": the orbits of $\xi$ 
with the respect to the action of $W$ on $R'$ and
$\hW'$ on $R^\dag$.  The representatives 
of $R^\dag$ are described in $(ii)$. 
}
\vskip 0.2cm

\section{\sc Residual subtori and points} 
Informally, they are those that can {\sf potentially} occur
in the meromorphic 
continuation of the functional 
$I^{im}(f)=
\int _{\imath\R^n} f(x)\mu(q^x; q,t)dx$ from $\Re k_\nu>0$ to all
complex $k_\nu$ or for $I^{\imath\aa}(f)$. 
 If they can be obtained from each other
by the action of  $\hw\in \hW$, we say that they are in the same
{\sf packet}. However not all $\mu$-residual subtori and points 
defined below really occur
in the integral formulas; finding them  is a combinatorial problem,
which can be involved. After they are found,
the count of the 
corresponding coefficients, the residues for the points, is an 
entirely algebraic procedure. The following definition
is a double affine
extension of Definition 2.1 
 from \cite{HO} coupled with Theorem 2.2  to the 
$\mu$-function.  Also, see
\cite{O20} (Theorem 7.1, Remark 7.3).

\begin{definition}\label{DEF:restori}
We continue to assume that $0\!\le\! q\!<\!1$ and $t_\nu$ are
sufficiently general. The {\sf double affine residual
subtori}, called $\mu$-residual below,  
are the affine tori $T$ of codimension $m$
given by the equations $1\!-\!t_{\tbe_i} X_{\tbe_i}\!=\!0$
for $1\le i\le m$ and linearly independent $\tbe_i\in \tR_+$,
provided the following condition.
The number $\kapp_1$  of the binomials $(1\!-\!t_{\tal} X_{\tal})$ 
for $\tal\in \tR_+$
vanishing at $T$ must be $\ge \kapp_0+m$ for the number $\kapp_0$
of $(1\!-\!X_{\tal})$ for $\tal\in \tR_+$ vanishing at $T$. 
The {\sf $\mu$-residual points} are for $m=n$.\sq
\end{definition}

The residual points play a key role in the $q,t$-case.
They alone are sufficient to obtain the meromorphic continuation
of $I^{im,\imath\aa}(f)$ for $|t_\nu|>1$ ($\Re k_\nu<0$) provided
the integrability and the convergence of
$f(x)$. The convergence conditions depend on $\Re k_\nu$
and the order of iterated integrations. Any
analytic functions $f(x)$ integrable in the
imaginary directions of no greater than exponential
growth in the real directions can be taken here when 
$\Re k<0$ is sufficiently large.
We will provide a reasonably complete general
description of residual {\sf points} for sufficiently 
general $t_\nu$ in the case of ``equal parameters".
 The calculation of the corresponding residues
is straightforward when $\kapp_1-\kapp_0=n$.

We note that a direct  
affine generalization of the AHA residual subtori from
in \cite{HO} 
is more restrictive. In our context,
it would be  $\kapp_1-\tilde{\kapp}_0\ge m$, 
where $\tilde{\kapp}_{0}$
is the number of binomials $(1-X_{\tal})$ vanishing at
$0$ for $\tal\in R_-\cup \tR_+$. This is basically the switch
to $\de$, the symmetrization of $\mu$,
and  $W$-invariant functions $f(x)$;
we will not discuss this possibility in the paper.
\vskip 0.2cm

Following the proof of part $(ii)$ of
Theorem \ref{thm:residues}, we
obtain the following claim, which reduces the
description of $\mu$-residual points 
to some combinatorial analysis of the corresponding
root system. Any residual $\mu$-point $\xi$ can be obtained
by the following construction, though we do not claim that
they occur in some integral formulas and that the
corresponding residues are nonzero. 

\begin{theorem}\label{thm:resid}
As in $(ii,iii)$ of Theorem \ref{thm:residues}, 
  $0<q<1$, let 
$q^m\neq t_\nu^l$ for any integer $l,m\neq 0$,
$|t_\nu|>1$ for any $\nu$ and either $t_{\lng}=t_{\sht}$ or
$t_{\lng}=t_{\sht}^{\nu_{\lng}}$. Also, we assume
that $q$ is sufficiently  general by imposing the condition
from $(iii)$ there.

Given a
{\sf closed root subsystem}  $R^\flat\subset R$  of rank $n$,
we begin with a subset  $\{\al^\flat_i, i\in I\}\subset 
\{\al^\flat_i\}$ 
of simple roots of $\tR^\flat_+=\tR_+\cap \tR^\flat$ such that
exactly one simple root is
removed from  $\{\al^\flat_i\}$ for every connected 
component of $R^\flat$. We follow Theorem \ref{thm:residues}, (iii).
Then we fix $\hw\in \hW$. Let $\tbe_i=\hw(\al^\flat_i)$ and
$R^\dag$ be the (closed) root subsystem with simple roots 
$\{\tbe_i=[\be_i,\cdots] \text{ for }
i\in I\}$, which are assumed from  $\tR_+$. 

Next, let $\{i\in I'\}$  be a subdiagram of the
Dynkin diagram $\{i\in I\}$ of $R^\dag$ and $R^\ddag$ be
the corresponding closed root subsystem of $R^\dag$.
Then we define  $\xi\in \C^n$ such that 
$q^{(\xi,\tbe_i)}=1$ for $i\in I'$
and $q^{(\xi,\tbe_i)}=t_{\nu^\bullet_i}^{-1}$ for 
$i\in I\setminus I'$,
 where $\nu^\bullet_i=\nu_{\tbe_i}$,\, 
$(\xi,[\be,j\nu_\be])
=(\xi,\be)+j\nu_\be.$

Then the numerical condition for 
$\mu$-residual points $\xi$ becomes
$$
|\,\{\tbe=\tbe_m+\sum_{i\in I'} c_i\tbe_i\subset R^\dag_+ \text{ s.t. }
m\in I\setminus I',\, c_i\in \Z_+\}\,|-|\,R^\ddag_+\,|\ge n.
$$
Any $\mu$-residual points
occur in this way for proper 
$R^\flat,\hw,
R^\dag, R^\ddag$. Moreover, 
$\hw$ can be assumed from $\tW'$ if 
$\tR^\flat=\tR'$ for the closed root subsystem $R'$ generated by
$\be_i$ for $i\in I$. 
\end{theorem}
{\it Proof.}
The direct statement follows from the definition of
$\mu$-residual points. We need to check that any 
$\mu$-residual $\xi$
can be represented in this way.
Let $R'\subset R$
be a {\sf closed root subsystem}  of the same
rank as $R$ such that its standard (full) affine extension
$\tilde{R}'\subset \tR$ contains the set
$\r^1=\{\tbe\in \tR_+ \mid
X_{\tbe}(q^\xi)=t^{-1}_\be\}$ and the closed root subsystem
$\r^0=\{\tbe\in \tR_+ \mid
X_{\tbe}(q^\xi)=1\}$. 

We take simple roots of $\r^0$ and
add to them {\sf primitive roots} from $\r^1$ defined
as $\tbe$ there such that $\tbe\neq \tbe'+\tal$ for
$\tbe'\in \r^1$ and $\tal\in \r^0_+$. Let this set be
$\{\tbe_i, i\in I\}$, where $\tbe_i$ for $i\in I'$ are all simple
roots from $\r^0_+$.
This set linearly generates $\R^n$ and
satisfies the conditions 
$(\tbe_j,\tbe_i)\le 0$ for $i\in I, j\in I\setminus I'$ because
$\tbe_j$ 
are assumed primitive. Similarly, $(\tbe_i,\tbe_j)\le 0$ for primitive
ones, i.e. for 
$i,j\in I\setminus I'$, because
otherwise $\tbe_i-\tbe_{j}\in \r^0$ and one of them cannot be
{\sf primitive}. Thus, $\{\tbe_i, i\in I\}$
are linearly independent and $|I|=n.$ 


Then we impose the
inequality for $q$ from Theorem \ref{thm:residues}, $(iii)$
for the system $R^\dag$. 
Following the reasoning there, we introduce the positivity
condition for $\tbe\in \tR'$ by 
$|X_{\tbe}(q^{\xi'})|<1$, 
 where $\xi'$ is a small deformation such that 
$|X_{\tbe}(q^{\xi'})|\neq 1$ for $\tbe\in \r^0$.
One has: $X_{\tbe_i}(q^{\xi'})\approx t_{\nu_i^\bullet}^{-1}$,
where $i\in I\setminus I', \nu_i^\bullet=\nu_{\tbe_i},$
so they are still positive. This positivity may 
result in different simple roots in $\r^0$: let us take them
as $\tbe_i$ for $i\in I$ 
instead of the initial ones. 

Following the proof of $(iii)$,
$\tbe_i\, 
(i\in I)$ are simple roots for the positivity condition
above for sufficiently
small $q$ (under the inequality we imposed). 
 Thus, they become simple in $\tR'$ upon the action
of some $\tw'\in \tW'$ for the positivity condition
there induced from that 
in $\tR_+$. We obtain that
 $R^\flat$ is $W$-conjugated to $R'$. \sq

\section{\sc Integral presentations} For $\Re k_\nu>0$
the following inner products
in $\mathscr{X}$ induce $\Diamond_{\pm l}$ for $l>0$: 
\begin{align}\label{iminner}
&\lan f,g\ran_l^{im}=
\int_{\imath \R^n} f T_{w_0}(g^\varsigma)q^{-lx^2/2}
\mu(q^x;q,t)dx \text{ \  induces } \Diamond_l,\\
&\lan f,g\ran _{-l}^{re}=
\int_{\ \R^n} f T_{w_0}(g^\varsigma)\,q^{lx^2/2}
\,\mu(q^x;q,t) dx\, \text{ \  induces } \Diamond_{-l}.
\label{reinner}
\end{align}
Here $f,g\in \mathscr{X}$, but this can be extended
to any completions of $\mathscr{X}$ provided the 
analyticity of $f,g$ and the integrability.
We use here that $q^{lx^2/2}$ is $W$-invariant,
commutes with $T_{w_0}$ (considered as an
operator of multiplication),
and is preserved 
by $\Diamond$; see (\ref{diamondef}).
One has:\,

$
\lan f,H(g)\,\ran_l^{im}\,=\,\lan f,q^{-lx^2/2}H(g)\,\ran_0^{im}\,=\,
\lan\, \Diamond\bigl(q^{-lx^2/2}H\bigr)(f),\,g\,\ran_0^{im}$

$\ \ \ \ =\,\lan\, \bigl(q^{lx^2/2}\circ
\Diamond\circ q^{-lx^2/2}\bigr)(H)\bigr)f,\,g\,\ran_l^{im}$
for $H\in \HH$. 

\noindent
Here
$q^{lx^2/2}\circ
\Diamond\circ q^{-lx^2/2}(H)=q^{lx^2/2}\circ
\Diamond(H)\circ q^{-lx^2/2}=\tau^l_+\bigl(\Diamond(H)\bigr).
$

For $l=0$, the following integral
replaces (\ref{iminner}): 
\begin{align}\label{int-zero}
&\lan f,g\ran_0^{\imath\aa}=\frac{1}{(2\pi \imath \aa)^n}
\int_{\imath \R^n/2\pi\imath \aa P^\vee} f T_{w_0}(g^\varsigma)
\mu(q^x;q,t) dx=\\
\frac{1}{(2\pi \imath \aa)^n}&
\int_{-\imath \pi \aa}^{\imath \pi \aa}\!\!\cdots\!
\int_{-\imath \pi \aa}^{\imath \pi \aa}
f T_{w_0}(g^\varsigma)
\mu(q^x;q,t)\ dx_{\al_1}\!\cdots dx_{\al_n} 
\text{ for } q=e^{-1/\aa}.\notag
\end{align} 

Here the order of integration can be arbitrary, though the
meromorphic continuation to negative $\Re k_\nu$ depends
on this order. 
This integral 
coincides with the {\sf constant term} ct$( f T_{w_0}(g^\varsigma))
\mu$ for $\Re k_\al>0$ and
provided  the inequalities $|t_\al|^2<q_\al$. Indeed, $\mu(q^x)$ is 
analytic in the
annulus $t_\al q_\al^{-1}< |X_\al|<t_\al^{-1}$ for $\al\in R_+$.
Therefore
we can replace $\mu$ with the corresponding Laurent series:
its expansion in terms
of $q^i$ for $i\ge 0$ and $ t_\al X_{\tal}$ for $\tal\in \tR_+$. 

The fact that the imaginary integrals give the 
$\Diamond_l$-invariant
DAHA inner products for $l\ge 0$ does require the conditions 
$\Re k_\nu >0$. The give that there are
no singularities of $\mu$ between the initial contour of
integration and its translations by $b\in P$
when $\Re k_\nu$ is sufficiently large. Then the
analytic continuation to any $\Re k_\nu >0$ is used. 

Indeed, the poles of $\mu$ modulo the imaginary periods
are at
$x_{\al^\vee}=-k_{\al}-i$ and $x_{\al^\vee}=k_\al+i+1$
for $\al\in R_+, i\ge 0$. 
Thus, the ``gap" between $-\Re k_\al-1$ and $\Re k_\al+\nu$ 
gives the required  when $\Re k_\nu\gg 0$;
 Stokes' theorem is used.  
The integrals over $\imath \R^n$ 
make sense of course for any sufficiently general $k_\nu$ 
but the corresponding pairings are only
$\h_X$-invariant (not $\HH$-invariant for $\Re k_\nu<0$).  
\vskip 0.2cm
\vfil

{\bf Comments.} 
Making $g=1$, $\lan f, 1\ran^{im}_l$ is a coinvariant of level $l$,
i.e. that for $\Diamond_l$. For $f\in \mathscr{X}$, 
one can switch here from the imaginary integration
to $\lan\cdots\ran^{\imath\aa}$. Namely, we  
replace $q^{-lx^2/2}$ in the integrand 
with the sum of its translations 
by $2\pi\imath\aa P^\vee$ and 
use that $\mu$ is in terms of $X_\al$.

Alternatively, let $\lan f, g\ran^{\imath\aa}_l\equal
\lan f, g\, \Th(q^x)^l\,\ran^{\imath\aa}_0$,
where we can use the theta-function  
$\Th(q^x)\equal\sum_{b\in P}\,X_{b}\,q^{b^{2}/2}$
for $\tR$ instead of $q^{-x^2/2}$ because 
$q^{x^2/2}\Th(q^x)$ is $\hW$-invariant.
Then 
$\lan f, 1\ran^{\imath\aa}_l=
\lan f, \Th(q^x)^l\ran^{\imath\aa}_0$ is a coinvariant of
level $l$ too. 
Note that for $\Re k_\nu>1$,  the 
integral $\lan f, 1\ran^{\imath\aa}_l$ for $f\in \mathscr{X}$
is reduced to taking the corresponding constant term.

Generally, we have two different approaches,
which result in the coinciding  (proportional) formulas
only for $l=1$. This is because the space 
of coinvariants is one-dimensional for $\Diamond_{l}$
only for $l=0,\pm1$. 
 For $l=1$ the explicit connection is established
via the functional equation
for $\Th(q^x)$; see, e.g., \cite{Kac} and
 Lemma 4.6 from \cite{ChW}. Actually, this  is how the
functional equation for $\Th$ can be  justified. 
\vskip 0.2cm

For the sake of completeness, let us state 
Theorem 4.9 from \cite{ChD} in this context.
One has for $l=1$:
\begin{align}\label{evinner}
&\lan f,g\ran _1^{im}=\lan 1,1\ran _1^{im}\bigl(\tau_-^{-1}(f)
T_{w_0}(\tau_-^{-1}(g^\varsigma))\bigr)(q^{-\rho_k}). 
\end{align}
Here we use that $\tau_-$ acts in $\mathscr{X}$; 
the nonsymmetric Macdonald polynomials are its eigenvectors. 
As above: $X_a(q^b)=q^{(a,b)}$ and for any functions here.
In particular, $\lan f,1\ran _1^{im}=
\lan 1,1\ran _1^{im}\,t^{\frac{l(w_0)}{2}}
\bigl(\tau_-^{-1}(f)\bigr)(q^{-\rho_k}). $
\vskip 0.2cm
\vfil

{\sf The space of coinvariants.}
More generally, let us consider 
$\lan f, \Th[l]\ran^{\imath\aa}_0$
for any theta-functions $\Th[l]$ of level $l$. They are
by definition are analytic in terms of $q^x$ such that
$\Th[l]/\Th^l$ are $\hW$-invariant. 
These functionals
are {\sf coinvariants} of level $l$ for $\Diamond_l$. 
This approach actually gives
that the dimension of the space of such coinvariants
coincides with the number of
the integrable irreducible Kac-Moody modules of level $l>0$
for the root system $\tR$. This is an algebraic fact: 
Theorem 2.13 from \cite{CMa}. The proof there was  
based on the deformation argument. 
Equivalently, 
this number is the dimension of the space of inner products in
$\mathscr{X}$ associated with $\Diamond_l$; cf.
Theorem \ref{DIAZERO}. 

We note that a certain
$q,t$-generalization
of affine Demazure characters of any level $l>0$
was suggested  in 
\cite{CMa}; a connection is expected 
with paper \cite{Ka} upon the limits $t\to 0,\infty$. 

Given
any $f\in \mathscr{X}$ and using the constant term functional,
 the coinvariants 
ct$(f\,\Th[l]\,\mu)$ for any $l\ge 0$ and theta-functions
$\Th[l]$ of level $l$
are meromorphic functions in terms of  $k_\nu$. 
The formulas are explicit 
for $l=0,1$ and the nonsymmetric Macdonald
polynomials taken as $f$:  some products of binomials. 
 They are the generalized difference
Macdonald-Mehta identities.
 Also, one can use that $\Diamond_{l=1}$
is a Shapovalov anti-involution, which provides that 
the coinvariants for $l=1$
 are actually analytic upon some normalization. 
Employing ``picking up the residues" 
we arrive at ``the DAHA trace formulas" for any $l>0$.
\vskip 0.2cm

{\bf Non-compact theories.}
Let us briefly discuss the real integration. Here $k$ is
arbitrary complex and  there is no problem with an
analytic continuation to $\Re k<0$ for $k$ sufficiently close to
the real axis. The integration is 
$I^{re}_\pm(f)=
\int_{\pm \imath \ep \varrho+ \R^n} f(x) \mu(q^x)\, dx$ for $\ep>0$
and regular $\varrho \in \R^n$; 
the poles of $\mu$ at $\R^n$ must be avoided.
We can
set $\lan f,g\ran _{-l}^{re}=I^{re}_{\pm} 
\bigl(f\, T_{w_0}(g^\varsigma)\, q^{l x^2/2}\bigr)$,
where the Gaussian
ensures the convergence. 

We note that the Jackson
integration $J(f;\xi)$ is related to $I^{re}_+ - 
I^{re}_-$. In its turn, $J(f;\xi)$ is 
related to the imaginary integration, so we have some
connection between the imaginary and real integrations
via the Jackson integration. The latter is related 
to $\hat{\mathscr{I}}_+$. For instance, 
the  {\sf Jackson integration} of $f q^{x^2/2}$
for $l=-1$ and $\xi=-\rho_k$ is basically 
$\hat{\mathscr{I}}_+(f q^{x^2/2})/\Th(q^x)$, which is a constant
for any Laurent polynomial $f$.

Let us provide the adjustment of
the identity from (\ref{evinner}) to the real integration:
$\lan f,g\ran _{-1}^{re}=\lan 1,1\ran _{-1}^{re}
\Bigl(\tau_-(f)
T_{w_0}\bigl(\tau_-(g^\varsigma)\bigr)\Bigr)(q^{-\rho_k})$.
The formula
for $\lan 1,1\ran _{-1}^{re}$ is quite interesting.
For $A_1$, it is in terms of Appel functions due to Etingof;
see Section 2.3.5 of \cite{C101}. This is an indication that
we can try to replace $q^{-x^2/2}$ by $1/\Th(q^x)$ and connect
$\lan f,g\ran _{-1}^{re}$ with  
$\lan f,g\ran _{-1}^{im}$. 
The series for $1/\Th$ is of fundamental importance; 
see e.g. \cite{Car}.
\vskip 0.2cm

An important feature of the real (noncompact) theory
is that  $\mu(q^x)$ can be replaced by $\tilde{\mu}=
\mu^{-1}(q^x;q,t_\nu^{-1})$
from (\ref{reinner}). Everything in the real 
theory is up to {\sf quasi-constants}, which are $\hW$-periodic 
functions. Using this feature,
we can replace 
the denominator of $\hat{\mu}$ by 
the Gaussian with some corrections
ensuring the proper multiplicators upon the action of $P$.
This will ``eliminate" the denominator of $\mu$ and therefore
we can make
$\ep=0$ in the contour shift $\imath\ep\varrho$ above.    

\begin{theorem}\label{thm:noncom}
Let $h=(\rho^\vee,\th)+1$ be the dual Coxeter number, and
$$ 
M(x)=\sin\bigl(\pi (2\rho^\vee,x)\bigr)q^{h\frac{x^2}{2}}
X_{\rho}^{-1}
\prod_{\tilde{\al}>0}(1-t_\al^{-1}X_{\tilde{\al}}).
$$
Then the  pairing
$\int_{\R^n} f\, T_{w_0}(g^\varsigma)\, q^{l x^2/2}
M(x)\, dx$ is well defined for any
$q,t_\nu$ and real  $l>0$; it induces in $\mathscr{X}\ni f,g$ 
the anti-involution $\Diamond_{-l}$ for $l\in \N$.
\end{theorem}
{\it Proof.} We set $x_\al^\vee\equal
x_{\al^\vee}=(x,\al^\vee)$; recall that $X_a=q^{x_a}$.
Let us calculate explicitly the multipliers
of the functions under consideration
 upon the translations by $\om_j$.
For $1\le j\le n$, one has:
\begin{align*}
\om_j^{-1}\Bigl(\,&h\,x^2/2-(x,\rho)
\Bigr)=
h(x+\om_j)^2/2-(x+\om_j,\rho)\\
=\,
&h\,x^2/2-(x,\rho)
+h(x,\om_j)+h\om_j^2/2-(\om_j,\rho).
\end{align*}
The change is  $h(x,\om_j)+h\om_j^2/2-(\om_j,\rho)$. Next,
using $l(\om_j)=(2\rho^\vee,\om_j)$:
\begin{align*}
&\om_j^{-1}\Bigl(
\sin\bigl(\pi(2\rho^\vee,x)\bigr)\Bigr)\!=\!
\sin\bigl(\pi 
(2\rho^\vee,x+\om_j)\bigr)
\!=\!(-1)^{l(\om_j)}\sin\bigl(\pi(2\rho^\vee,x)\bigr).
\end{align*}

For the denominator $\De(q^{x})=\prod_{\tal>0}(1-X_{\tal})$ 
of $\tilde{\mu}$, which is
basically the denominator of the twisted Kac-Moody character
formula, one has:
$$\om_j^{-1}\bigl(\De(q^x)\bigr)\De(q^x)^{-1}\!\!=\!
\De(q^{x\!+\!\om_j})\De(q^x)^{-1}\!\!=\!
\prod_{\al>0}(-X_{\al}^{-1})^{(\al^\vee\!\!,\,\om_j)}
q^{-\frac{\nu_\al\de_j^{\al}(\de^{\al}_j-1)}{2}},
$$
where $\de_j^{\al}\!=\!(\al^\vee,\,\om_j)$.
It equals $
(-1)^{l(\om_j)}q^{-\sum_{\al>0}(\al,\,\om_j)x^\vee_\al
-\frac{(\al,\om_j)((\al^\vee,\,\om_j)-1)}{2}}$.

Then we use the standard identity:
 $\sum_{\al>0}(\al^\vee,u)\al=hu$, which holds for any $u\in \C^n$. 
For the sake of completeness, let us provide its proof. 
Setting $\sum_{\al>0}(\al^\vee,u)\al=\hat{u}$, $(\hat{u},v)$ is a
$W$-invariant symmetric form. We obtain that 
$(\hat{u},v)=c(u,v)$ for some constant $c$ due to the 
irreducibility of
$R$, and $(\hat{\th},\th^\vee)\!=\!c(\th,\th^\vee)\!=\!2c$.
Let us use now that  $(\al,\th^\vee)\!=\!1$ unless 
$(\al,\th^\vee)\!=\!0$ and $\al\!=\!\th$,
when it is $2$. Thus,
$(\hat{\th},\th^\vee)=\sum_{\al>0}(\al^\vee,\th)(\al,\th^\vee)=
(2\rho^\vee,\th)+(\th^\vee, \th)=2(\rho^\vee,\th)+2$ and $c=h$. 
Using the same identity,
$\sum_{\al>0}(\al,\om_j)(\al^\vee,\om_j)/2=h\,\om_j^2/2$ and:
$$
\om_j^{-1}\bigl(\De(q^x)\bigr)\De(q^x)^{-1}\!=\!
(-1)^{(2\rho^\vee,\,\om_j)}X_{\om_j}^{-h}
q^{-h\om_j^2/2+(\rho,\,\om_j)}.
$$

The convergence for $l>0$ is the same as it was for $\mu$. 
Integer levels $l>0$  are
needed here for $\Diamond_{-l}$ to serve the inner product. \sq
\vskip 0.2cm

The convergence holds here for $l=0$ when $f,g$ are of
sufficiently small degrees depending on $\Re k_\nu<0$.
This is exactly as in Theorem \ref{LEVZERO}, $(i)$.
Thus, such $f,g$ can be served by both,
the imaginary and real integrations, when $l=0$. The inner
product will be 
the same up to proportionality.
\vskip 0.2cm

{\bf The case of $A_1$ and $q$-zeta}. Let
 $x=x_{\om}, x_{\al}=2x$; the Gaussian is $q^{x^2}$
in terms of such $x$. We will omit $1$ in $\al_1,\om_1$. Then 
we obtain:  
$M(x)=\sin(2\pi x) q^{2x^2-x}(1-q^{2x})\prod_{i=1}^\infty
(1-t^{-1}q^{2x+i})(1-t^{-1}q^{-2x+i})$.

Accordingly, the pairings for $l\ge 0$ are
$\int_\R f(x)T(g(x))q^{lx^2} M(x)dx$ for any $k$ or 
$\int_{\imath\R} f(x)T(g(x))q^{-lx^2} \mu(x)dx$
 for $\Re k>0$ (subject
to the meromorphic continuation to $\Re k<0$). 

The integrals $\int_{\R}\frac{q^{x^2}}
{1\pm q^{x^2}}M(x)\,dx$ and similar ones 
lead to the definition
of the ``real" $q$-zeta function and 
Dirichlet $q$-$L$-functions studied
in \cite{ChZ}. The imaginary integration
results in their ``imaginary counterparts". Such integrands
ensure the convergence for $\int_{\R}$ and for $\int_{\imath \R}$,
but there will be now poles due to their denominators. 
\vfil

Upon some symmetrization needed for the
functional equation, the one for Dedekind's zeta,
 they conjecturally
 satisfy the Riemann hypothesis in terms of $s=k+\frac{1}{2}$
(Conjecture 6.3 at the end of \cite{ChZ}). There is another
version of {\sf RH} there without the symmetrization:
all ``interesting" zeros belong to one half-plane
with respect to  $\Re s=1/2$. 
\vfil

The analytic continuation
to $s<1/2$ is {\sf not} needed for the real integration 
and the Jackson-type summation. In the case of imaginary
integration, this continuation   
can be achieved using the pole decomposition and integral
formulas: the ones  we will do below, but with the contributions
of zeros of $1\pm q^{x^2}$.

The limit to the classical
$\ze(s)$ and the corresponding $L$-functions
$L(s)$ times some $\Gamma$ is when $q\to 1$. This is
generally for any $s$. However,
$\int_{\R}\frac{q^{x^2}}
{1- q^{x^2}}M(x)\,dx$ (for the minus sign here) will converge
to  $\sim \tan(\pi k)\Ga(k)^2$ for $\Re k< 1/2$ (i.e. for
$s<1$). It will converge to the (modified)
zeta for $\Re s>1$ in this (exceptional) case. Actually, it
will be like this even for $\Re s<1$ unless $\aa$ become very
large. 

The convergence is generally fast when $q<1$, even for 
$\aa\sim 1000$ or so (for reasonably small $\Im k$),
which  is  thanks to the Gaussians. 
\vskip 0.2cm

{\sf Six major DAHA theories.}
To summarize, we
mainly have two theories: the one based on the imaginary 
integration and that for the
real integration. In the Harish-Chandra theory,  
they are the  so-called {\sf compact\,} and 
{\sf noncompact\,}
cases. Totally, we have 
$6$ major theories by now, corresponding
to different choices of ``integrations":
from $(i)$ to $(vi)$. Namely, 
$(i)$ the usage of the constant term,
$(ii)$ imaginary integrations,
$(iii)$  real  integrations,
$(iv)$ Jackson integrations, $(v)$ the theory at
roots of unity, and $(vi)$ the theory at $|q|=1$
when Barnes' Gamma functions are needed.
Basically, any  $\hW$-invariant
integration can be taken for the corresponding
$\mu$-measure. 

The $\hW$-invariance
of the initial integration is immediate for 
$(i)$ and $(iii)-(v)$. We mostly stick to the 
imaginary integration (case $(ii)$) 
in this paper, which requires ``picking up the residues"
and integral formulas below for $\Re k_\nu<0$. 
This one  matches the $p$-adic theory and can potentially
admit some adelic version. 

Also, there is 
``DAHA-Satake theory", which is based on the  usage of
the {\sf affine} symmetrizers, $\hat{\mathscr{P}}_+$
and  $\hat{\mathscr{I}}_+$.  
The latter operator
is the summation over extended affine Weyl group ``twisted" by
the $\mu$-function, which 
is closely related to the Jackson integration. The former is
the {\sf affine $t$-symmetrizer}, which does not require 
any integration
(and the $\mu$-function), and certainly has some adelic 
generalization.

Let us mention  2 more directions:
DAHA theory  over finite fields 
and (related) theory
when $q$-Gamma and the corresponding $\mu$ 
are  replaced by those in terms of the $p$-adic Gamma. 

\comment{
\vskip 0.2cm
\centerline{
DAHA INTEGRATIONS:}

\begin{center}
\noindent
\begin{tabular}{|c||c|}
\hline
 imaginary ($|q|\neq 1)$            & real ($|q|\neq 1$)\\ 
 $\Downarrow$                 & $\Downarrow$  \\
 constant term ($\forall q$)  & Jackson sums \\ 
$\Uparrow$                    & $\Downarrow$  \\
 the case $|q|=1$ \hfill  $\Rightarrow$ 
                              & $\Rightarrow$ \hfill
                                 roots of unity\\
\hline
\end{tabular}  
\end{center} 
\vskip 0.2cm
}


\vskip 0.2cm
{\bf The p-adic limit}.
The $p$-adic limit is considered in detail in \cite{C101,CMa}.
Basically, $t^{-1}$ becomes the cardinality
of the residue field and $q\to 0$,
 but we need to be more exact here.

Let us use the homomorphism $\ze: \h_X\to \h_Y$ sending
$X_b\mapsto Y_b$ for $b\in P$,
 $T_i\mapsto T_i^{-1}$ for $1\le i\le n$, and
$t^{1/2}_\nu\mapsto t_{\nu}^{-1/2}$. We will extend it
elements in  $\h_X$ and  $\mathscr{X}$ depending
of $q$ by making $q\to 0$; the notation
will be  $\ze_0$.  Then
$\ze_0(\e_b)=\psi_b\equal t^{-l(b)/2}Y_b\p_+$ for $b\in P$,
 which are {\sf Matsumoto spherical functions}
in $\h \p_+$. Accordingly, the {\sf Satake-Macdonald $p$-adic 
spherical functions} 
are $\,\p_+\psi_b (b\in P_-)$, the images of the symmetrizations
of $\e_b$. 
 
Here, as above, $\e_b$ are nonsymmetric Macdonald
polynomials under the normalization $\e_b(q^{-\rho_k})=1$,
$\p_+=\sum_{w\in W}t^{-l(w)/2}T_w^{-1}/ 
\sum_{w\in W}t^{-l(w)}$.

For the AHA of type $A_1$, 
$$\psi_n\!\equal\! t^{-\frac{|n|}{2}}T_{n\om}\p_+,\ 
\p_+\!=\!(1+t^{1/2}T)/(1\!+\!t) \for n\in \Z.$$

Let $\mu^0\equal\lim_{q\to 0}\mu=
\prod_{\al>0}\frac{1-X_\al}{1-tX_{\al}}$. We set
$ \lan f,g\ran^0=
\text{ct}(f\,T_{w_0}(g^\varsigma)\mu^0)$, where $q\to 0$.
Recall that the ``$p$-adic trace" and  the
standard anti-involution in $\h_Y$ are as follows: 
$\lan T_{\hw} \ran_{reg}= \de_{id, \hw}$ and 
$T_{\hw}^\ast\equal T_{\hw^{-1}}$. We omit the complex conjugation
of the coefficients. 

We arrive at the following nonsymmetric
spherical AHA-Plancherel 
formula for any Laurent polynomials $f,g$ in terms of $X_b$:
\begin{align}\label{ahaplan}
&\lan f,g\ran^0
\mid_{t^{1/2}_\nu\!\mapsto t_{\nu}^{-1/2}}\ =
\sum_{w\in W}t^{\frac{l(w_0)-l(w)}{2}} \Bigl\langle (\ze_0(f)\p_+) 
(\ze_0(g)\p_+)^\ast\Bigr\rangle_{reg}.
\end{align}
It includes
the presentation of the Matsumoto spherical functions as nonsymmetric
Hall polynomials. 

The Gaussian and the action of projective $PSL_2(\Z)$ collapse
as $q\to 0$, and the definition of the Fourier transform
requires the characters of the (unitary) irreducible 
representations. In the $q,t$-setting,  
(\ref{ahaplan}) is direct from the action of
 {\tiny 
$\left(\begin{array}{cc}0& 1 \\
                       -1 & 0 \\
\end{array} \right)$} in DAHA.

\vskip 0.2cm

\section{\sc Meromorphic continuations}
For  $0<q=e^{-1/\aa}<1$ as above and $\upsilon\in \R$, we
set: 
$I^{\imath\aa}_{\upsilon}(f)\equal
\int_{\upsilon-\imath \pi \aa}^{\upsilon+
\imath \pi \aa}\!\!\cdots\!
\Bigl(\int_{\upsilon-\imath \pi \aa}^{\upsilon+\imath \pi \aa}
f(q^x)
\mu(q^x) dx_{\al_1}\Bigr)\!\cdots dx_{\al_n}.$ Here the order
of simple roots $\al_i$ can be arbitrary. Functions
$f(q^x)$ are assumed series in terms of 
$X_a (a\in Q)$ convergent 
in sufficiently large strips $|\Re x|<C$;
 the norm is the standard one in $\R^n$.
For the sake of convenience of notations we restrict ourselves
with $a\in Q$. If the whole polynomial representation is
considered, i.e. $X_b$ for $b\in P$, then the integrals 
$\frac{1}{M}\int_{\upsilon-\imath M \pi \aa}^{\upsilon+
\imath M \pi \aa}$ must be considered for proper $M\in 1+\Z_+$.

For the corresponding $I^{im}_\upsilon$ we integrate over  
for $\upsilon\rho^\vee +\imath \R^n$; 
 $f(x)$ can be any entire functions in sufficiently
large strips $|\Re x|<C$
provided the integrability of $|f|$ in the imaginary
directions. The notation $I^{im}$ used above is for $\upsilon=0$.

\vskip 0.2cm

\begin{theorem}\label{thm:iterint}
(i) Let
$0\!<\!\upsilon\!\le\!\frac{1}{h}$ for the dual Coxeter number
$h\!=\!(\rho^\vee,\vth)\!+\!1$, $\{\al_i\}$ be a fixed set of
simple roots taken in any order. Then
the corresponding iterated integral for
 $I^{im}_\upsilon(f)\equal\int_{\imath \R^n+\upsilon
\rho^\vee}f\mu dx$  is a meromorphic continuation of 
$I^{im}(f)=I^{im}_0(f)$ from
$\Re k_\nu>0$ to $\Re k_\nu>-\ep$ for some $\ep>0$.
 The same holds for $I^{\imath\aa}_\upsilon(f)$ assuming that
functions $f$ are in terms of $X_a, a\in Q$.
The meromorphic continuation
of $I^{im,\imath\aa}_\upsilon(f)$ to any $\Re k_\nu<0$
can be presented as a {\sf finite} linear combination
 of integrals
over certain {\sf $\mu$-residual subtori}, 
with the leading term $I_0^{im,\imath\aa}$. 
The number of such integrals 
grows as  $|\Re k_\nu|$ increase. 

(ii) We define
$\Si^{\imath\aa}(f)$ as the sum of
the residues of  $\mu(q^x)f(q^x)$ over 
{\sf $\mu$-residual points} $\xi$ subject to the
consecutive inequalities  
 $\Re x_{\al_i}\!>\!\upsilon$ for $1\le i\le n$
imposed  when taking the
iterated integrals $\int_{\upsilon-\imath\pi\aa}^
{\upsilon+\imath\pi\aa}
\{\cdots\} dx_{\al_i}$. 
The points $\xi$ that occur here  depend 
on  the order of $\{\al_i\}$, but not on $\Re k_\nu$.
This (infinite) sum is
convergent for sufficiently small $\Re k_\nu \le 0$ and
extends meromorphically the 
analytic function  $\frac{I^{\imath\aa}_0(f)}
{(2\imath\pi\aa)^n}$ from $\Re k_\nu>0$
to {\sf any} $\Re k_\nu\le 0$ provided the 
convergence of $\Si^{\imath\aa}(f)$. The residues in this sum
are essentially the values of $\mu$ upon the deletion of the 
binomials vanishing at the corresponding $\xi$ in the
setup of Theorems \ref{thm:residues}, \ref{thm:resid}.
\end{theorem}
{\it Proof.} The fact that $I_{\upsilon}^{im,\imath\aa}$ 
extend $I_{0}^{im,\imath\aa}$ analytically to small negative $\Re k_\nu$
is straightforward. Generally, we determine the corrections when
moving the contours of integration
by $\upsilon\rho^\vee$; they are
iterated integrals over $\mu$-residual subtori of smaller dimensions.
Let us take  $I_{\upsilon}^{\imath\aa}$ for the sake of
concreteness. First, we replace
$\int_{\upsilon-\imath \pi \aa}^{\upsilon+
\imath \pi \aa}(\cdots)
d_{x_{\al_i}}$  by 
$\int_{-\imath \pi \aa}^{
\imath \pi \aa}(\cdots)
d_{x_{\al_i}}$ for every $1\le i\le n$. 
The corresponding correction will be 
$$
\int_{\upsilon-\imath \pi \aa}^{\upsilon+
\imath \pi \aa}\cdots\Bigl(
\int_{\upsilon-\imath \pi \aa}^{\upsilon+
\imath \pi \aa}-\int_{-\imath \pi \aa}^{
\imath \pi \aa}\Bigr)
\cdots\int_{\upsilon-\imath \pi \aa}^{\upsilon+
\imath \pi \aa} f(q^x)
\mu(q^x)\, dx_{\al_1}\cdots dx_{\al_n},$$
where the difference is at place $i$. It is 
a finite  sums of integrals over the
proper (imaginary) contours of dimension $(n-1)$. 
The integrands will be some (partial) 
residues for the corresponding $x_{\al_i}$. 
Then we continue inductively: replace all remaining 
$\int_{\upsilon-\imath \pi \aa}^{\upsilon+
\imath \pi \aa}$ by $\int_{-\imath \pi \aa}^{
\imath \pi \aa}$ in the same way. The final output will be
a finite sum of integrals over certain $\mu$-residual subtori.
It will depend on the  order of integrations.
The coefficients in this sum  will be
the corresponding (partial) residues of $\mu$. 

\vskip 0.2cm

{\sf (ii).}
Taking the iterated integrals in terms of the residues in the 
corresponding right
half-planes requires explanations. We will provide the exact algorithm
for finding the set of all $\xi$ that occur
in the pole decomposition of $I_{\upsilon}^{\imath\aa}$;
all of them are $\mu$-residual points,
but not all will occur, which significantly
depends on the order of integrations.

\vfil
The description below is purely combinatorial;
it suffices to assume that
$1\gg \upsilon\gg -k_\nu>0$. We will set $\bar{\al}_i=
(\al,z)$ for any $\al\in R$, where $z\in \R^n$.
For instance, $\bar{\mu}=\mu(q^z)$.
Given a pole $\xi$ of $\mu(q^z)$,
let $\{\tbe_i=[\be_i, m_i]\}$ be a {\sf sequence} of consecutive 
binomials that result from the iterated integrations, where 
 Generally, they can be different as (unordered) sets for
different sequences and the same set can occur more than once.
Here $i=1,\ldots, n$,
$m_i\in \nu(\be_i)\Z_+$,  $m_i>0$ for $\be_i<0$. The corresponding
$\xi$ modulo the periods of $X_\al$
will be a unique solution of the  system of equations
$\bar{\be}_i+m_i+k_{\be_i}=0$ for $1\le i\le n$.

We will treat in the following algorithm  
$\bar{\al}_j$ as undetermined variables, which
will be eliminated one by one until we obtain their values
at $\xi$. 
\vfil 

One has for $i=1$:\,
 $\be_1=\sum_{j=1}^n c_j^1\al_j$, where $c_1\neq 0$, 
 $\bar{\be}_1=\sum_{j=1}^n c_j^1\bar{\al}_j$. We impose then
the equation
$k_{\be_1}+\bar{\be}_1+m_1=0$ and obtain that $\bar{\al}_1=
-(m_1+k_{\be_1}+\sum_{j=2}^n c_j^1\bar{\al}_j)/c_1^1>0$.
Here $\Re \bar{\al}_j$ will become $\upsilon$ in the following
integrations: for $j=2,3,\ldots, n$. Let us use that 
$k_\nu$ and $\upsilon$ are assumed very small. Then we arrive at 
$-(m_1+ k_{\be_1}+C\upsilon)/c_1^1>0$ for some $C$ with the upper
bound depending only on the root system $R$ and the terms
with $k$ and $\upsilon$ can be disregarded.  We obtain that
the initial inequality 
holds if and only if $c_1^1<0$ and $m_1>0$. The former means
that $\be_1<0$ (then $m_1>0$ anyway). Equivalently, 
$(1-t_{\be_1}X_{\tbe_1})$ belongs 
to the ``negative
half" of the denominator of $\mu$. 
This is so only for the 1{\small st} integration;
the ``positive half" of $\mu$ may contribute too. 

To go to the second step, 
we set $\al_1^\bullet=\bar{\al}_1(k_\nu\mapsto 0)=
\bar{\al}_1\!\mid_{k_\nu\mapsto 0}$
and replace  $\bar{\al}_1$ by the relation above
in all $\bar{\be}_i,\bar{\al}_i$ for  $i>1$ and
all binomials of $\mu$. 
The one with 
$t_{\be_1}X_{\tbe_1}$ in the denominator of $\mu$
 will be deleted and
we perform the reduction of coinciding or proportional
binomials in the numerator and denominator of $\bar{\mu}$. 
The binomials with $t_{\be_i}X_{\tbe_i}$ for $i>1$ will not be 
reduced  by this construction since they are among 
the defining relations for $\xi$.
\vfil

We arrive at {\sf new} 
$\bar{\be}_i$ and
$\bar{\al}_i$ for  $i>1$ and $\bar{\mu}$ in terms
of $\bar{\al}_i (i>1)$ and $m_1$.  
By construction,  $(\bar{\al_i},\xi)= 
(\al_i,\xi)$ and  $(\bar{\be_i},\xi)= 
(\be_i,\xi)$, where 
 $(\al\!+\!c,\xi)\!=\!(\al,\xi)\!+\!c$ here and below for $c\in \Q$.

\vskip 0.2cm

Then, we represent 
$\bar{\be}_2=\sum_{j=2}^n c_j^2\bar{\al}_j$ for
{\sf new} $\bar{\be}_2$ and $\bar{\al}_j$ ($j\ge 2$), where $c_2^2\neq 0$,
and obtain: $\bar{\al}_2=
-(m_2+k_{\be_2}+\sum_{j=3}^n c_j^2\bar{\al}_j)/c_2^2$. 
The 2{\small nd} positivity condition is:
$
-(\,m_2+\sum_{j=3}^n c_j^2\bigl(\bar{\al}_j\bigr)^{k\to 0}
_{\al\to 0}\,)/c_2^2>0,
$
where $\bigl(\bar{\al}_j\bigr)^{k\to 0}_{\al\to 0}$ means that we
delete all $\al$ and $k$ from $\bar{\al}_j$, i.e. keep 
only constants, which are in terms of $\{m_i\}$. 
Note that $c_2^2$ is not the coefficient of $\al_2$
in $\be_2$.
Then we switch to {\sf new}
 $\bar{\al_i},\bar{\be}_i$ for $i>2$
and $\bar{\mu}$ as above using the formula for $\bar{\al}_2$ 
and continue by induction.

Finally, we obtain the complete list of
substitutions  $\bar{\al}_i\mapsto
\al_i^\bullet=\sum_{j>i}C_{ij}\bar{\al}_j+M_i+K_i$ with some $M_i$
in terms of $m_j$ for $j\le i$ and $K_i$ in terms of
$k_\nu$. This  gives 
the formulas for 
$\{\bar{\al}_i\}$ in terms of $\{m_i\}$ and $k_\nu$,
and the list of inequalities for $\{m_i\}$. The latter
are $\bigl(\bar{\al}_i\bigr)_{\al\to 0}^{k\to 0}>0$
for the corresponding substitution formulas. 
These inequalities
are  necessary and sufficient
for $\xi$ to occur in $\Si^{im,\imath\aa}$ for a given
order of integrations.  However, 
the corresponding 
residue can be $0$ or there can be cancelations of the terms. 

For instance, the  $n${\footnotesize th}
 step (the last) gives that $\bar{\be_n}=c_n^n\bar{\al}_n$, 
$\bar{\al}_n=\bar{\be_n}/c_n^n=
-(k_{\be_n}+m_n)/c_n^n$ and the inequality is
$\bigl(\bar{\al}_n\bigl)_{\al\to 0}^{k\to 0}\,=
-m_n/c_n^n>0$. 
Thus, $c^n_n<0,m_n>0$ and we have
$(\al_n,\xi^\bullet)=-m_n/c^n_n>0$,
where $\xi^\bullet\equal(\xi)^{k\to 0}
\in \sum_{i=1}^n \Q\om_i$ is obtained when we solve the
system above with 
$\bar{\al}_i\mapsto (\bar{\al}_i)^\bullet$. 
Actually, $(\al_n,\xi)>\upsilon$ by construction, which
gives $(\al_n,\xi^\bullet)>0$.
Generally, we arrive at the
following description of $\xi$.

\begin{lemma}\label{lem:xi-cond}
Let $\be_i=[\be_i]+\lan\be_i\ran$, where 
$[\be_i]=\sum_{j=i}^n c_j^i\al_j$,
$\lan \be_i\ran \in \sum_{j<i} \Q\be_j$. Accordingly,
$(\bar{\be}_i)^\bullet+m_i=[\be_i]+M^i$ at step $i$, where
$M^i=\lan\be_i\ran_{\be_j\to -m_j}$ is expressed in terms of $m_j$ 
for $j\le i$. Then
$(\xi^\bullet,[\be_i])=-M^i$ and the defining inequalities 
for $\xi$ become  $(\xi^\bullet,[\be_i]/c_i^i)=-M^i/c_i^i>0$ 
for $i=1,\ldots,n$. 
\sq
\end{lemma}



This procedure and the inequalities for $\xi^\bullet$ 
depend on the order of integrations. These inequalities do not 
guarantee that such $\xi$ 
occur only once  and with nonzero coefficients;
there can be
some cancelations even for $A_3$ (see below). 
\sq

\vfil

{\bf Concluding remarks.}
{\sf (a).} The $p$-adic limit of the sums over
residual subtori 
 from $(i)$ for $I^{\imath\aa}_{\upsilon}(f)$
is as follows. We assume that $0>k_\nu>-\ep$
for small $\ep$, take $f\in \C[X_a,a\in Q]$ and
replace 
the integrals 
$\frac{1}{2\pi \imath\aa}\int_{-\pi\imath\aa}^
{\pi\imath\aa}$ by
$\frac{1}{2M\pi \imath\aa}\int_{-M\pi\imath\aa}^
{M\pi\imath\aa}$ for $M\in 1+\Z_+$. By the way,
the usage of $M$ here allows us to incorporate $f=X_b$ for $b\in P$
instead of $X_a$ with $a\in Q$. 
 
Then we set 
 $\aa=\frac{1}{M}$, $k_\nu=c_\nu \aa/\nu$ for $c_\nu<0$
 and make  $M\to \infty$. This results in 
$q\!\to\! 0$, $\Re k_\nu\to 0_-$.
$t_\nu\!=\!e^{-\frac{\nu k_\nu}{a}}\!\to\! e^{c_\nu}<1$. We
arrive at the integrals over AHA residual subtori 
and the formulas from \cite{HO}. Actually, there is one more step
here: the $W$-symmetrization. 

Recall that $t\mapsto 1/t$ when we go from
 DAHA to AHA with the {\sf standard} meaning 
of the parameter $t$ there, which is $|\mathbb F|$ classically. So 
the range $\Re k_\nu<0$ or $t_\nu\!>\!1$ in terms of $t$ for DAHA 
corresponds to  
$t_\nu\!<\!1$ in the {\sf standard} AHA setting, when the
discrete series occurs. 

\vskip 0.2cm

{\sf (b).}
As we already discussed, $I^{im}(f)$ can be generally 
reduced to $I^{\imath\aa}(f)$. Namely, we replace
$f\!\mapsto\! \sum_{b\in 2\pi\imath\aa P^\vee}\! f(z+b)$ if this
sum converges. Then, the integral formulas in terms of
integrals $I^{\imath\aa}(f)$ over residual subtori
for Laurent polynomials $f$ and series coincide
with  ct$(f\mu)$; to be exact, they are proportional in
the corresponding range of $k_\nu$.  Then, 
given $f\in \mathscr{X}$,  the constant term ct$(f\mu)$
is meromorphic for any $k_\nu$, 
which can 
seen directly using the formulas for the $E$-polynomials.

There are  3 more {\sf algebraic} ways to calculate ct$(f\mu)$.
One can use 
$(i)$ the coinvariants for $\Diamond$,
$(ii)$  affine symmetrizers 
$\hat{\mathscr{I}}_+$ and  
$\hat{\mathscr{P}}_+$, and $(iii)$ the 
Jackson integrals $J_\xi(f)=J(f;\xi)$. Here $(ii-iii)$ require
$\Re k_\nu<0$. The Jackson integrals are
the closest to 
$\Si^{\imath\aa}(f)$; they are related  to
the formulas via real integrations. The latter 
provide 
another tool for obtaining  ct$(f\mu)$ and 
result in {\sf non-compact trace formulas.} The existence
of the affine symmetrizer in $\mathscr{X}$ for sufficiently
small $\Re k<0$ and in $\mathscr{X}q^{lx^2/2}$ for $l>0$ and
$\Re k < 1/h$ is remarkable; these modules behave as
discrete series representations in AHA theory.
\vfil

The
problem with the usage of the series
$\Si^{\imath\aa}(f)$ 
is that its convergence holds for 
restricted classes of $f$ and heavily depends on $\Re k$. This 
is similar to Jackson integrals. For instance, 
only constants
can be taken as $f$ for $\Si^{\imath\aa}(f)$
among $W$-invariant Laurent polynomials $f\in \mathscr{X}$ when 
$\Re k_\nu<0$ are close to $0$; cf. \cite{Ma3}.
 Also, the Gaussians $q^{-lx^2/2}$
for $l>0$ and their Laurent
expansions  diverge in the real
directions and result in divergent $\Si^{\imath\aa}(f)$. The
$\HH$-modules  $q^{- lx^2/2}\mathscr{X}$ are important in the
DAHA theory, but $\Si^{im,\imath\aa}$ cannot be used for them. 
\vskip 0.2cm

By contrast, the {\sf finite} sums over residual subtori
from $(i)$ can be used for practically arbitrary analytic 
functions provided
the integrability. For the real integration, i.e. in
the {\sf non-compact case}, even a single integral can
be used. 
We did not state the
uniqueness of the integral formulas 
explicitly. The $\HH$-invariance is one
way to fix them uniquely. 

The other way is by combinatorial 
collecting the $\mu$-residual points from $\Si^{\imath\aa}$
in the families corresponding to
$\mu$-residual subtori, which is a canonical process. Fig.
\ref{ct-a2-2}  for $A_2$
demonstrates this. One needs to find the pole decompositions
for the integrals over the residual subtori for $\Re k_\nu <0$
for this, starting with $I_0^{im,\imath\aa}$. 
Obtaining explicit integral
formulas for imaginary integrations is not simple:   
 the number of residual tori and
points that occur there
grows when $|\Re k_\nu|$ increases. 
\vskip 0.2cm
\vfil

{\sf (c).}
As an application to $\mathscr{X}$, one can consider {\sf singular}
$t_\nu=q_\nu^{k_\nu}$ 
such that the coefficients in the integral formulas 
have poles; we renormalize the integral formula making some
coefficients (measures) without $k$-poles. 
 Then functions $f$ vanishing
at all subtori with singular coefficients
form  an $\HH$-submodule, and 
 $I_0^{\imath\aa}\bigl(f(q^x) \bar{g}(q^{-x})\bigr)$
will induce a {\sf positive definite} inner product there for the
complex conjugation $g\mapsto \bar{g}$ of the coefficients
of $g\in \mathscr{X}$. The corresponding 
anti-involution of $\HH$ was calculate
 in \cite{C101}; it 
is not $\Diamond$ since $T_{w_0}$ is omitted.
 
It was 
proven under some technical conditions in \cite{C6} that
a certain ``smallest" submodule of $\mathscr{X}$ is $Y$-semisimple
for any singular $q,t_\nu$. The technique of
intertwiners was used; the inner products were not involved. 
 The integral formulas 
provide an alternative approach to this and similar facts, 
including generalizations 
to other spaces of functions. 

\vskip 0.2cm

{\sf (d).}
The usage of $\upsilon=1/h$ for the translation
 $\upsilon \rho^\vee$ of the contour of integration 
in the theorem 
gives the greatest analyticity range of $k_\nu$ in 
$I^{im,\imath\aa}_{\upsilon}(f)$. Recall that 
arbitrarily small negative $\Re k_\nu$ in
$(ii)$ were used to define $\Si^{\imath\aa}$, but making
the analyticity range ``optimal" for $I^{im, \imath\aa}_\upsilon(f)$
is of importance. The main fact is that if 
$\Si^{\imath\aa}$ is known, where only small $\Re k_\nu<0$ are 
sufficient, the corresponding series provides the required
meromorphic continuation to {\sf any} $\Re k_\nu <0$ 
provided the convergence. 
 will be the required meromorphic 
continuation
for all negative $\Re k_\nu$ assuming the convergence.
Also, 
$\rho^\vee/h$
 is invariant under the affine action of $\Pi=P/Q$
in $\R^n$, which provides additional symmetries of 
 $I^{im,\imath\aa}_\upsilon(f)$. After we establish
the pole decomposition of  $I^{\imath \aa}$
for {\sf arbitrarily} small negative $k_\nu$, the same 
formula will work for {\sf any} negative $k_{\nu}$
ensuring the convergence of the resulting series. 
\vskip 0.2cm

{\sf (e).}
The conditions from $(i)$ of Theorem \ref{LEVZERO}
for $f=X_a$
are sufficient for
the convergence in $(ii)$ of Theorem \ref{thm:iterint}
 but they are not necessary. For instance,
it  converges for $f=X_{ca}$ for certain $a\in Q\cap P_+$ and 
{\sf any} 
$c\in \Z_+$.
These cones are 
nonempty for any orders of integrations. For $A_n$, any direction $a$
can be made such for a proper order of integrations.  

In the case of $I^{im}$ for $\Re k_\nu<0$, the
convergence is granted for {\sf Paley-Wiener functions} 
(for the Laplace transform). Namely, it suffices to assume 
that $f(x)$ is analytic in $x\in \C^n$ such that
for every positive $N$ there exists some $C_N>0$ such that 
 $|f(w(x))| \leq C_N(1+|x|)^{-N}  e^{B|\Re(x)|}$ for
every $w\in W$, where 
 $0<B< A|\Re(\rho_k)|$ for some $A>0$.
\vskip 0.2cm

Let us emphasize that if the meromorphic continuation
of $I_{0}^{im,\imath\aa}(f)$ from $\Re k_\nu <0$ to $\Re k_\nu\ge 0$
exists, then it is unique and 
does not depend on the specific choice of variables and
their order of integrations.  However, the 
sums $\Si^{im,\imath\aa}(f)$ 
and the corresponding growth conditions for $f$ 
depend on the order of integrations.
This leads to some nontrivial identities. 
Our integral formulas generally
depend on the order of integrations too. Their uniqueness 
is under some assumptions.

\section{\sc Pole expansion for 
\texorpdfstring{{\mathversion{normal}$A_n$}}{An}}
The combinatorial algorithm for finding $\Si^{im,\imath\aa}$
becomes relatively
simple for $A_n$. We will provide 
$\Si^{\imath\aa}$ only for the standard order
of $\al_i=\vep_i-\vep_{i+1}$ in the iterated integral,
though see below an example for $A_3$ with $\al_1,\al_3,\al_2$.
We will describe $\xi$ and $\xi^\bullet=\xi(k\to 0)$ following  
the proof
of Theorem \ref{thm:iterint}. 

\begin{theorem} \label{thm:Anres} 
For $A_n$ and the standard sequence of $\al_i$ as above,
let  $X=\{\xi\}$ be the set of the 
$\mu$-residual points in $\Si^{\imath\aa}$. Then
the relations for $\xi^\bullet=b\in P$ are
$(b,\al_n+\cdots+\al_i)>0$ for $i=1,\ldots,n$.
The corresponding $\xi$ are $\pi_b(-k\rho)=b-u_b^{-1}(k\rho)$,
denoted simply by $\pi_b$ in 
Theorem \ref{thm:residues},$(i)$ for the initial
$\xi=-k\rho$. One has:
 $\Si^{\imath\aa}(f)\!=\!\sum_{b\in X} Res_{\pi_b}f(\pi_b)$, where
$f(bw)\!=\!f(q^{b\!-\!w(k\rho)})$.
The residues here are 
as in Theorem  \ref{thm:residues},(iv).

In the integral formulas from 
Theorem \ref{thm:iterint}, the integrands and residues
are obtained by deleting the binomials from
the denominator of $\mu$ vanishing at the
corresponding $\mu$-residual torus followed by the
evaluation  of  $f(q^x)$ 
and the rest of $\mu(q^x)$ at the corresponding tori and points. 
\end{theorem}
{\it Proof.}
The description of $\{\xi\}$ follows directly
from the explicit algorithm
from Theorem \ref{thm:iterint}. Due to
Theorem \ref{thm:resid}, the $\mu$-residual points $\xi$ can be
potentially with zeros in the numerator of $\mu$ for $A_{n>2}$. 

For $A_3$,
such $\xi$ is as follows up to the action of $\hW$:
$\bar{\al}_2=\bar{\vep}_2-\bar{\vep}_3=0$ and 
$\bar{\vep}_i-\bar{\vep}_1+1=-k=
\bar{\vep}_4-\bar{\vep}_i+1$ 
 for $i=2,3$. One has
$\kapp_1-\kapp_0=4-1=3=n$; so it is $\mu$-residual. 

We claim that $\mu$-residual $\xi$ with binomials
in the numerator of $\mu$ vanishing at $\xi$ do 
not contribute to $\Si^{\imath\aa}$.
Let us outline the justification. 

Such $\xi$  have
nontrivial stabilizers $W_\xi$ in $W$;
for instance, $s_2(\xi)=\xi$ for the $\xi$ above. 
The group $W_\xi$ will permute the corresponding $\tbe_i$
and these permutations are non-trivial unless for id$\in W_\xi$. 
Accordingly, such $\xi$ will occur $|W_\xi|$  times 
in the procedure of finding $\Si^{\imath\aa}$ from the proof
of Theorem \ref{thm:iterint}.
Following Lemma
\ref{lem:stabxi}, which states that the Jackson integrals $J_\xi(f)$
vanish if $|W_\xi|>1$,
we check that the corresponding residues will
cancel each other in this orbit.

This argument is generally
applicable to other cases in Theorem \ref{thm:iterint}, but we claim
the absence of $\xi$ with non-trivial stabilizers 
only for $A_n$. Note that Lemma \ref{lem:stabxi} cannot be applied
if $\kapp_1-\kapp_0>n$. However, such points will not 
occur for $A_n$,
which is straightforward to verify.

Concerning the inequalities for $b$, $\bar{\be}_n=
c_n^n\al_n$ for the last step and $(b,\al_n)>0$, as 
it was checked in the proof of 
Theorem \ref{thm:iterint} for any root systems. One has  
$c_n^n=-1$ for $A_n$.  Similarly,
$\bar{\be}_{n-1}^\bullet$ can be $-\al_{n-1}$
or $-\al_{n-1}-\al_n$ for the $(n-1)${\footnotesize th}
step, which gives $(b,\al_{n-1}+\al_n)>0$.
Then we continue by induction using the following general lemma.

\begin{lemma}\label{lem:redmin}
For a minuscule $\om_r$ such that $r$ is an endpoint
of the corresponding
Dynkin diagram of $R$, let $\al\mapsto \al'$ be 
the deletion of $\al_r$ if it is present in $\al\in R$. 
$(\al+\be)'=\al'+\be'$ for any $\al,\be\in R$
and $\al'\neq 0$ is a root in the root system $R'$ with
simple roots $\{\al_i, i\neq r\}$.
\sq
\end{lemma}

According to Lemma \ref{lem:xi-cond},
we need to find $[\be_i]$ for $i=1,\ldots,n$
in the following decompositions:
$\be_i=[\be_i]+\lan\be_i\ran$, where 
$[\be_i]=\sum_{j\ge i} c_j^i\al_j, $
$\lan \be_i\ran \in \sum_{j<i} \Q\be_j$ and
the relations ensuring that $\xi^\bullet$
occurs in $\Si^{\imath\aa}$ are 
$(\xi^\bullet,[\be_i]/c_i^i)>0$ and $\tbe_i=[\be_i,m_i]\in \tR_+.$

The following procedure is generally applicable to any roots
system $R$ with at least one minuscule weight. Namely,
$\al_i$ must be in the form
$\al_r$ from Lemma \ref{lem:redmin} for the system $R^{(n-i+1)}$
with simple roots $\{\al_i,\cdots,\al_n\}.$
It really gives that $\xi^\bullet=b\in P$. 
However, it can result in $b$ that do not actually occur
as $\xi^\bullet$ due to cancelations. Also,
the conclusion that $\xi=\pi_b(-k\rho)$ is for $A_n$ only.
We will check that  $\xi=\pi_b(\xi^0)$ for some $\xi^0$, but the
latter can be not of Steinberg type.    
\vskip 0.2cm

We obtain that $\xi^\bullet=b\in P$ since $c_i^i=\pm 1$,
which gives that $\xi=\pi_b(-k\rho)$ for $A_n$. 
Then 
all $c_i^j$ for fixed $i$ and $j\ge 0$
have the same sign (if not zero)
because $\al_i,\cdots,\al_n$ are simple roots in the
corresponding system $R^{(n-i+1)}$. 
Next, let us check that $c_i^i=-1$.

For the last step, $[\be_n]=c_n^n \al_n$ and 
$-m_n=c_n^n(\al_n,b)$, where $(\al_n,b)>0$. We obtain
that $c_n^n<0$, which we already know (for any root systems).
For the $(n-1)${\tiny st} step: $[\be_n]=c_{n-1}^{n-1}(\al_{n-1}+c\al_n)$
for $c\ge 0$ and $-m_{n-1}/c_{n-1}^{n-1}=
\bigl((\al_{n-1},b)+c(\al_n,b)\bigr)>0$. This  gives that 
$c_{n-1}^{n-1}=-1$. Then we go to $c_{n-2}^{n-2}$ and so on. 

Finally, we obtain for $A_n$ and the standard order
of the integrations that 
the all possible sequences
are  $\{[\be_n],[\be_{n-1}],[\be_{n-2}],\ldots\}$ are:

\noindent 
$\{\al_n,\, \al_{n-1} \text{ or } \al_{n-1}+\al_n,\,
 \al_{n-2} \text{ or } \al_{n-2}+\al_{n-1}
\text{ or } \al_{n-2}+\al_{n-1}+\al_n,\,\ldots\}.$
They all occur and 
result in the statement of
the theorem.
\sq

\vskip 0.2cm
{\bf Examples.}
{\sf (a).} Let us consider the standard order of $\al_i$ for  $A_3$.
For example, let  $\tbe_1=[-\al_1-\al_2-\al_3,n_1+1],
\tbe_2=[\al_1,n_2], \tbe_3=[\al_2,n_3]$, where $n_i\ge 0$. Then   
the consecutive substitutions for $k=0$ are 

\noindent
$\bar{\al}_1\mapsto 
-\bar{\al}_2-\bar{\al}_3+n_1\!+\!1, 
\bar{\al}_2\mapsto -\bar{\al}_3+n_1\!+\!n_2\!+\!1,
\bar{\al}_3=n_1\!+\!n_2\!+\!n_3\!+\!1.$

 Accordingly, 
$b=(n_1+n_2+n_3+1)\om_3+(-n_3)\om_2+(-n_2)\om_1$ and
$\xi=b+u_b^{-1}(k\rho)=
(n_1+n_2+n_3+1+3k)\om_3+(-n_3-k)\om_2+(-n_2-k)\om_1$.

For the sequence
$\tbe_1=[-\al_1-\al_2,n_1+1], \tbe_2=[-\al_2-\al_3,n_2+1],
\tbe_3=[\al_2,n_3]$ the substitutions (under $k=0$) are: 
$$\bar{\al}_1\mapsto -\bar{\al}_2+n_1\!+\!1,
\bar{\al}_2\mapsto -\bar{\al}_3+n_2\!+\!1,
\bar{\al}_3=n_2\!+\!n_3\!+\!1.$$
 Finally,
$b=(n_2+n_3+1)\om_3+ (-n_3)\om_2+ (n_1+n_3+1)\om_1$ and 
 $\xi=(n_2+n_3+1+2k)\om_3+ (-n_3-k)\om_2+ (n_1+n_3+1+2k)\om_1$.
Recall that  $\tbe_i$ are from the binomials that are taken
for the corresponding integration: $d x_{\al_1}$, $d x_{\al_2}$ and
$d x_{\al_3}$.  
Here and below $n_i\in \Z_+$.
\vskip 0.2cm

Let us provide
the whole set of $\xi$ for $A_3$. The first three
numbers in the list below give the types of 
the $\tbe_i (i=1,2,3)$, their corresponding 
numbers in the following sequence of 12 types:

\noindent
{\small  
$[\al_1,m_1],\, [\al_1+\al_2,m_2],\, 
[\al_1+\al_2+\al_3,m_3], \,[\al_2,m_4],\, [\al_2+\al_3,m_5], \,
[\al_3,m_6],$}\\
{\small $[-\al_1,m_7], [-\al_1-\al_2,m_8], 
[-\al_1-\al_2-\al_3,m_9], [-\al_2,m_{10}], [-\al_2-\al_3,m_{11}]$}
and {\small $[-\al_3,m_{12}]$} (number 12) in $\tR_+$; they
are from the denominator 
of $\mu$. Here $m_i=n_1$ for $1\le i\le 6$ and $m_i=n_1+1$ for
$i>6$, i.e. in the ``negative half" of $\mu$.
The examples above correspond to $\{9,1,4\}$ and $\{8,11,4\}$. 
We have:
{\small
$
\{7,10,12,1+k+n_1 ,1+k+n_2 ,1+k+n_3 \},
\{7,11,2,1+k+n_1 ,-1-2 k-n_1\! -\!n_3 ,2+3 k+n_1\! +\!n_2\! +\!n_3 \},
\{8,1,12,-k-n_2 ,1+2 k+n_1\! +\!n_2 ,1+k+n_3 \},
\{8,11,4,1+2 k+n_1 +n_3 ,-k-n_3 ,1+2 k+n_2\! +\!n_3 \},
\{9,1,4,-k-n_2 ,-k-n_3 ,1+3 k+n_1 +n_2\! +\!n_3 \},
\{9,10,2,-1-2 k-n_2\! -\!n_3 ,1+k+n_2 ,1+2 k+n_1\! +\!n_3\}.
$
}  

{\sf (b).}
For the order $dx_{\al_1}dx_{\al_3}dx_{\al_2}$ of integrations,
let $\tbe_i=[\be_i,m_i], m_i=n_i+1$ for $\be_i<0$ and $m_i=n_i$
otherwise. The corresponding
families of $\xi^\bullet$ are:
{\small $
\{7,11,6,1+n_1,1+n_2\!+\!n_3,-n_2\},
\{7,12,10,1+n_1,1+n_2,1+n_3\},
\{8,11,3,-1-n_2\!-\!n_3,2+n_1\!+\!n_2\!+\!n_3,-1-n_1\!-\!n_2\},
\{8,12,1,-n_2,1+n_1\!+\!n_2,1+n_3\},
\{9,1,6,-n_3, 1+ n_1\!+\!n_2\!+\!n_3,-n_2\},
\{9,2,10,-1-n_2\!-\!n_3,1\!+\!n_2,1+n_1\!+\!n_3\},
$
}\ where $n_i\in \Z_+$ as above and we transpose $n_2$ and $n_3$
to match $n_i$ used in the families with the ones
for $\{\al_1,\al_2,\al_3\}$ above.


Let $b\!=\!\sum_{i=1}^3 b_i \om_i \in P$.
Then all families above satisfy the inequalities 
$b_2>0,b_2+b_3>0$. 
Imposing them,
families $\{7,11,6\}$ and $\{7,12,10\}$ are given by $b_1>0$,
families $\{8,12,1\}$ and $\{9,2,10\}$ are given by $b_1\le 0, b_3>0$,
and families $\{9,1,6\}$ and $\{8,11,13\}$ are given by
$b_1\le 0, b_3\le 0$ and $b_1+b_2>0$. Finally, $b$ are all such
that $b_2>0,b_2+b_3>0$, where the sector 
$\{b \mid b_1\le 0,b_3\le0,b_1+b_2\le 0\}$
is excluded. 
\vskip 0.2cm

Via Lemma \ref{lem:redmin}, 
the sequences $\{-[\be_i],i\!=\!1,2,3\}$ in this case
are
{\small
\begin{align*}
&\{\al_1,\al_3,\al_2\}, \{\al_1,\al_2+\al_3,\al_2\}, 
\{\al_1+\al_2,\al_3,\al_2\}, \{\al_1+\al_2,\al_2+\al_3,\al_2\},\\ 
&\{\al_1+\al_2+\al_3,\al_3,\al_2\}, 
\{\al_1+\al_2+\al_3,\al_2+\al_3,\al_2\}.
\end{align*}
} 
The inequalities 
$(b,-[\be_i])=-M^i>0$ for $i=1,2,3$ hold but they can give $b$
that do not actually occur.
For instance, family $\{9,1,6\}$ 
with $\bar{\al}_1=-n_2-k,\bar{\al}_3=-n_3,
\bar{\al_2}=1+n_1+n_2+n_3+3k$ results in $\{[-\be_i]\}=$
$\{\al_1+\al_2+\al_3,\al_2+\al_3,\al_2\}$. However, the
same $\{[-\be_i]\}$ are for $\{9,1,8\}$ with 
 $\bar{\al}_1=-n_3-k,\bar{\al}_3=n_1-n_2,
\bar{\al_2}=1+n_1+n_2+2k$. The numerator of $\mu$ 
vanishes at $\bar{\al}_3=n_1-n_2$, which results 
in the cancelation.
\vskip 0.2cm

{\bf Relation to Jackson integrals.}
The following lemma verifies explicitly  that $\Si^{\imath\aa}(f)$
for the standard order of integrations
is proportional to the Jackson integral $J(f;-k\rho)=
\sum_{b\in P}f(\pi_b)
\mu(\pi_b)/\mu(0)$ for functions $f$ invariant
with respect to the (affine) action of $\Pi=\{\pi_i=\pi_{\om_i}=
\pi_1^i,\,
0\le i\le n\}$. The latter condition is not really a restriction, since 
one can replace $f\mapsto \sum_{i=0}^n \pi_i(f)$
because $\pi_i(\mu)=\mu$
and we integrate over $\rho/h+\imath \R^n$, which is $\Pi$-invariant.
However, the $\Pi$-symmetrization of
$f$ can generally worsen the convergence of $\Si^{\imath\aa}(f)$,
as well as $f\mapsto f\!+\!f^\varsigma$ for  $\varsigma\!=\!-w_0$.
The Jackson integral, when it exists, is $\Diamond$-invariant.
 
\begin{lemma} \label{lem:Punion}
For the set $X\subset P$ in the theorem, $P$ is a disjoint union
of $\pi_i(X)$ for $0\le i\le n$, where $\pi_{\om_0}=$id.
I.e. $X$ is a fundamental domain for the action
of $\Pi$ in $P$ 
for the affine action
$bw(z)=b+w(z)$ in 
$\R^n\ni z$; we need 
$\pi_{b}(c) =u_b^{-1}(c)+b$ for $b=\pi_i$ and $c\in P.$
\end{lemma}
{\it Proof.} Let $R^{(n)}$ be $R$ for $A_n$, $S_{n+1}$ 
the corresponding Weyl group.
Recall that $\pi_b=bu_b^{-1}$ has the smallest length
in $\{bw,\,w\in W\}$ (it is unique such); equivalently, 
$u_b\in W$ is of minimal possible length
such that $u_b(b)\in -P_+.$ Then 
$v_1=u_{\om_1}^{-1}$ equals
$s_1\cdots s_n$ (the Coxeter element). It sends
 $\vep_{1}\mapsto\vep_2,\cdots, \vep_{n+1}\mapsto \vep_1$,
$$
v_1:\,\al_1\mapsto \al_2\mapsto \cdots\al_n\mapsto -\th=\al_1+\cdots+
\al_n\mapsto \al_1 \text{ and }$$
$$v_1:\,
\om_1\mapsto \om_2-\om_1\mapsto \cdots \om_n-\om_{n-1}
\mapsto -\om_n\mapsto \om_1.$$

Let us check that:
 
$(i)$\, $R^{(n)}=\medcup_{\,\,m=0}^{\,\,n}\, B_m$ for  
$B_m\equal v_1^m( R^{(n)}_+\setminus R^{(n-1)}_+)$, 
where the union is disjoint and $R^{(n-1)}$ is the root system 
for $\al_1,\ldots, \al_{n-1}$, and

$(ii)$\, given
$b\in P$ such that $(b,\al)\neq 0$ for any $\al\in R^{(n)}$,
there exists a unique $v_1^m$ such that 
$v_1^m (b)$ has positive inner products with all 
$\be\in B_0$.

The set $R^{(n)}_+\setminus R^{(n-1)}_+=
\{\vep_i-\vep_{n+1}, 1\le i\le n\}$
is invariant with respect to $S_{n}$ (for $R^{(n-1)}$)
and only for such $w$. Then,
$B_m= \{\vep_j\!-\!\vep_m, j\neq m\}$ and  their union is the
whole $R^{(n)}$. Explicitly:  
$v_1^m(\vep_i-\vep_{n+1})=\vep_{i\!+\!m \!\!\!\mod\!\! 
(n\!+\!1)+1}-\vep_{m}$ 
for $1\le m \le n$.  It 
contains $\pm \vep_{n+1}$ only for $i\!+\!m=n$. 
However,
the root $\vep_{i+m+1}\!-\!\vep_{m}=\vep_{n+1}\!-\!\vep_m$ 
 is negative for this $i$. Thus, $v_1^m(B_0)\cap B_0=\emptyset$
for any $1\le m\le n$, which proves $(i)$.
To go from $(i)$ to  $(ii)$,  
$b=\sum_{i}c_i\vep_i \in P$ belongs to $B_m$ if and
only if $c_m=\min\{c_i\}$.

Let us switch from $v_1$ to $\pi_1$.  Due to the
inequalities $(b,\vep_j-\vep_{n+1})>0$
for $j\le n$, we have: $X=\{b=\sum_{i}c_i\vep_i 
\in P \mid c_{n\!+\!1}<c_i
\for i\neq n+1\}$, where $c_i\in \Z$ and $\sum_i c_i=0$. Using 
that $\pi_1^m= \om_m v_1^m$ for $1\le m\le n$:
\begin{align*}
\pi_1^m(X)&= \{b=\sum_{i}c_i\vep_{i+m \!\!\!\!\mod\!\! (n+1)+1}+
(\vep_1+\cdots+\vep_m)-\frac{m}{n+1}\sum_{i=1}^{n\!+\!1}\vep_i
\}\\
&=\{b=\sum_{i}b_i\vep_{i} \mid b_m<  b_i \text{ for } i<m \,\text{ and }
\, b_m\le b_{i}\text{ for } i>m\}.
\end{align*}
Here $b_i\in \Z$ and $\sum_{i=1}^n b_i=0$. 
We see that any $b\in P$ can be represented
in this form for a unique $m$: it is such that  
$b_m=\min\{b_i\}$ for 
the smallest index $i$ when this minimum is reached. \sq
\vskip 0.2cm

\section{\sc Integral formulas for 
\texorpdfstring{{\mathversion{normal}$A_2$}}{A2}}\label{sec:a2}
As we see, explicit formulas for $\Si^{\imath\aa}$ can be
obtained in relatively simple way for $A_n$ for the standard order of
$\al_i$ (and the corresponding iterated integration). 
However, 
the problem of finding explicit {\sf finite} sums from
$(i)$ in Theorem \ref{thm:iterint} is subtle even in
this case for
arbitrary $\Re k<0$. We will solve it only 
for $A_2$ and provide the answer for $A_1$. 

For $A_2$, we denote  $x=x_{\al_1}, y=x_{\al_2}, X=q^x, Y=q^y$. 
As above,
the residues are obtained by deleting the binomials of $\mu$
vanishing at $\pi_b$ and evaluating the rest at $t^{-\rho}$.
More generally, the notation $\mu_\bullet$ will be used for this
procedure at any $\xi$ when the numerator of $\mu$ has no zeros.
Due to the $\Pi$-invariance of $\mu$ and the symmetry
$\mu(q^x,q^y)=\mu(q^y,q^x)$:
$$
\mu_\bullet(tq^m,tq^n)\!=\!\mu_\bullet(tq^n,tq^m)\!=\!
\mu_\bullet(tq^n,t^{-2}q^{1-m-n})\!=\!
\mu_\bullet(t^{-2}q^{1-m-n},tq^n).
$$
We will set:\ 
 $\varpi_1(m,n)\!\equal\!\mu_\bullet(tq^m,tq^n),\ 
\varpi_2(m,n)\!\equal\!\mu_\bullet(t^{-1}q^{m\!-\!n},t^2q^n).$ 
They are connected as follows:
\begin{align*}
&\varpi_2(m,m+n)=t^{-1}\frac{1-t^2q^n}{1-q^n}
\varpi_1(m,n) \hbox{\,\, for\,\,} n>0, \hbox{\,\, and\,}\\
&\varpi_2(m,m)=t^{-1}\varpi_1(m,0) \hbox{\,\, due to\,\,}
s_1(\mu)=t^{-1}\frac{1-tX}{1-t^{-1}X}\mu.
\end{align*}
The following are explicit formulas for $\varpi_{1,2}$\,:
\begin{align*}
&\varpi_1(m\!>\!0,n\!>\!0)=t^{3\!-\!2(m+n)}
\prod_{j=1}^{m-1}\frac{(1\!-\!t^2q^j)}{(1\!-\!q^j)}
\prod_{j=1}^{n-1}\frac{(1\!-\!t^2q^j)}{(1\!-\!q^j)}
\prod_{j=1}^{m+n-1}\frac{(1\!-\!t^3q^j)}{(1\!-\!tq^j)},\\
&\varpi_2(m\!>\!0,n\!\ge\! m)\!=\!t^{2\!-\!2n}
\prod_{j=1}^{m-1}\!\frac{(1\!-\!t^2q^j)}{(1\!-\!q^j)}
\frac{(1\!-\!q^{n-j})}{(1\!-\!t^2q^{n-j})}
\prod_{j=1}^{n-1}\!\frac{(1\!-\!t^2q^j)}{(1\!-\!q^j)}
\frac{(1\!-\!t^3q^j)}{(1\!-\!tq^j)}.
\end{align*}

We set
\begin{align*}
&\varrho_0=  
\prod_{i=0}^\infty\frac{(1-t^{-1}q^i)(1-tq^{i+1})}
{(1-q^{i+1})(1-t^2q^{i+1})},\ \ 
\varrho=\mu_\bullet(X=t^{-1},Y=t^{-1})=  \\
&=\prod_{i=0}^\infty\frac{(1\!-\!t^{-1}q^i)(1\!-\!tq^{i+1})^2
(1\!-\!t^{-2}q^i)}
{(1\!-\!q^{i+1})^2(1\!-\!t^2q^{i+1})(1\!-\!t^{3}q^{i+1})}=
\varrho_0\prod_{i=0}^\infty\frac{(1\!-\!tq^{i+1})(1\!-\!t^{-2}q^i)}
{(1\!-\!q^{i+1})(1\!-\!t^{3}q^{i+1})}.
\end{align*}

As above, $q=\exp(-1/\aa), t=q^k=\exp(-k/\aa)$ for $\aa>0$.
 Let
\begin{align*}
Int_{\upsilon}(f)&=\frac{I^{\imath\aa}_\upsilon(f)}
{(2\pi \imath \aa)^2}=
\frac{1}{(2\pi \imath \aa)^2}
\int_{\upsilon-\imath \pi \aa}^{\upsilon+\imath \pi \aa}
\int_{\upsilon-\imath \pi \aa}^{\upsilon+\imath \pi \aa}
f(q^x,q^y)\mu(q^x,q^y) dx dy.
\end{align*}

\begin{proposition}\label{prop:a2si}
(i) Let $\Si_1(f)\!\!=\!\varrho\sum_{m=1,n=1}^{\infty}\,\varpi_1(m,n)
f(tq^m,tq^n)$ and 
$\,\Si_2(f)\!=\!\varrho\sum_{m=1,n=m}^{\infty}\varpi_2(m,n)
f(t^{-1}q^{m-n},t^2q^n).$
Then  $\Si^{\imath\aa}(f)=\Si_1(f)+\Si_2(f)$ provided 
the convergence of\, $\Si_{1,2}(f)$. 

(ii) Let $\Re k <-m/2$ for $m\in \Z_+$. Then $\Si_{1,2}(f)$
converge absolutely for $f=X_a$ with 
$a\in \sum_{i=1}^m \al[i]+\Z_+\al_2+\Z_+(\al_1+\al_2)$,
where $\al[i]$ is either $\pm\al_1$ or $\pm\al_2$. Moreover,
such $X_a$ can be divided by any number of binomials
$(1-X_b)$ for $0\neq 
b\in \Z_+\al_2+\Z_+(\al_1+\al_2)$.

(iii) In particular, ct$(f\mu)$, which is a meromorphic function
of $t=q^k$ for any given Laurent polynomial $f$, coincides with 
$\Si_1(f)+\Si_2(f)$ for any $\Re k<0$
assuming the convergence of $\Si_{1,2}(f)$. This sum coincides
with $Int_{1/3}(f)$ for $-1/3<\Re k<0$; here $Int_{1/3}(f)$
extends
analytically $Int_{0}(f)$ from $\Re k>0$ to  $\Re k> -1/3$ assuming
the integrability. \sq
\end{proposition}

\vskip 0.1cm
Let us provide the corresponding
integral formulas for ct$(f\mu)$ for $A_2$. They are based
on the pole decomposition for $Int_{0}(f)$ from $\Re k>0$ to
$\Re k\le 0$
combined with the formulas from Proposition  \ref{prop:a2si},(i).
We arrange the corresponding (infinite) sum for
$Int_{1/3}(f)-Int_{0}(f)$ as a sum of one-dimensional
integrals and the sum of the remaining residues of
$f\mu$. The integrands for the one-dimensional integrals will be

$
\ze^1_m(q^z)=\mu_\bullet(t^{-1}q^{-m},q^z)=
\mu_\bullet(q^z,t^{-1}q^{-m})=\mu_\bullet(tq^{-z+m+1},q^z)$,

$
\ze^2_m(q^z)=\mu_{\bullet} (t q^{m+1},q^z)=
\mu_\bullet(q^z,tq^{m+1})=
\mu_\bullet(t^{-1}q^{-z-m},q^z).
$
\vskip 0.1cm

One has: $\ze_{m}^2 (q^z)
= \frac{(1-t q^{-z})}{t(1-t^{-1}q^{-z})} \ze^1_m(-z)$. Explicitly:
{\small
\begin{align*}
    \ze^1_m (q^z) &= \varrho_0 t^{-2m} \prod_{j=1}^{m} 
\frac{(1-t^2 q^{j})(1-t^2 q^{j-z})}{(1- q^{j})(1-q^{j-z})} 
\prod_{j=0}^{\infty} \frac{\left(1-q^{-z+j+1}\right)
    \left(1-t^{-1}q^{z+j}\right)}
{\left(1-tq^{z+j}\right)\left(1-t^2 q^{-z+j+1}\right)}, \\
    \ze_{m}^2 (q^z)&= \varrho_0 t^{-2m-1} \prod_{j=1}^{m} 
\frac{(1-t^2 q^{j})(1-t^2 q^{j+z})}{(1- q^{j})(1-q^{j+z})}
\prod_{j=1}^{\infty}\frac{ \left(1-q^{z+j}\right) 
\left(1-t^{-1}q^{-z+j}\right)}{\left(1-tq^{-z+j}\right)
\left(1-t^2 q^{z+j}\right)}. 
\end{align*}
}

\begin{proposition}\label{prop:a2}
For $0\!>\!\Re k\!>\!-0.5$ and any $f\in \C[X^{\pm 1},Y^{\pm 1}]$,
 $\text{ct}(f\mu)=\frac{I_0^{\imath\aa}(f)}
{(2\pi \imath \aa)^2}$
$+\frac{1}{2 \pi \imath \aa} 
\int_{-\imath \pi \aa}^{\imath \pi \aa} 
\left(f(t^{-1},q^y)\!+\!f(q^y,t^{-1})\right)  \ze_0^1(q^y)\, dy$\\
$+\frac{1}{2 \pi \imath \aa} \int_{-\imath \pi \aa}^{\imath \pi \aa} 
f(t^{-1}q^{-y},q^y) \ze_{0}^2 (q^y)\, dy + \varrho f(t^{-1},t^{-1}).
$

For $-0.5\!>\!\Re k\!>\!-1$, the term
$\varrho f(t^{-1},t^{-1})$ here must be
replaced by 
{\small
$$\varrho\Bigr(f(t^{-1},t^{-1})\!+\!\varpi_1(1,1)f(tq,t^{-2}q^{-1})
\!+\!
\varpi_2(1,1)\bigl(f(t^2q,t^{-1})\!+\!f(t^{-1},t^{2}q)\bigr)\Bigr).$$
}

Also, 
the term $\varpi_1(1,1)f(tq,t^{-2}q^{-1})$ in the latter 
sum must be omitted 
when  $\Re k=-0.5$. The functions $f(q^x,q^y)$ here are arbitrary
analytic provided the convergence of the integrals. \sq
\end{proposition}

We note that the formula for $\Re k <-0.5$ contains
the integrand $
f(t^{-1}q^{-y},q^y) \ze_{0}^2 (q^y)$
and the term $\varpi_2(1,1)\bigl(f(t^2q,t^{-1})$
for $\Re k<-0.5$  that are 
not invariant with respect to the symmetry
$\varsigma: x\leftrightarrow y$. The meromorphic continuation 
must by $\varsigma$-invariance, i.e.
the same for $f$ and $f^\varsigma$. Some
symmetries of  $\ze_{0}^2$ and the corresponding cancelation
of residues in this range of $k$ ensure this.

\begin{figure}[htbp]
\centering
\includegraphics[scale=0.3]{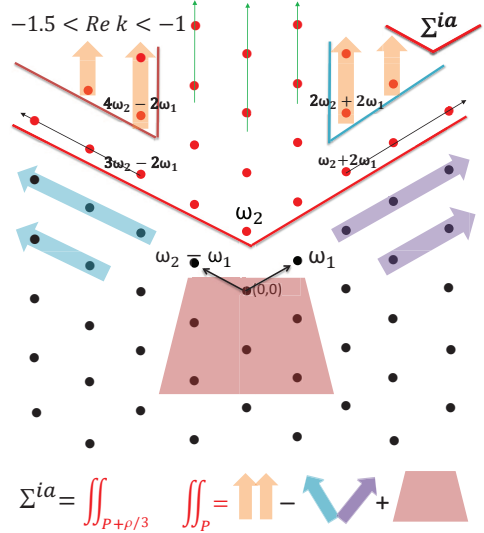}
\vskip -0.5in
\caption{Support of $\Si^{\imath\aa}$ and 
$Int_0$ for $-1.5<\Re k<-1$.} 
\label{ct-a2-1}
\end{figure}
\noindent

\begin{figure}[htbp]
\centering
\includegraphics[scale=0.4]{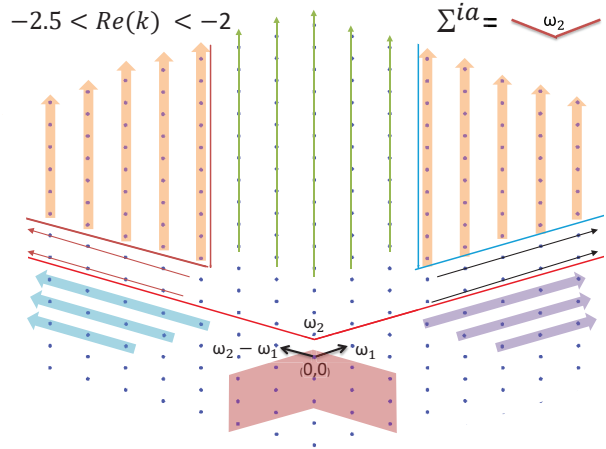}
\vskip -2in
\caption{Support of $\Si^{\imath\aa}$
 and $Int_0$ for $-2.5<\Re k<-2$.} 
\label{ct-a2-2}
\end{figure}

The figures below give  the points $(b_1,b_2)$ 
$b=b_1\om_1+b_2\om_2$
such that $b-ku_b^{-1}(\rho)$ are the corresponding $\mu$-residual
points 
that occur in $\Si^{\imath\aa}(f)$. 
These vectors $b$ form the upward 
sector (angle) with its vertex at $\om_2$.  
It is clearly 
$1/3${\tiny rd} of the total lattice $P$ under the
rotations by $2m\pi/3$ for $m=0,1,2$
with the center at $\rho/3$. 

Recall that this sum (when it converges) 
is proportional
to\, $(a)$ ct$(f\mu)$,\, $(b)$ $\lan f \mu \ran$ for the coinvariant
$\lan \cdot \ran_\xi$ for the anti-involution $\Diamond$
and $\xi=-k\rho,$ \, $(c)$ to the Jackson integral $J_{-k\rho}(f)$
and\,  $(d)$ $\hat{\mathscr{P}}_+(f)$ for the affine symmetrizer
 $\hat{\mathscr{P}}_+.$

These figures show  the set of 
$b$ that occur in  the pole
decomposition of $\frac{I^{ia}_0(f)}
{(2\pi \imath a)^2}$ for the corresponding $\Re k<0$.
They are those belonging to  
the {\sf thick arrows}  and inside the 
polygon containing $(0,0)$. This set is obtained from
the sector describing $\Si^{\imath\aa}$
by removing finitely many
lines and points  and 
adding some lines and points below this sector.

The directions of the lines
that are removed or added give the corresponding integrals
over one-dimensional $\mu$-residual subtori. They are
shown by {\sf thin arrows}, but the residual points due to the 
corresponding one-dimensional 
integrals are not exactly those belonging to these arrows:
some must be added to the corresponding arrows.
 
We note that the presentation of the residual
points of $\Si^{\imath\aa}$ as  $I_0^{\imath\aa}$ plus
those in one-dimensional integrals and 
the remaining points 
is unique (in this picture). We assume that the integrands are
``standard": the partial residues of $\mu$ upon the
restriction to residual tori. Then 
the exact ``thin arrows" are canonically determined by their 
directions.
 
The integral
formulas from Proposition \ref{prop:a2} and their generalizations
to any $\Re k<0$ result combinatorially
from the description of the ``support"
of $I_0^{\imath\aa}$. Recall that the pole expansions
of this integral can be calculated for any
$\Re k\not\in -\Z_+$ but the corresponding analytic functions
will be not connected with each other in different strips. 

The $\mu$-residual points (residual subtori of $dim=0$)
from the integral formula
$\Si^{\imath\aa}(f)=\Int_0(f)+\ldots$
are expected to correspond to square integrable modules
that occur in the regular
DAHA representation, but this is a subject of
some future theory.

Let us provide the integral formulas for $-\ell-0.5<\Re k<-\ell$. 
We omit those for  $-\ell-1<\Re k\le -\ell-0.5$.
Figures \ref{ct-a2-1} and \ref{ct-a2-2} 
are for  $\ell=1,2.$ 

\begin{theorem}
Let $P=[-\imath \pi \aa, \imath \pi \aa]$. For $-\ell\!-\!0.5\!<\!
 \Re k \!<\! 
-\ell$, $\ell\!\in\! \Z_+$:
{\footnotesize
\begin{align*}
&\text{ct}(f\mu)\!=\!
Int_0(f) + \frac{1}{2 \pi \imath \aa}\Bigl(\sum_{m=0}^\ell 
\int_P \left( f(t^{-1}q^{-m},q^y)+
f(q^y,t^{-1}q^{-m})\right) \ze^1_m(q^y) \, dy \\
&+\!\sum_{m=0}^\ell 
\int_P f(t^{-1}q^{-y-m},q^y) \ze^2_m(q^y) \, dy 
\!+\!
\sum_{m=1}^\ell \!\int_P f(t q^{-y+m},q^y) \ze^1_{m-1} (q^y) \, dy \\
&+\!\sum_{m=1}^\ell  \int_P 
\left( f(t q^{m},q^y)\!+\!
f(q^y,t q^{m})\right) \ze^2_{m-1} (q^y) \, dy \Bigr)
\!+\! \varrho\Bigl(\sum_{m,n=1}^\ell \varpi_1(m,n) f(t q^m, t q^n)+ \\
&\sum_{m=1}^{2 \ell} \sum_{n=m\!+\!1}^{2\ell+1} \varpi_1(m,n\!-\!m) 
f(tq^{m},t^{-2}q^{-n\!+\!1})
\!+\! \sum_{m=1}^{2\ell} \sum_{n=m}^{2\ell} \varpi_2(m,n)
f(t^{-1}q^{m\!-\!n},t^2 q^n)+\\ 
&\sum_{m=0}^\ell \sum_{n=m\!+\!1}^{2\ell} 
\!\!\!\varpi_2(n\!-\!m,n) 
f(t^2 q^{n},t^{-1}q^{-m})\!+\!
\sum_{m=1}^{\ell\!+\!1} \sum_{n=m}^{m\!+\!\ell} \!\varpi_2(m,n)
f(t^{-1}q^{m\!-\!n},t^{-1} q^{-m\!+\!1}) \\
&+\!\sum_{m=1}^\ell \sum_{n=m+1}^{2\ell+1} \varpi_1(n-m,m)
f(t^{-2} q^{-n+1},t q^m)\Bigr).
\end{align*}
}
Vectors $b=b_1\om_1+b_2\om_2=\xi^\bullet$ associated with the 
terms in the
double sums can be seen from the corresponding values of $f$, which
are $f(t^{\cdots} q^{b_1},t^{\cdots}q^{b_2})$. For instance,
only the vector with  $m\!=\!1\!=\!n$  from 
$\sum_{m=1}^{\ell\!+\!1} \sum_{n=m}^{m\!+\!\ell}$ occurs
for $\ell=0$; its contribution is   
$\varrho f(t^{-1},t^{-1})$.  \sq
\end{theorem}

Notice that all terms in the integral formula have the
coefficient $1$ in this presentation. We expect this to hold
for $A_n$ and the standard order of $\al_i$, but the evidence
is limited beyond $A_2$. 

\comment{
The points in the upper sector with the vertex at  $\om_2$  
are $b=n_1\om_1+n_2\om_2$ from the support of $\Si^{\imath\aa}$. It is
$1/3${\tiny rd} of the total lattice under the
rotation at $\rho/3$. The points in the 
{\sf thick} arrows and the 
polygon are from the pole expansion of $I^{ia}_0$ for small $\Re k<0$.
The presentation  of supp$(\Si^{ia})$ as supp$(I^{ia}_0)$ plus the
union of those from 
the 1-dim integrals and the residual points is a
combinatorial challenge. The arrows show the 1-dim integrals,
but not exactly the corresponding (full) supports.
The remaining residual points must be 
square integrable representations in some future theory.
}

{\bf The case of \mathversion{normal}{$A_1$}.} For the sake of
completeness, let us provide the integral formula 
from \cite{ChA} in the case of $A_1$. 
As above, $q=e^{-1/\aa}, t=q^k$\  and we set
$x=x_{\al_1}$. 

\begin{proposition} \label{prop:a1}
Let 
$\{j^+,j^-\}\equal 
\{j-1,j\}$ and $\ell\ge 0$ be the
integral part of $-\Re k>0$.
Then for $\mu$ for $A_1$ and $f(q^x)\in\C[q^{\pm x}]$:
\begin{align*}
\text{ct}(f\mu)
=&\,\frac{1}{2\pi \imath \aa }
\int_{-\pi \imath\aa}^{+\pi \imath \aa}
f(q^x) \mu(q^x)\, dx\\
+\,\mu_\bullet(q^{-k})&\,
\Bigl(f(q^{-k})+ \sum_{\ep=\pm}\sum_{j=1}^\ell 
f(q^{\ep(k+j)})
\ t^{-j^\ep}\prod_{i=1}^{j^{\ep}}
\frac{1\!-\!t^2q^i}{1\!-\!q^i}
\Bigr).
\end{align*}
Also, 
$
\text{ct}(f\mu)=
\mu_\bullet(q^{-k})\,
\Bigl(\sum_{j=1}^{\infty} 
f(q^{k+j})
\ t^{1-j}\prod_{i=1}^{j-1}
\frac{1\!-\!t^2q^i}{1\!-\!q^i}
\Bigr),$
where $f(q^x)\!\in 
q^{-\ell x}\C[q^{+x}]$, which provides  the
convergence of this sum.\sq
\end{proposition}

\vskip -3cm
\bibliographystyle{unsrt}


\end{document}

%% file: macro.tex
\renewcommand{\tilde}{\widetilde}
\renewcommand{\hat}{\widehat}

\newcommand{\Z}{{\mathbb Z}}
\newcommand{\Q}{{\mathbb Q}}
\newcommand{\N}{{\mathbb N}}
\newcommand{\C}{{\mathbb C}}
\newcommand{\R}{{\mathbb R}}

\def\HH{\mbox{${\mathcal H}$\kern-5.2pt${\mathcal H}$}}

\newtheorem{theorem}{Theorem}[section]

\newtheorem{proposition}[theorem]{Proposition}
\newtheorem{definition}[theorem]{Definition}
\newtheorem{lemma}[theorem]{Lemma}

\newtheorem{theorem }{Theorem}[section]
\newtheorem{maintheorem }[theorem]{Main Theorem}
\newtheorem{proposition }[theorem]{Proposition}
\newtheorem{definition }[theorem]{Definition}
\newtheorem{lemma }[theorem]{Lemma}
\newtheorem{corollary }[theorem]{Corollary}
\newtheorem{notation }[theorem]{Notation}
\newtheorem{remark }[theorem]{Remark}
\newtheorem{example }[theorem]{Example}

\newtheorem{ maintheorem }[theorem]{Main Theorem}
\newtheorem{ theorem}{Theorem}[section]
\newtheorem{ proposition}[theorem]{Proposition}
\newtheorem{ definition}[theorem]{Definition}
\newtheorem{ lemma}[theorem]{Lemma}
\newtheorem{ corollary}[theorem]{Corollary}
\newtheorem{ notation}[theorem]{Notation}
\newtheorem{ remark}[theorem]{Remark}
\newtheorem{ example}[theorem]{Example}

\def\for{\  \hbox{ for } \ }
\def\iif{ \ \hbox{ if } \ }

\def\and{\  \hbox{ and } \ }
\def\and{\  \hbox{ and } \ }

\def\equal{\stackrel{\,\mathbf{def}}{= \kern-3pt =}}

\def\la{\lambda}
\def\La{\Lambda}
\def\om{\omega}

\def\Th{\Theta}
\def\th{\theta}
\def\al{\alpha}
\def\be{\beta}

\def\ep{\epsilon}

\def\de{\delta}
\def\De{\Delta}

\def\kapp{\hbox{\bf \ae}}
\def\si{\sigma}
\def\Si{\Sigma}
\def\Ga{\Gamma}
\def\ze{\zeta}

\def\vph{\varphi}

\def\vep{\varepsilon}

\def\vth{{\vartheta}}

\def\tal{\tilde{\alpha}}
\def\tbe{\tilde{\beta}}

\def\tu{\tilde{u}}

\def\tw{\widetilde w}
\def\tW{\widetilde W}

\def\tz{\tilde z}
\def\tb{\tilde b}

\def\tR{\tilde R}

\def\hw{\widehat{w}}
\def\hW{\widehat{W}}
\def\hu{\hat{u}}

\def\S{\mathbf{S}}

\def\0{\mathbf{0}}
\def\H{\mathbf{H}}

\def\f{\mathcal{F}}
\def\çF{\mathcal{F}}

\def\r{\mathcal{R}}

\def\p{\mathcal{P}}

\def\h{\mathcal{H}}

\def\y{\mathcal{Y}}
\def\e{\mathcal{E}}

\def\x{\mathcal{X}}

\def\u{\mathcal{U}}

\def\i{\mathcal{I}}
\def\j{\mathcal{J}}

\def\lan{\langle}

\def\ran{\rangle}

\def\lng{\hbox{\rm{\tiny lng}}}
\def\sht{\hbox{\rm{\tiny sht}}}

\newcommand{\Int}{\operatorname{Int}}

\newcommand{\sq}{\phantom{1}\hfill$\qed$}

\def\HH{\mathfrak{H}}

\def\HH{\hbox{${\mathcal H}$\kern-5.2pt${\mathcal H}$}}
\def\HHH{\hbox{${\mathbb H}$\kern-4.2pt${\mathbb H}$}}

\font\smm=msbm10 at 12pt 
\def\symbol#1{\hbox{\smm #1}}
\def\lsmash{{\symbol n}}

\def\#{\sharp}
